\documentclass[11pt, a4paper]{article}
\usepackage[a4paper, total={6.7in, 9.75in}]{geometry}

\title{Effective Darmon's program for the generalised Fermat equation}
\author{\Large Martin \textsc{Azon}}
\date{}

\usepackage[utf8]{inputenc}
\usepackage{amsmath}
\usepackage{amsthm}
\usepackage{amssymb}
\usepackage{amsfonts}
\usepackage{graphicx}
\usepackage{tikz-cd}
\usepackage{enumerate}
\usepackage{hyperref}
\usepackage{xcolor}
\usepackage{pbox,tkz-graph}
\usetikzlibrary{snakes,positioning}
\usepackage{float}
\usepackage{caption}
\usepackage[colorinlistoftodos]{todonotes}
\usepackage{mathrsfs}
\usepackage{mathtools}
\usepackage{multicol}
\usepackage{comment}
\usepackage{array}
\usepackage{stmaryrd}
\usepackage{xparse}
\usepackage{enumitem}
\usepackage{adjustbox}
\usepackage{graphics}
\usepackage{forest}
\usepackage{multirow}
\usepackage{bbm}
\usepackage{tocloft}

\hypersetup{colorlinks, linkcolor={red!70!black}, citecolor={orange!90!yellow},}

\newtheorem{lemma}{Lemma}[section]
\newtheorem{theorem}[lemma]{Theorem}
\newtheorem{proposition}[lemma]{Proposition}
\newtheorem{corollary}[lemma]{Corollary}

\newtheorem{example}[lemma]{Example}
\newtheorem{remark}[lemma]{Remark}
\newtheorem{hypothesis}{Hypothesis}
\newtheorem*{main_idea*}{Main idea}

\theoremstyle{definition}
\newtheorem{definition}[lemma]{Definition}

\usepackage{pbox}
\usepackage{tkz-graph}
\usetikzlibrary{snakes,fit,positioning,calc,shapes}
\usepackage{graphicx}
\usepackage{cleveref,tikz-cd,pbox,wrapfig}

\definecolor{amethyst}{rgb}{0.6, 0.4, 0.8}
\definecolor{atomictangerine}{rgb}{1.0, 0.6, 0.4}
\definecolor{deeppeach}{rgb}{1.0, 0.8, 0.64}
\definecolor{eggshell}{rgb}{0.94, 0.92, 0.84}
\definecolor{lightapricot}{rgb}{0.99, 0.84, 0.69}
\definecolor{lemonchiffon}{rgb}{1.0, 0.98, 0.8}
\definecolor{roundabout}{rgb}{1.0, 0.91, 0.75}
\definecolor{atomictangerine}{rgb}{1.0, 0.6, 0.4}
\definecolor{ruby}{rgb}{0.88, 0.07, 0.37}
\definecolor{sapphire}{rgb}{0.03, 0.15, 0.4}

\def\rootsep{0.05}               % horizontal space between roots in a cluster picture
\def\clustersep{0.06}            % horizontal and vertical space between parent & child cluster
\def\cnamescale{0.8}             % cluster names font size 
\def\cdepthscale{0.7}            % cluster depths and signs font size 
\def\cltopskip{1pt}              % space above a cluster picture
\def\clbottomskip{1pt}           % space below a cluster picture

\def\rootscale{0.5}   \def\rootcolor{gray}
\def\rootscaleA{0.7}  \def\rootcolorA{yellow}
\def\rootscaleB{0.5}  \def\rootcolorB{green}
\def\rootscaleC{0.4}  \def\rootcolorC{sapphire}
\def\rootscaleD{0.45}  \def\rootcolorD{ruby}
\tikzset{
	clA/.style = {very thick,black},
	clB/.style = {thick,purple}
	clC/.style = {thick, red}
}

\tikzset{
	root/.style = {circle,scale=\rootscale,fill=\rootcolor},
	rc/.style 2 args = {right=#1*1.5*\clustersep of {#2.east|-first},root}, rr/.style = {right=\rootsep of {#1.east|-first},root},
	roott/.style = {circle,inner sep=-2pt,minimum size=5pt,black,font=\ttfamily\footnotesize},
	rct/.style 2 args = {right=#1*1.5*\clustersep of {#2.east|-first},roott}, rrt/.style = {right=\rootsep of {#1.east|-first},roott},
	rootA/.style = {circle,scale=\rootscaleA,ball color=\rootcolorA},
	rcA/.style 2 args = {right=#1*1.5*\clustersep of {#2.east|-first},rootA}, rrA/.style = {right=\rootsep of {#1.east|-first},rootA},
	rootB/.style = {circle,scale=\rootscaleB,ball color=\rootcolorB},
	rcB/.style 2 args = {right=#1*1.5*\clustersep of {#2.east|-first},rootB}, rrB/.style = {right=\rootsep of {#1.east|-first},rootB},
	rootC/.style = {diamond,scale=\rootscaleC,ball color=\rootcolorC},
	rcC/.style 2 args = {right=#1*1.5*\clustersep of {#2.east|-first},rootC}, rrC/.style = {right=\rootsep of {#1.east|-first},rootC},
	rootD/.style = {circle,scale=\rootscaleD,ball color=\rootcolorD},
	rcD/.style 2 args = {right=#1*1.5*\clustersep of {#2.east|-first},rootD}, rrD/.style = {right=\rootsep of {#1.east|-first},rootD},
	rootE/.style = {circle,scale=0.1,fill=black},
	rcE/.style 2 args = {right=#1*1.5*\clustersep of {#2.east|-first},rootE}, rr/.style = {right=\rootsep of {#1.east|-first},rootE},
	cluster/.style = {draw=black!90,thick,rounded corners,inner sep=22*\clustersep,outer xsep=22*\clustersep,fit=#1},
	clusterc/.style = {draw=white,thick,rounded corners,inner sep=22*\clustersep,outer xsep=22*\clustersep,fit=#1},
	clabel/.style  = {anchor=west,scale=\cdepthscale,black,inner sep=0,outer xsep=1,outer ysep=0},
	clabelL/.style = {above right=-\clustersep of #1t.north east,clabel},
	clabelD/.style = {below right=-\clustersep of #1t.south east,clabel},
	clouter/.style = {inner sep=0,outer sep=0,fit=#1}
}

\def\Cluster #1 = #2;{\node[cluster=#2] (#1) {};}
\def\ClusterL #1[#2] = #3;{
	\node[cluster=#3] (#1t) {}; \node[clabelL=#1] (#1l) {$#2$}; \node[clouter=(#1t)(#1l)] (#1) {};}
\def\ClusterD #1[#2] = #3;{
	\node[cluster=#3] (#1t) {}; \node[clabelD=#1] (#1d) {$#2$}; \node[clouter=(#1t)(#1d)] (#1) {};}
\def\ClusterLD #1[#2][#3] = #4;{
	\node[cluster=#4] (#1t) {}; \node[clabelL=#1] (#1l) {$#2$}; 
	\node[clabelD=#1] (#1d) {$#3$}; \node[clouter=(#1t)(#1l)(#1d)] (#1) {};}
\def\ClusterLDName #1[#2][#3][#4] = #5;{
	\node[cluster=#5] (#1t) {}; \node[clabelL=#1] (#1l) {$#2$}; 
	\node[clabelD=#1] (#1d) {$#3$}; 
	\node[scale=\cnamescale,above=\clustersep/3 of #1t,inner sep=0, outer sep=0] (#1n) {$#4$}; 
	\node[clouter=(#1l)(#1d)(#1t)] (#1) {};}
\def\ClustercLDName #1[#2][#3][#4] = #5;{
	\node[clusterc=#5] (#1t) {}; \node[clabelL=#1] (#1l) {$#2$}; 
	\node[clabelD=#1] (#1d) {$#3$}; 
	\node[scale=\cnamescale,above=\clustersep/3 of #1t,inner sep=0, outer sep=0] (#1n) {$#4$}; 
	\node[clouter=(#1l)(#1d)(#1t)] (#1) {};}

\newcommand{\Root}[4][]{
	\ifx\relax#2\relax\node[rr#1=#3] (#4) {};\else\node[rc#1={#2}{#3}] (#4) {};\fi}
\newcommand{\RootT}[5][]{
	\ifx\relax#2\relax\node[rrt#1=#3] (#4) {#5};\else\node[rct#1={#2}{#3}] (#4) {#5};\fi}

\def\frob(#1)(#2){\path[draw,thick,shorten <=-22*\clustersep,shorten >=-22*\clustersep](#1.east)--(#2.west|-#1){};}

\usepackage{pbox}
\def\pb#1{\pbox[c]{\textwidth}{\hfil #1\hfil}}

\long\def\clusterpicture#1\endclusterpicture{\pb{\vbox to \cltopskip{\vfill}\\%
		\begin{tikzpicture}\node[coordinate] (first) {};#1\end{tikzpicture}\\[-11pt]\vbox to \clbottomskip{\vfill}}}   
\long\def\clusterpictureopt#1#2\endclusterpicture{\pb{\vbox to \cltopskip{\vfill}\\%
		\begin{tikzpicture}[#1]\node[coordinate] (first) {};#2\end{tikzpicture}\\[-11pt]\vbox to \clbottomskip{\vfill}}}

\def\pb#1{\pbox[c]{\textwidth}{\hfil #1\hfil}}

\newcommand{\cf}{\textit{cf.}}
\newcommand{\ie}{\textit{i.e.,\ }}

\newcommand{\ord}{\mathrm{ord}}

\newcommand{\st}{^{\times}}

\renewcommand{\mod}{\operatorname{\, mod}}
\newcommand{\lcm}{\operatorname{lcm}}
\newcommand{\charact}{\operatorname{char}}

\newcommand{\abs}[1]{\left\arrowvert #1 \right\arrowvert }

\newcommand{\sd}{\operatorname{sd}}

\newcommand{\C}{\mathbb{C}}
\newcommand{\Z}{\mathbb{Z}}
\newcommand{\Q}{\mathbb{Q}}

\newcommand{\Ff}{\mathbb{F}}
\newcommand{\Qq}{\Q_q}
\newcommand{\Qp}{\Q_p}
\newcommand{\Qr}{\Q_r}
\newcommand{\Zq}{\Z_q}
\newcommand{\Zr}{\Z_r}
\newcommand{\Qbar}{\overline{\Q}}
\newcommand{\Zbar}{\overline{\Z}}
\newcommand{\Qlbar}{\overline{\Q_{\ell}}}
\newcommand{\Qqbar}{\overline{\Q_{q}}}
\newcommand{\Qrbar}{\overline{\Q_{r}}}
\newcommand{\Fpbar}{\overline{\Fp}}
\newcommand{\Fqbar}{\overline{\Fq}}
\newcommand{\Kqbar}{\overline{\Kq}}

\newcommand{\Kglbar}{\overline{\Kgl}}

\newcommand{\Fp}{\Ff_p}

\newcommand{\Tr}[2]{\operatorname{Tr}_{#1 / #2}}
\newcommand{\Nrm}[2]{\operatorname{N}_{#1 / #2}}

\newcommand{\Trace}{\operatorname{Tr}}

\newcommand{\sep}[1]{#1^{\mathrm{sep}}}
\newcommand{\unr}[1]{#1^{\mathrm{unr}}}
\renewcommand{\O}{\mathcal{O}}

\newcommand{\Ql}{\Q_{\ell}}

\newcommand{\Gal}{\operatorname{Gal}}

\newcommand{\Pic}{\operatorname{Pic}}

\newcommand{\Spec}{\operatorname{Spec}}

\newcommand{\Jac}{\operatorname{Jac}}

\newcommand{\radst}{\operatorname{rad}^{*}}
\newcommand{\zr}{\zeta_r}
\newcommand{\zp}{\zeta_p}
\newcommand{\eps}{\varepsilon}

\newcommand{\Cms}{C_{r}^{-}}
\newcommand{\Cpl}{C_{r}^{+}}
\newcommand{\CmssQ}{\Cms(s_0)^{(\deltaQ)}}

\newcommand{\Cmsabc}{C_{r}^{-}(a, b, c)}
\newcommand{\Cplabc}{C_{r}^{+}(a, b, c)}
\newcommand{\Cpmabc}{C_{r}^{\pm}(a, b, c)}
\newcommand{\Jpmabc}{J_{r}^{\pm}(a, b, c)}
\newcommand{\Jmsabc}{J_{r}^{-}(a, b, c)}
\newcommand{\Jplabc}{J_{r}^{+}(a, b, c)}
\newcommand{\Crabc}{C_{r}(a, b, c)}
\newcommand{\Cpm}{C_{r}^{\pm}}
\newcommand{\Cmstw}{(\Cms)^{(\deltaK)}}
\newcommand{\Cpltw}{(\Cpl)^{(\deltaK)}}
\newcommand{\Cpmtw}{(\Cpm)^{(\deltaK)}}
\newcommand{\Jms}{J_{r}^{-}}
\newcommand{\Jpl}{J_{r}^{+}}
\newcommand{\Jpm}{J_{r}^{\pm}}
\newcommand{\Jmstw}{(J_{r}^{-})^{(\deltaK)}}
\newcommand{\Jpltw}{(J_{r}^{+})^{(\deltaK)}}
\newcommand{\Jpmtw}{(J_{r}^{\pm})^{(\deltaK)}}
\newcommand{\Wrst}{\mathcal{W}}
\newcommand{\Wms}{\Wrst_{r}^{-}}
\newcommand{\Wpl}{\Wrst_{r}^{+}}
\newcommand{\Wpm}{\Wrst_{r}^{\pm}}
\newcommand{\Kgl}{\mathcal{K}}
\newcommand{\Qquad}{\mathcal{Q}}
\newcommand{\Qone}{\mathcal{Q}_{1}}
\newcommand{\Qtau}{\mathcal{Q}_{\tau}}
\newcommand{\OK}{\O_{\Kgl}}
\newcommand{\OF}{\O_{F}}

\newcommand{\frakP}{\mathfrak{P}}
\newcommand{\frakQ}{\mathfrak{Q}}
\newcommand{\frakp}{\mathfrak{p}}
\newcommand{\frakq}{\mathfrak{q}}
\newcommand{\frakm}{\mathfrak{m}}

\newcommand{\frakr}{\mathfrak{r}}
\newcommand{\Ffp}{\Ff_{\frakp}}
\newcommand{\Ffq}{\Ff_{\frakq}}
\newcommand{\Ffm}{\Ff_{\frakm}}
\newcommand{\Ffr}{\Ff_{\frakr}}
\newcommand{\Fr}{\Ff_{r}}
\newcommand{\Frbar}{\overline{\Ff_{r}}}
\newcommand{\Kp}{\Kgl_{\frakp}}
\newcommand{\Kq}{\Kgl_{\frakq}}
\newcommand{\Kr}{\Kgl_{\frakr}}
\newcommand{\Cmod}{\mathcal{C}}
\newcommand{\Xmod}{\mathcal{X}}

\newcommand{\Gmod}{\mathcal{G}}

\newcommand{\Leg}{\mathcal{L}}
\newcommand{\Cmodm}{\mathcal{C}_{\frakm}}

\newcommand{\Jmodpm}{\Jmod_{r}^{\pm}}
\newcommand{\Gmodpm}{\Gmod_{r}^{\pm}}
\newcommand{\Gmodms}{\Gmod_{r}^{-}}
\newcommand{\Gmodpl}{\Gmod_{r}^{+}}

\newcommand{\Jmodms}{\Jmod_{r}^{-}}
\newcommand{\Jmodpl}{\Jmod_{r}^{+}}
\newcommand{\Jmodmstw}{(\Jmod_{r}^{-})^{(\deltaK)}}
\newcommand{\Jmodpltw}{(\Jmod_{r}^{+})^{(\deltaK)}}
\newcommand{\Jmodpmtw}{(\Jmod_{r}^{\pm})^{(\deltaK)}}
\newcommand{\Wpmtw}{\Wrst_{r}^{\pm \, (\deltaK)}}
\newcommand{\Wpltw}{\Wrst_{r}^{+ \, (\deltaK)}}
\newcommand{\Wmstw}{\Wrst_{r}^{- \, (\deltaK)}}

\newcommand{\gminus}{g_{r}^{-}}
\newcommand{\gplus}{g_{r}^{+}}
\newcommand{\gminuss}{g_{r, s}^{-}}
\newcommand{\gpluss}{g_{r, s}^{+}}
\newcommand{\gpm}{g_{r}^{\pm}}
\newcommand{\deltatilde}{\widetilde{\delta}}
\newcommand{\stilde}{\widetilde{s}}

\newcommand{\Aut}{\operatorname{Aut}}
\newcommand{\Res}{\operatorname{Res}}
\newcommand{\GL}{\operatorname{GL}}
\newcommand{\SL}{\operatorname{SL}}

\newcommand{\End}{\operatorname{End}}

\newcommand{\TlA}{\operatorname{T}_{\ell}(A)}
\newcommand{\Vl}{\operatorname{V}_{\ell}}
\newcommand{\VlA}{\operatorname{V}_{\ell}(A)}
\newcommand{\Fl}{F_{\ell}}
\newcommand{\Flam}{F_{\lambda}}
\newcommand{\Klam}{\Kgl_{\lambda}}
\newcommand{\Elam}{E_{\lambda}}
\newcommand{\V}{\operatorname{V}}
\newcommand{\T}{\operatorname{T}}
\newcommand{\Vlam}{\operatorname{V}_{\lambda}}
\newcommand{\VlamA}{\operatorname{V}_{\lambda}(A)}

\newcommand{\rholA}{\rho_{A, \, \ell}}

\newcommand{\rholam}{\rho_{\lambda}}
\newcommand{\rholamA}{\rho_{A, \, \lambda}}

\newcommand{\F}{\mathcal{F}}
\newcommand{\Fq}{\mathcal{F}_{\frakq}}
\newcommand{\Oq}{\O_{\frakq}}
\newcommand{\Om}{\O_{\frakm}}
\newcommand{\OP}{\O_{\frakP}}

\newcommand{\disc}{\operatorname{disc}}
\newcommand{\vq}{v_{\frakq}}
\newcommand{\vr}{v_{\frakr}}
\newcommand{\vm}{v_{\frakm}}
\newcommand{\vP}{v_{\frakP}}
\newcommand{\vQ}{v_{\frakQ}^{ }}
\newcommand{\vF}{v_{F}}
\newcommand{\vgam}{v_{\Qr(\gamma_0)}}
\newcommand{\RM}{$\mathsf{RM}$}
\newcommand{\CM}{$\mathsf{CM}$}
\newcommand{\WD}{$\operatorname{WD}$}
\newcommand{\W}{\operatorname{W}}
\newcommand{\Wi}{\operatorname{W}_{i}}
\newcommand{\Jmod}{\mathcal{J}}
\newcommand{\Amod}{\mathcal{A}}
\newcommand{\Amodm}{\mathcal{A}_{\frakm}}

\newcommand{\Rroots}{\mathcal{R}}
\newcommand{\s}{\mathfrak{s}}

\newcommand{\Kqroots}{\Kq(\Rroots)}
\newcommand{\Krroots}{\Kr(\Rroots)}

\newcommand{\Qqroots}{\Qq(\Rroots)}

\newcommand{\sqt}{\sqrt{s_0(s_0-1)}}
\newcommand{\sqts}{\sqrt{s_0}}
\newcommand{\sqtsm}{\sqrt{s_0-1}}
\newcommand{\deltaQ}{\delta_{\Q}}
\newcommand{\deltaK}{\delta_{\Kgl}}

\newcommand{\Hhyp}{\mathcal{H}}
\newcommand{\Eppr}{\mathcal{E}_{p, p, r}}
\newcommand{\Errp}{\mathcal{E}_{r, r, p}}
\newcommand{\Epqr}{\mathcal{E}_{p, q, r}}
\newcommand{\Ppr}{\mathbb{P}^{1}}

\newcommand{\Frob}{\operatorname{Frob}}
\newcommand{\rhobar}{\overline{\rho}}

\newcommand{\rhopm}{\rho_{\Jpm, \, \lambda}}
\newcommand{\rhoms}{\rho_{\Jms, \, \lambda}}
\newcommand{\rhopl}{\rho_{\Jpl, \, \lambda}}
\newcommand{\rhomsK}{\rho_{\Jmstw, \, \lambda}}
\newcommand{\rhoplK}{\rho_{\Jpltw, \, \lambda}}
\newcommand{\rhopmK}{\rho_{\Jpmtw, \, \lambda}}
\newcommand{\Rep}{\operatorname{Rep}}
\newcommand{\phim}{\varphi_{\frakm}}
\newcommand{\Cond}[1]{\N \left( #1 \right) }
\newcommand{\condexp}[1]{\mathfrak{n} (#1)}
\newcommand{\condtame}[1]{\mathfrak{n}_{\operatorname{tame}} \left(#1 \right)}
\newcommand{\condwild}[1]{\mathfrak{n}_{\operatorname{wild}} \left(#1 \right)}

\newcommand{\codim}{\operatorname{codim}}
\newcommand{\B}{\mathcal{B}}
\newcommand{\SQ}[1]{\mathsf{SQ}( #1) }

\newcommand{\rhojtw}{\rho_{\Jpmtw, \, \lambda}}
\newcommand{\rhojlam}{\rho_{\Jpm, \, \lambda}}

\newcommand{\rhopmbar}{\overline{\rho}_{\Jpm, \, \frakr}}
\newcommand{\rhomsbar}{\overline{\rho}_{\Jms, \, \frakr}}
\newcommand{\rhoplbar}{\overline{\rho}_{\Jpl, \, \frakr}}
\newcommand{\rhopmpbar}{\overline{\rho}_{\Jpmtw, \, \frakp}}
\newcommand{\rhomspbar}{\overline{\rho}_{\Jmstw, \, \frakp}}
\newcommand{\rhoplpbar}{\overline{\rho}_{\Jpltw, \, \frakp}}

\newcommand{\cdmsK}{\mathfrak{n}_{\frakq}^{-}}
\newcommand{\cdplK}{\mathfrak{n}_{\frakq}^{+}}
\newcommand{\cdpmK}{\mathfrak{n}_{\frakq}^{\pm}}
\newcommand{\ctamemsK}{\mathfrak{n}_{\operatorname{tame}, \, \frakq}^{-}}
\newcommand{\ctameplK}{\mathfrak{n}_{\operatorname{tame}, \, \frakq}^{+}}
\newcommand{\ctamepmK}{\mathfrak{n}_{\operatorname{tame}, \, \frakq}^{\pm}}
\newcommand{\cwildmsK}{\mathfrak{n}_{\operatorname{wild}, \, \frakq}^{-}}
\newcommand{\cwildplK}{\mathfrak{n}_{\operatorname{wild}, \, \frakq}^{+}}
\newcommand{\cwildpmK}{\mathfrak{n}_{\operatorname{wild}, \, \frakq}^{\pm}}
\newcommand{\cdmsrK}{\mathfrak{n}_{\frakr}^{-}}
\newcommand{\cdplrK}{\mathfrak{n}_{\frakr}^{+}}
\newcommand{\cdpmrK}{\mathfrak{n}_{\frakr}^{\pm}}
\newcommand{\ctamemsrK}{\mathfrak{n}_{\operatorname{tame}, \, \frakr}^{-}}
\newcommand{\ctameplrK}{\mathfrak{n}_{\operatorname{tame}, \, \frakr}^{+}}

\newcommand{\cwildmsrK}{\mathfrak{n}_{\operatorname{wild}, \, \frakr}^{-}}
\newcommand{\cwildplrK}{\mathfrak{n}_{\operatorname{wild}, \, \frakr}^{+}}

\newcommand{\ntwopm}{\mathfrak{n}_{2}^{\pm}}
\newcommand{\ntwoms}{\mathfrak{n}_{2}^{-}}
\newcommand{\ntor}{\mathfrak{n}_{\operatorname{tor}}}
\newcommand{\ntornotp}{\mathfrak{n}_{\operatorname{tor}, \, \not \equiv \, 0 \, (p)}}
\newcommand{\nunip}{\mathfrak{n}_{\operatorname{unip}}}

\newcommand{\rhiotmsqK}{\rho_{\iota, \frakq}^{-}}
\newcommand{\rhiotplqK}{\rho_{\iota, \frakq}^{+}}
\newcommand{\rhiotpmqK}{\rho_{\iota, \frakq}^{\pm}}
\newcommand{\rhiotpmqprK}{\rho_{\iota, \frakq'}^{\pm}}

\newcommand{\rhiotpmrK}{\rho_{\iota, \frakr}^{\pm}}

\newcommand{\N}{\mathcal{N}}
\newcommand{\Snew}[1]{\mathcal{S}_{2}(#1)}
\newcommand{\aqJ}{a_{\frakq}(J)}
\newcommand{\aqg}{a_{\frakq}(g)}

\newcommand{\Magma}{{\fontfamily{lmtt}\selectfont Magma}}
\newcommand{\Ngood}{N_{\operatorname{good}}(\frakq, g)}
\newcommand{\Mtoric}{M_{\operatorname{toric}}(\frakq, g)}

\let\emptyset\varnothing

\newcolumntype{P}[1]{>{\centering\arraybackslash}p{#1}}
\renewcommand{\arraystretch}{1.1}

\newcommand{\boxppr}{\fbox{$\boldsymbol{(p, p, r)}$} {\normalfont{:}} }
\newcommand{\boxrrp}{\fbox{$\boldsymbol{(r, r, p)}$} {\normalfont{:}} }

\makeatletter
\newcommand{\leqnomode}{\tagsleft@true\let\veqno\@@leqno}
\newcommand{\reqnomode}{\tagsleft@false\let\veqno\@@eqno}

\makeatother

\DeclarePairedDelimiterX{\Iintv}[1]{\llbracket}{\rrbracket}{\iintvargs{#1}}
\NewDocumentCommand{\iintvargs}{>{\SplitArgument{1}{,}}m}
{\iintvargsaux#1} 
\NewDocumentCommand{\iintvargsaux}{mm} {#1\mkern1.5mu..\mkern1.5mu#2}

\let\originalleft\left
\let\originalright\right
\renewcommand{\left}{\mathopen{}\mathclose\bgroup\originalleft}
\renewcommand{\right}{\aftergroup\egroup\originalright}

\newcommand\restr[2]{#1\raisebox{-.5ex}{$|$}_{#2}}

\numberwithin{equation}{section}

\DeclareFontFamily{U}{wncy}{}
\DeclareFontShape{U}{wncy}{m}{n}{<->wncyr10}{}
\DeclareSymbolFont{mcy}{U}{wncy}{m}{n}
\DeclareMathSymbol{\Sh}{\mathord}{mcy}{"58}

\begin{document}
	
	\maketitle
	\vspace{-1.5em}
	\noindent\hfill\rule{7cm}{0.5pt}\hfill\phantom{.}
	\begin{center}
		\parbox{5.75in}{\vspace{-1em}\small\paragraph{Abstract --}
			We follow the ideas of Darmon's program for solving infinite families of generalised Fermat equations of signatures $(p,p,r)$ and $(r,r,p)$, where $r$ is a fixed prime and $p$ is varying. We do so by introducing a common framework for both signatures, allowing for a uniform treatment for the two families of equations. We analyse in detail the geometry of Frey hyperelliptic curves, and the reduction types of the Néron models of their Jacobians. We then study the associated $2$-dimensional Galois representations: modularity, irreducibility, and level lowering. We provide a Magma package that performs the elimination step for many choices of coefficients and the exponent $r$. In order to illustrate the effectiveness of our results, we solve several examples of families of equations of signatures $(p,p,5)$ and $(5, 5, p)$.

			\vspace{1em}
			\noindent {\it Keywords:} Generalised Fermat equation, Darmon's program, Frey hyperelliptic curves.

			\smallskip
			\noindent {\it 2020 Math.\ Subj.\ Classification:}  
			11D41	%Higher degree equations; Fermat's equation
			11G30	%Curves of arbitrary genus or genus $\ne 1$ over global fields 
			11F80	%Galois representations
			11Y50	%Computer solution of Diophantine equations	
		}
	\end{center}

	\noindent\hfill\rule{7cm}{0.5pt}\hfill\phantom{.}

\maketitle

{\hypersetup{hidelinks} \tableofcontents}

% -- % -- % -- % -- % -- % -- % -- % -- % -- % -- % -- % -- % -- % -- % -- % -- % -- % -- % -- % -- % -- % -- % -- % -- % -- % -- % -- % -- % -- % -- % -- % -- % -- % -- % -- % -- % -- % -- % -- 

\section[Introduction]{Introduction}

\subsubsection*{The generalised Fermat equation}

Wiles' proof of  Fermat's last theorem \cite{Wiles95} was one of the major breakthroughs in the recent history of number theory. Since then, there has been great progress in the study of modularity of Galois representations, leading to new Diophantine results. Building upon ideas of Frey, Serre, Ribet among others, many mathematicians have explored how to apply techniques from the study of Galois representations to the resolution of Diophantine equations. In particular, of great interest is the \textit{generalised Fermat equation} (GFE)
\leqnomode
\begin{equation*}\label{eq: GFE pqr}\tag{$\Epqr$}
	 A x^p + B y^q = C z^r.
\end{equation*}
Here we denote by $A, B, C$ three non-zero pairwise coprime integers, and by $p, q, r$ three integers greater or equal to $2$. The triple $(p, q, r)$ is called the \textit{signature} of the equation \eqref{eq: GFE pqr}. A solution $(a, b, c) \in \Z^3$ to \eqref{eq: GFE pqr} is called \textit{non-trivial} if $abc \neq 0$, and \textit{primitive} if $\gcd(a, b, c) = 1$.
\medskip

Darmon and Granville proved that, if we fix $A, B, C, p, q, r$ and we further assume that $\frac{1}{p} + \frac{1}{q} + \frac{1}{r} < 1$, then equation \eqref{eq: GFE pqr} has only finitely many non-trivial primitive solutions (see \cite{DarmonGranville}). In the last few decades, several instances of \eqref{eq: GFE pqr} have been solved for fixed numerical values of every considered parameter (we refer the reader to \cite{RatcliffeGrechuk25} for an overview of the current state of the art). On a different direction, various authors have studied infinite families of equations of the shape \eqref{eq: GFE pqr}, where one allows some of the parameters $p, q, r$ in the signature to vary.
\medskip

In this paper, we study the infinite families of generalised Fermat equations \eqref{eq: GFE ppr}$_p$ and \eqref{eq: GFE rrp}$_p$, of respective signatures $(p,p,r)$ and $(r,r,p)$. We fix, once and for all, a prime number $r$ and three non-zero pairwise coprime integers $A, B, C$. The parameter $p$ ranges through an infinite family of prime numbers, \textit{e.g.,} the prime numbers greater than some fixed quantity (in this case we talk about an \textit{asymptotic result}). Even if the roles of $r$ and $p$ might seem to be symmetric, the families of equations \eqref{eq: GFE ppr}$_p$ and \eqref{eq: GFE rrp}$_p$ are distinct. For any fixed value of $p$, the equations above are given by
\begin{align}
	\label{eq: GFE ppr}\tag{$\Eppr$}  &  A x^p + By^p = C z^r, \\
	\label{eq: GFE rrp}\tag{$\Errp$}  &  A x^r + By^r = C z^p.
\end{align}

Our goal in this paper is to develop the necessary theoretical background to solve the families of equations \eqref{eq: GFE ppr}$_p$ and \eqref{eq: GFE rrp}$_p$ for given numerical values of $r, A, B, C$. We develop the ideas of Darmon's program (which will be explained below) to establish effective results regarding the solutions to such equations. We provide a \Magma \ package to perform the elimination step for several choices of exponent $r$ and coefficients $A, B, C$, which therefore allows to solve many families of GFEs. To illustrate the explicitness of our discussion, we specify $r = 5$ and solve different instances of the families of equations $(\mathcal{E}_{p, p, 5})_p$ and $(\mathcal{E}_{5, 5, p})_p$. More precisely, we prove the following asymptotic results.

\begin{theorem}[= Theorems~\ref{thm: Example 1 pp5} + \ref{thm: Example 2 pp5}]\label{thm: Example pp5}
	Let $p > 71$ be a prime number. 
	\begin{enumerate}[itemsep=2pt, topsep=5pt]
		\item There are no primitive non-trivial solutions $(a, b, c) \in \Z^3$ to the generalised Fermat equation 
		\begin{equation*}
			7 x^p +  y^p = 3 z^5 \vspace{-0.5em}
		\end{equation*}
		that satisfy $10 \mid ab$.
		
		\item For any $i \in \Iintv{1,4}$ and $j \in \lbrace 3, 4 \rbrace$, there are no primitive non-trivial solutions to the generalised Fermat equation
		\begin{equation*}
			7 x^p + 2^i 5^j y^p = 3 z^5.
		\end{equation*}
	\end{enumerate}
\end{theorem}

~ \newpage 
\begin{theorem}[= Theorems~\ref{thm: Example 1 55p} + \ref{thm: Example 2 55p}]\label{thm: Example 55p}
	Let $p > 41$ be a prime number. 
	\begin{enumerate}[itemsep=2pt, topsep=5pt]
		\item There are no primitive non-trivial solutions $(a, b, c) \in \Z^3$ to the generalised Fermat equation 
		\begin{equation*}
			x^5 + 7 y^5 = z^p \vspace{-0.5em}
		\end{equation*}
		that satisfy $10 \mid c$. 
		
		\item For any $i \in \Iintv{1,4}$ and $j \in \Iintv{2,4}$, there are no primitive non-trivial solutions to the generalised Fermat equation
		\begin{equation*}
			x^5 + 7 y^5 = 2^i 5^j z^p.
		\end{equation*}
	\end{enumerate}
\end{theorem}

We stress out the fact that the equations above are the first examples of GFEs of signature $(r, r, p)$ and $(p, p, r)$ with $r \geq 5$ and $|AB| \neq 1$ that have been solved (even under $2$-adic and $5$-adic conditions). In section~\ref{sect: Solving instances GFEs}, the reader will find explanations, both theoretical and computational, on why we considered such examples.

\subsubsection*{The modular method}

The approach that we follow here for establishing Theorems~\ref{thm: Example pp5} and \ref{thm: Example 55p} is the so-called \textit{modular method}, which was pioneered by Frey, Serre, Ribet, and Wiles, among others. It allows to solve instances of infinite families of equations with one varying parameter in the signature, and may be summarised as follows:
\begin{enumerate}[leftmargin=1.75cm]
	\item[\textbf{Step 1:}] \textbf{Construct a Galois representation.} Attach to any putative solution a $2$-dimensional representation $\rho : \Gal(\Kglbar / \Kgl) \rightarrow \GL_2(\overline{\Qp})$. Here, $\Kgl$ is a totally real number field, and $p$ is the varying parameter in the signature. We require to explicitly describe the ramification of $\rho$, which should behave ``well" in terms of the putative solution.

	\item[\textbf{Step 2:}] \textbf{Modularity / Level lowering.} Prove modularity of $\rho$, and irreducibility of the residual representation $\rhobar$. Applying level lowering results, show that $\rhobar$ arises from a Hilbert newform over $\Kgl$ of parallel weight $2$ whose level is independent of the solution.

	\item[\textbf{Step 3:}] \textbf{Elimination.} Compute all newforms predicted in the previous step, and prove that none of them gives rise to $\overline{\rho}$.
\end{enumerate}

In the case of Fermat's last theorem, Wiles' main contribution was to show modularity of semistable elliptic curves defined over $\Q$, which concerns Step $2$ above. His work opened the door for the era of modularity lifting theorems, a domain that has seen remarkable progress in the last few decades. Even if Step $2$ uses a heavy machinery, many modularity results are well understood nowadays, and they are often sufficient for Diophantine applications. 
\medskip

On the other hand, Steps $1$ and $3$ are quite challenging to establish in full generality. In most of the already treated examples, the representation $\rho$ given by Step $1$ arises from the action of $\Gal(\Kglbar / \Kgl)$ on the $p$-adic Tate module of an elliptic curve $E(a, b, c) / \Kgl$, usually called a \textit{Frey} curve. There is no generic recipe to construct Frey elliptic curves, and most of the known cases in the literature are obtained using \textit{ad hoc} arguments. On top of that, the lack of theoretical techniques for distinguishing residual Galois representations makes Step $3$ difficult in general. Moreover, the latter requires, among other things, numerical computations of spaces of Hilbert newforms, whose dimensions grow very quickly with the size of the parameters in the equation. This algorithmic challenge is one of the main bottlenecks for applying the modular method in many situations.

\subsubsection*{Darmon's program: strengthening the modular method}

In \cite{Darmon00}, Darmon described a program that aims at solving some of the difficulties raised by the modular method. We refer the curious reader to \cite{ChenKoutsianas25} for a detailed survey on this program. Among other questions, he studied how to successfully perform Step~$1$ in the modular method. One of his main contributions was to define what a suitable Frey object should be, and which properties should be satisfied by the representation $\rhobar$ introduced above. Its core idea is the following:

\begin{main_idea*}
	For each signature $(p, q, r)$, construct a so-called \textnormal{Frey representation} $\rhobar_s : G_{\Kgl(s)} \rightarrow \GL_{2}(\overline{\Fp})$ having ``nice" ramification properties. Here $\Kgl(s)$ is the rational function field over a totally real field $\Kgl$. The representation $\rhobar$ from above is then obtained by specialising the indeterminate $s$ at some algebraic number $s_0$ whose value depends on the considered equation \eqref{eq: GFE pqr} and its putative solution. 
\end{main_idea*}

Darmon classified and constructed Frey representations for every signature. As in the classical modular method, his examples arise from the geometry of abelian varieties. More specifically, he introduced curves of genus $\geq 1$, defined over function fields, whose Jacobians have \textit{real multiplication} (hereafter abbreviated as \RM). On a different direction, Darmon also envisioned how to use different Frey objects to propagate modularity from one another. We will give further details about this in subsection~\ref{sect: Modularity}.
\medskip

For the signature $(p,p,r)$, Darmon defines two curves $\Cpm(s) / \Q(s)$ that are \textit{hyperelliptic}. For the signature $(r,r,p)$, he constructed two \textit{superelliptic} curves defined over $\Q(s)$. Billerey, Chen, Dieulefait and Freitas introduce in \cite{BCDF23} a \textit{hyperelliptic} curve $C_r(s) / \Q(s)$ whose Jacobian gives rise to a Frey representation of signature $(r, r, p)$. It is a quadratic twist of Darmon's hyperelliptic curve $\Cms(s)$ for the signature $(p,p,r)$ (see \S \ref{sect: Frey curves rrp} for further details). In \cite{Pacetti25}, Pacetti uses the theory of hypergeometric motives to explain and give heuristics about the intimate relationship between all these geometric objects. We note that, when $r \geq 5$ (and away from the case $|AB| = 1$ for the signature $(r,r,p)$), these are the only known Frey objects, and, in particular, there are no known available elliptic curves. 
\medskip

The above discussion illustrates that the construction of Frey representations is not unique. Curves having different geometric natures may give rise to the same Frey representation. Most of the arithmetic properties of the latter can be read off from the geometry of the underlying curve, so understanding in detail this geometry is essential. After specialising the indeterminate $s$ at some $s_0 \in \Qbar$, one obtains a curve defined over a number field. Not all kinds of curves are equally understood. For some of them (like hyperelliptic ones), many effective and computational techniques have appeared in the last few years. For instance, the theory of cluster pictures \cite{M2D2}, developed by Dokchitser, Dokchitser, Maistret and Morgan provides combinatorial tools to understand the local behaviour of a hyperelliptic curve at odd places of bad reduction. We refer the curious reader to \cite{ACIKMM}, where Curc\'{o}-Iranzo, Khawaja, Maistret, Mocanu and the author use the machinery of clusters to compute conductor exponents for certain Frey hyperelliptic curves.

\subsubsection*{Our contribution and results} 

A major contribution of this paper is the introduction of a common framework to simultaneously study the generalised Fermat equations \eqref{eq: GFE ppr} and \eqref{eq: GFE rrp}, of respective signatures $(p,p,r)$ and $(r,r,p)$. The core idea is the following: to obtain Frey hyperelliptic curves defined over $\Kgl$ for these signatures, we consider a single curve defined over $\Kgl(s)$, that we specialise at different algebraic numbers. 
\medskip 

Following \cite{Darmon00}, we attach two Frey hyperelliptic curves $\Cpmabc / \Q$ to a putative solution $(a, b, c)$ to \eqref{eq: GFE ppr}. Similarly, following \cite{BCDF23} and generalising a construction by Kraus, we attach a Frey hyperelliptic curve $\Crabc / \Q$ to a putative solution to \eqref{eq: GFE rrp}. The former ones are quadratic twists of specialisations of the curves $\Cpm(s) / \Q(s)$. The latter one is a quadratic twist of a specialisation of $C_r(s) / \Q(s)$. But as explained above, the curve $C_r(s)$ is a quadratic twist of $\Cms(s)$, so the curve $\Crabc$ is also a quadratic twist of a specialisation of $\Cms(s)$. This elementary observation, which is made precise in the proposition below, is the starting point for the whole discussion in this paper.

\begin{proposition}[= Proposition~\ref{prop: right choices s0 deltaQ} + Lemma~\ref{lem: s0 deltaQ fill hypoth}]\label{prop: right choices intro}
	The following properties hold.	
	\begin{enumerate}[itemsep=0pt, topsep= 5pt]
		\item Assume that $(a, b, c)$ is a primitive non-trivial solution to \eqref{eq: GFE ppr}. There exists some $s_0 \in \Q$ and some $\deltaQ \in \Z$ depending on $(a, b, c)$ and \eqref{eq: GFE ppr}, such that $\Cpmabc$ is the quadratic twist by $\deltaQ$ of the specialisation $\Cpm(s_0)$.
		
		\item Assume that $(a, b, c)$ is a primitive non-trivial solution to \eqref{eq: GFE rrp}. There exists some $s_0 \in \Qbar$ and some $\deltaQ \in \Zbar$ depending on $(a, b, c)$ and \eqref{eq: GFE rrp}, such that $\Crabc$ is the quadratic twist by $\deltaQ$ of the specialisation $\Cms(s_0)$.
	\end{enumerate}
	In Table~\ref{table: Main values}, we display the corresponding values of $s_0$ and $\deltaQ$ for each of the signatures. Moreover, there are various hypotheses commonly satisfied for both choices of $s_0$ and $\deltaQ$ (\cf \ \S \ref{sect: Common framework}).
\end{proposition}

Proposition~\ref{prop: right choices intro} provides a common framework for a uniform study of the two GFEs \eqref{eq: GFE ppr}, \eqref{eq: GFE rrp}. In subsection \ref{sect: Common framework}, we give axiomatic definitions for $s_0, \deltaQ$, and we list the mentioned hypotheses that they should satisfy. Throughout this paper, we state all our results in terms of $s_0$ and $\deltaQ$. Nevertheless, we leave our own \textit{``breadcrumb trail"}: after every important statement, we present the corresponding result for both generalised Fermat equations, a putative solution and the associated Frey object. Even if the construction of $\Cpm(s_0)$ depends (tautologically) on $s_0$, most of its properties are encoded in the behaviour of $s_0(s_0-1)$. This is beneficial for the sake of effectiveness: for instance, the choice of $s_0$ for the signature $(r, r, p)$ depends in a complicated way on the coefficients $A, B, C$ and the solution $(a, b, c)$. Nevertheless, Table~\ref{table: Main values} shows that $s_0(s_0-1)$ has a simpler expression in terms of the mentioned parameters.
\medskip

Throughout the paper, we let $\Kgl \coloneqq \Q(\zr)^{+}$ be the maximal totally real subfield of the cyclotomic field $\Q(\zr)$. Consider the base-changed curves $\Cpm \coloneqq \Cpm(s_0)^{(\deltaQ)} \times_{\Q} \Kgl$: their Jacobians $\Jpm \coloneqq \Jac(\Cpm)$ have real multiplication \RM \ by $\Kgl$ (see Theorem~\ref{thm: Jpm have RM}). Following the work of Ribet \cite{Ribet76}, we construct $2$-dimensional representations $\rhopm : \Gal(\Kglbar / \Kgl) \rightarrow \GL_{2}(\Klam)$, where $\lambda$ is a finite place of $\Kgl$ (see \S \ref{sect: Decomp Tate module}). When $\lambda$ divides the prime number $p$ appearing in the signature, we will use the reductions of these representations (and their twists) to perform the modular method as depicted above.
\medskip 

Serre, Tate \cite{SerreTate}, and Grothendieck \cite{SGA7} illustrated that much of the arithmetic information of $\rhopm$ is given by the geometric behaviour of $\Jpm / \Kgl$. With this in mind, we perform a careful analysis of the reduction types of their Néron models. Moreover, we introduce a new parameter $\deltaK \in \OK$, and consider the twisted Jacobian $\Jpmtw$. The purpose of this is to obtain an abelian variety having better reduction properties (for some choices of $\deltaK$, $\Jpmtw$ is semistable, whereas $\Jpm$ is not). 

\begin{theorem}[= Theorems~\ref{thm: reduction types Jms} + \ref{thm: reduction types Jpl}]\label{thm: Red types intro}
	Let $q$ be a rational prime, $\frakq$ a finite place of $\Kgl$ dividing $q$, and define $\nu_q$ to be the $q$-adic valuation of $s_0 (s_0-1)$. The decision trees depicted in Figures~\ref{fig: Reduction Jmodmstw} and \ref{fig: Reduction Jmodpltw} describe the reduction type of the Néron models of $\Jpmtw$ at $\frakq$ in terms of $\nu_q, \deltaQ$ and $\deltaK$.
\end{theorem}

In Example~\ref{examp: Reduction types Fermat}, we illustrate how this theorem applies when considering the specific choices of $s_0, \deltaQ$ for the signatures $(p,p,r)$ and $(r,r,p)$. To establish Theorem~\ref{thm: Red types intro}, we first exhibit algebraic expressions for the roots of the defining polynomials of $\Cpm(s_0)$, and we deduce from this the discriminants of the models and theirs twists (see Example~\ref{examp: Discriminants Fermat} for the particular arithmetic statements). Knowing the places of bad reduction, we make use of the following ``geometric trichotomy": if $A$ is an abelian variety defined over a local field which has \RM, then its Néron model has either good, totally toric, or totally unipotent reduction (\cf \ Proposition~\ref{prop: reduction RM}). This phenomenon, satisfied in particular by $\Jpmtw$, simplifies the analysis of the reduction types of its Néron model (hence the proof of Theorem~\ref{thm: Red types intro}).
\medskip 

In order to understand these reduction types, we make use of the machinery of cluster pictures \cite{M2D2}. At even places, we study the geometry of different models of $\Cpmtw$, if they are semistable, and, if not, over which extension they become so. We also use the recent work of Gehrunger \cite{Gehrunger25} and Dokchitser--Morgan \cite{DokchitserMorgan23}, which help understanding the stable reduction of a hyperelliptic curve in residue characteristic $2$, yielding in certain situations the toric rank. Our analysis relies on some $2$-adic conditions on $s_0 (s_0-1)$: it covers the ``generic" cases, and is therefore sufficient for our Diophantine purposes. At odd places, we use the results on cluster pictures from \cite{M2D2}, in the style of \cite{ACIKMM}. We study the ramification properties of the splitting field of the defining polynomials, and we use the criteria on clusters (Theorems~\ref{thm: Crit clusters good} and \ref{thm: Crit clusters semist}), to check if $\Jpmtw$ have good or bad semistable reduction (hence totally toric): otherwise they have totally unipotent reduction. We note that our approach differs from the one in \cite{ACIKMM}, as we do not use clusters to compute tame conductors, but only the reduction types of the Néron models.
\medskip

Having understood the geometry of the (twisted) curves $\Cpmtw$ and their Jacobians, we move onto the study of the $2$-dimensional representations $\rhopmK$. These form a compatible system of Galois representations. Since we aim at performing the modular method with $\rhopm$ and its twist $\rhopmK$, it is essential to know that:

\begin{theorem}[= Theorems~\ref{thm: Modularity Jms} + \ref{thm: Modularity Jpl}]
	The system of Galois representations $(\rhoms)_{\lambda}^{ }$ arises from a Hilbert newform over $\Kgl$. If $v_r(s_0 (s_0-1)) > 2$, the same holds for $(\rhopl)_{\lambda}^{ }$.
\end{theorem}

Our proof of this result builds upon Darmon's panoramic view for propagating modularity (\cf \ \cite[Figure 1]{Darmon00}). We show that the particular representation $\rho_{\Jpm, \, \frakr}$ is modular, where $\frakr$ is the unique place of $\Kgl$ above $r$. We employ powerful tools, such as Serre's modularity conjecture and modularity lifting theorems. At different points, we make use of the knowledge of the reduction types from Theorem~\ref{thm: Red types intro} to establish desired properties about the considered representations.
\medskip

Modularity theorems ensure that $(\rhopmK)_{\lambda}^{ }$ arises from a Hilbert newform whose level is given by the conductor of the system $(\rhopmK)_{\lambda}^{ }$. It is therefore crucial to understand this conductor in detail, as it will be a key input for level lowering and the elimination step. We provide an explicit description of the Artin conductor of the system at every finite place:

\begin{theorem}[= Theorems~\ref{thm: Conductor rhomsK} + \ref{thm: Conductor rhoplK}]\label{thm: Conductor intro}
	Let $\frakq$ be a finite place of $\Kgl$. The Artin conductor of $\rhopmK$ restricted to the decomposition group $D_{\frakq} \simeq \Gal(\Kqbar / \Kq)$ is described in Tables~\ref{table: Cond Jmstw} and \ref{table: Cond Jpltw}.
\end{theorem}

In Example~\ref{examp: Conductor Fermat}, we illustrate how to recover from Theorem~\ref{thm: Conductor intro} the Artin conductors for the particular choices of $s_0, \deltaQ$ displayed in Table~\ref{table: Main values}. The proof of this theorem is largely based on the description of the reduction types given in Theorem~\ref{thm: Red types intro}. Indeed, Grothendieck \cite{SGA7} proved that the tame part of the Artin conductor is encoded in the reduction type of the Néron model. At even places, we restrict ourselves to the cases of tame potential semistable reduction, in which there is no wild conductor. At odd places, the theory of cluster pictures yields the wild conductor by means of algebraic number theory. The content of subsections~\ref{sect: Ramification} and \ref{sect: Cond & local types}, where this is discussed, is largely inspired by \cite{ACIKMM}, and provides a generalisation of some of its results.
\medskip

Along the way, we describe the inertial local types of the Weil--Deligne representations attached to $\rhopmK$. This local study of the action of inertia is a key datum for proving absolute irreducibility of the residual representations. Building upon the work of \cite[\S 7]{BCDF23} on finiteness of Galois representations, we finally establish our desired level lowering result (Theorem~\ref{thm: Level lowering}). For any prime number $p$ and any place $\frakp \mid p$, it asserts the existence of a newform $g$ such that $\rhopmpbar \simeq \overline{\rho}_{g, \frakP}$ for some place $\frakP \mid p$ in the field of coefficients $\Kgl_g$. The key feature of $g$ is that, when we choose $s_0, \deltaQ$ as in Table~\ref{table: Main values}, the primes dividing the level of $g$ do not depend on the putative solution, and are only expressed in terms of the parameters of the GFE. 
\medskip 

We finally discuss the elimination step, in order to prove Theorems~\ref{thm: Example pp5} and \ref{thm: Example 55p}. As explained above, this requires an algorithmic implementation to discard isomorphisms like $\rhopmpbar \simeq \overline{\rho}_{g, \frakP}$. We do so by comparing traces of Frobenius under these representations. We specify the values of the exponent $r$ and the coefficients $A, B, C$ as in the equations in Theorems~\ref{thm: Example pp5} and \ref{thm: Example 55p}, hence giving numerical values for the level of $g$. We then use the{\fontfamily{lmtt} \selectfont Magma} software \cite{Magma} to compute spaces of Hilbert newforms, and eliminate Galois orbits of newforms by discarding isomorphisms of representations. We provide a \Magma \ package that allows to perform this elimination step for many choices of parameters $r, A, B, C$, thus allowing to solve several instances of families of GFEs. This \Magma \ package is available at \cite{Git}.

\subsubsection*{Acknowledgments}

The present paper is based on the work that I have done during my doctoral studies at Université Clermont Auvergne. I wholeheartedly thank Nicolas Billerey for proposing this research topic, and for helpful comments on earlier versions of this manuscript. Special thanks go to Mar Curc\'{o}-Iranzo, Maleeha Khawaja, C\'{e}line Maistret and Diana Mocanu for fruitful discussions about clusters; to Enric Florit for helpful discussions about endomorphism algebras; and to Tim Gehrunger for taking the time to explain details about the reduction of hyperelliptic curves in residue characteristic $2$. I would also like to thank Ariel Pacetti: a part of this work answers a question he raised. Last, but not least, many thanks go to Lucas Villagra Torcomian and Pedro José Cazorla García for interesting conversations.

% -- % -- % -- % -- % -- % -- % -- % -- % -- % -- % -- % -- % -- % -- % -- % -- % -- % -- % -- % -- % -- % -- % -- % -- % -- % -- % -- % -- % -- % -- % -- % -- % -- % -- % -- % -- % -- % -- % -- 

\section[Background and notation]{Background and notation}

\subsection{Notation}\label{sect: Notation}

For any field $K$, we let $\overline{K}$ be an algebraic closure of $K$, $\sep{K}$ be a separable closure included in $\overline{K}$. We let $G_K \coloneqq \Gal(\sep{K} / K)$ be the absolute Galois group of $K$, endowed with the profinite topology. If $v$ is a valuation defined on $\sep{K}$, we will say that $v$ is normalised with respect to $K$ if $v(K\st) = \Z$.

For a fixed rational prime $q$, we let $\Qq$ be the field of $q$-adic numbers, and $v_q$ be the valuation on $\Qqbar$ normalised with respect to $\Qq$. Given a number field $\F$, we denote by $\O_{\F}$ its ring of integers. For a prime ideal $\frakq$ of $\O_{\F}$, we let $\Fq$ be the completion of $\F$ at $\frakq$, and we let $\vq$ be the valuation on $\Fqbar$ normalised with respect to $\Fq$. We say that $\frakq$ is even (resp. odd) if $\frakq \mid 2$ (resp. $\frakq \nmid 2$).

If $(K, \frakm)$ is a local field, endowed with a normalised valuation $\vm$, we write $\Om$ for the ring of integers of $K$, and $\Ffm$ for its residue field. We still denote by the symbol $\vm$ the valuation defined on $\sep{K}$ that matches $\vm$ over $K$. 
We let $\unr{K}$ be the maximal unramified extension of $K$, $I_K = \Gal(\sep{K} / \unr{K})$ its inertia group, and $I_K^{t}$ (resp. $I_K^{w}$) the tame (resp. wild) inertia group. We let $\Frob_{\frakm} \in G_{\Ffm}$ be the arithmetic Frobenius, and $\phim \in G_K$ be a geometric Frobenius lift: its reduction in $G_{\Ffm}$ is $\Frob_{\frakm}^{-1}$. We let $W_K$ be the Weil group of $K$, which consists of the elements of $G_K$ that act on $\Ffm$ as a power of $\Frob_{\frakm}$. The Weil group fits into the short exact sequence
\begin{equation*}
	1 \longrightarrow I_K \longrightarrow W_K \longrightarrow \Frob_{\frakm}^{\Z} \longrightarrow 1.
\end{equation*}
If $\F$ is a number field and $\frakq$ is a finite place of $\F$, we will simply denote by $D_{\frakq} \simeq G_{\F_{\frakq}^{ }}$ a decomposition group of $G_{\F}$ at $\frakq$ and by $I_{\frakq}$ its inertia subgroup. We denote by $W_{\frakq}$ for the Weil group of $D_{\frakq}$.
\medskip 
 
We fix a prime number $r \geq 5$. Let $\zr \in \Qbar$ be a primitive $r$-th root of unity, and let $\Q(\zr)$ be the $r$-th cyclotomic field. For any $0 \leq j \leq \frac{r-1}{2}$, we let $\omega_j \coloneqq \zr^{j} + \zr^{-j}$, $\tau_j \coloneqq \zr^j - \zr^{-j}$ and, for simplicity, we write $\omega \coloneqq \omega_1$, and $\tau \coloneqq \tau_1$.  We denote by $\Kgl \coloneqq \Q(\omega)$ the maximal totally real subfield of $\Q(\zr)$. The extension $\Kgl / \Q$ has degree $\frac{r-1}{2}$, and is defined by the irreducible polynomial
\reqnomode
\begin{equation}\label{eq: Factor h_r}
	h_r(x) \coloneqq \prod_{j = 1}^{\frac{r-1}{2}} (x - \omega_j) \in \Z[x].
\end{equation}
Moreover, the extension $\Kgl / \Q$ is totally ramified at $r$, and unramified at any other rational prime. We denote by $\frakr = (2 - \omega)$ the unique prime ideal of $\OK$ dividing $r$. Finally, we denote by $\phi_r(X, Y) \in \Z[X, Y]$ the polynomial defined by \vspace{-0.75em}
\begin{equation*}
	\phi_r(X, Y) \coloneqq \frac{X^r + Y^r}{X + Y} = \sum_{j = 0}^{r-1} (-1)^j X^{r-1-j} Y^j.
\end{equation*}

\subsection{$\lambda$-adic and Weil--Deligne representations}\label{sect: ladic WD reps}

We now recall essentials facts about $\lambda$-adic and Weil--Deligne representations that will be used throughout this article. We refer the reader to \cite[\S 8]{Deligne73} or \cite{Tate79} for details, and to \cite{Rohrlich94} for a motivation towards the study of elliptic curves and abelian varieties.
\medskip

Let $(K, \frakm)$ be a local field, we adopt the notation introduced in \S \ref{sect: Notation}. Fix a prime number $\ell \neq \charact(\Ffm)$. Let $\Elam$ be the completion of a number field $E$ at a finite place $\lambda$ dividing $\ell$.

\begin{definition}
	A $\lambda$-adic representation of $G_K$ is a continuous homomorphism $\rho : G_K \rightarrow \GL(V)$, where $V$ is finite dimensional $\Elam$-vector space endowed with the $\lambda$-adic topology. We denote by $\Rep_{\Elam}(G_K)$ the category of $\lambda$-adic representations of $G_{K}$.
\end{definition}

Examples of $\lambda$-adic representations naturally arise in the \'{e}tale cohomology of smooth proper varieties, but also associated to modular forms. We refer the curious reader to \cite{Taylor} for a survey.

\begin{definition}
	Let $V$ be a vector space over a field $F$ of characteristic zero. We define an $F$-linear Weil--Deligne (\WD) representation as a pair $(\rho, N)$, where $\rho : W_K \rightarrow \GL(V)$ is a homomorphism continuous with respect to the discrete topology on $V$, and $N$ is an endomorphism of $V$, called the monodromy operator, satisfying the equality $\rho N \rho^{-1} = \omega_K N$. Here $\omega_K : W_K \rightarrow F\st$ denotes the unramified character of $W_K$ such that $\omega_K(\phim) = \abs{\Ffm}^{-1}$.
	
	A morphism of \WD-representations is a morphism between the underlying representations that commutes with the monodromy operators. We denote by $\Rep_{F}(W_K')$ the category of $F$-linear \WD-representations.
\end{definition}

\begin{remark}
	The continuity of a \WD-representation $(\rho, N)$ implies that $\rho(I_K)$ is finite. Moreover, the equality satisfied by the monodromy operator $N$ implies that it is a nilpotent endomorphism.
\end{remark}

The following result establishes the link between $\lambda$-adic representations of $G_K$ and Weil--Deligne representations. Building on Grothendieck's $\ell$-adic monodromy theorem, one can attach a \WD-representation to any $\lambda$-adic representation of $G_K$. This construction is completely explicit, and the curious reader can find details in \cite[\S 8]{Deligne73}, \cite[Appendix]{SerreTate}. 

\begin{theorem}
	There is a fully faithful functor
	\begin{align*}
		\Rep_{\Elam}(G_K) & \longrightarrow \Rep_{\Elam}(W_K') \\
		\rholam & \longmapsto \left( \W(\rholam), N_{\rholam} \right)
	\end{align*}
\end{theorem}

\begin{remark}
	For any $\lambda$-adic representation $\rholam \in \Rep_{\Elam}(G_K)$, $\rholam(I_K)$ is finite if and only if $N_{\rholam} = 0$. Whenever this is the case, the \WD-representation associated to $\rholam$ simply is $(\rholam, 0)$.
\end{remark}

The interest of considering \WD-representations is that their definition uses only the discrete topology on $V$. Therefore, \WD-representations are well adapted to modify the field of coefficients of $V$, namely shifting from $\Ql$ or $\Elam$ to $\C$.

\begin{definition}
	For any $\rholam \in \Rep_{\Elam}(G_K)$ and any embedding $\iota : \Elam \hookrightarrow \C$, we denote by $\W_{\iota}(\rholam)$ the complex \WD-representation obtained by extending the scalars of $\W(\rholam)$ through $\iota$.
\end{definition}

In this article, we will be mainly interested in $2$-dimensional complex \WD-representations. These can be explicitly classified, and we refer the reader to \cite{Tate79, DembeleFreitasVoight} for a detailed discussion.

\begin{proposition}\label{prop: Classif WD}
	Let $(\rho, N) \in \Rep_{\C}(W_K')$ be a $2$-dimensional complex \WD-representation. Then $(\rho, N)$ arises up to isomorphism from one of the following three possibilities. 
	\begin{enumerate}[leftmargin=*, itemsep=0pt]
		\item {\normalfont Steinberg representations:} There is a continuous character $\chi : W_K \rightarrow \C\st$ with open kernel such that $\rho = \chi \otimes (\omega_K^{ } \oplus \mathbbm{1})$, and $N = \begin{pmatrix}
			0 & 1 \\
			0 & 0
		\end{pmatrix}$. 
	
		\item {\normalfont Principal series:} There are two continuous characters $\chi_1 , \chi_2 : W_K \rightarrow \C\st$ with open kernels satisfying $\chi_1 \chi_{2}^{-1} \neq \omega_K^{\pm 1}$, such that $\rho = \chi_1 \oplus \chi_2$, and $N = 0$. 
		
		\item {\normalfont Supercuspidal representations:} $\rho$ is an irreducible $2$-dimensional representation of $W_K$, and $N = 0$. We say that $\rho$ is non-exceptional if its projective image is dihedral. Otherwise, it is called exceptional, and it has projective image $A_4$ or $S_4$.
	\end{enumerate}
\end{proposition}

\begin{remark}
	Any non-exceptional supercuspidal \WD-representation is the induction of a character of a quadratic extension of $K$. For further details about properties of supercuspidal representations, we refer the reader to \cite[\S 2.3]{DembeleFreitasVoight}.
\end{remark}

When considering $\lambda$-adic and \WD-representations, one is particularly interested in the action of the inertia group $I_K$. For this purpose, we introduce:

\begin{definition}
	An inertial local type is an equivalence class $[\rho, N]$ of \WD-representations under the relation $(\rho, N) \sim (\rho', N')$ if and only if $(\restr{\rho}{I_K}, N) \sim (\restr{\rho'}{I_K}, N')$.
\end{definition}

\begin{remark}\label{rmk: Supercuspidal WD reps}
	There is an explicit classification of local inertial types (see \cite[Proposition 2.4.1]{DembeleFreitasVoight}), which is compatible with the classification given in Proposition~\ref{prop: Classif WD}. We will say that $[\rho, N]$ is Steinberg, principal series or supercuspidal if it is the restriction of such a \WD-representation.
\end{remark}

For any real number $u \geq -1$, we let $G_K^{u}$ be the absolute $u$-th ramification group of $K$ in upper numbering (see \cite[\S 3]{Ulmer16}). These groups produce a filtration of $G_K$: we have $G_K^{-1} = G_K$, $G_K^0 = I_K$ and $\bigcup_{u > 0} G_K^{u} = I_K^{w}$. For any $\lambda$-adic representation $\rholam : G_K \rightarrow \GL(V)$, and any subgroup $G < G_K$, we let $V^G$ be the subspace of $V$ fixed by $\rholam(G)$. 

\begin{definition}
	Let $\rholam : G_K \rightarrow \GL(V)$ be a $\lambda$-adic representation. We define the tame and wild (or Swan) conductor of $\rholam$ as
	\begin{equation*}
		\condtame{\rholam} \coloneqq  \int_{-1}^{0} \codim_{\Elam} V^{G_K^u}\, du \quad \text{ and } \quad \condwild{\rholam}\coloneqq\int_{0}^{\infty} \codim_{\Elam} V^{G_K^u}\, du.
	\end{equation*}
	We define the Artin conductor of $\rholam$ by $\condexp{\rholam} \coloneqq \condtame{\rholam} + \condwild{\rholam}$. 
\end{definition}

\begin{remark}\label{rmk: Cond WD lambda match}
	One can show that the Artin conductor $\condexp{\rholam}$ is an integer (see \cite[VI, \S2]{Serre79}). Similarly, one can define the conductor of a \WD-representation. The important feature of this definition, is that, for any $\rholam \in \Rep_{\Elam}(G_K)$, the conductor of $\W(\rholam)$ matches the Artin conductor $\condexp{\rholam}$ introduced above (see \cite[\S 8]{Ulmer16}).
\end{remark}

\begin{lemma}\label{lem: Behaviour nwild tame ext}
	Let $L / K$ be a tame extension with ramification index $e_{L/K}$. Let $\rholam : G_K \rightarrow \GL(V)$ be a $\lambda$-adic representation. The wild conductor of the restriction $\restr{\rholam}{G_L}$ is given by
	\begin{equation*}
		\condwild{\restr{\rholam}{G_L}} = e_{L/K} \, \condwild{\rholam}.
	\end{equation*}
\end{lemma}

\begin{proof}
	This is a straightforward generalisation of \cite[Lemma 2.2]{ACIKMM}.
\end{proof}

\subsection{Abelian varieties with real multiplication}\label{sect: AV with RM}

In this subsection we discuss various properties of abelian varieties with real multiplication. Our main references are \cite{Ribet76} and \cite{Melninkas}.
\medskip

The following discussion does not require the ground field to be local. That is why, we let $L$ be any field, and consider an abelian variety $A / L$ of dimension $g$.

\begin{definition}
	Let $F$ be a totally real field of degree $g$ over $\Q$. We say that $A$ has real multiplication by $F$ (or simply \RM) if there exists an isomorphism $F \simeq \End_{L}(A) \otimes \Q$.
\end{definition}

\begin{remark}
	Abelian varieties with \RM \ are instances of so-called {\normalfont abelian varieties of $\GL_{2}$-type}. Some of the properties below rely only on $A$ being of $\GL_{2}$-type. Nevertheless, $F$ being totally real has stronger consequences (see Remark~\ref{rmk: Counterexample GL2} below). Moreover, the abelian varieties considered in this article will have \RM, so we focus on this case. 
\end{remark}

\begin{lemma}\label{lem: quad twist has RM}
	Let $d \in L\st \setminus (L\st)^2$, and $\chi_d : G_L \rightarrow \lbrace \pm 1 \rbrace$ be the character corresponding to $L(\sqrt{d})$. If $A / L$ has \RM \ by $F$, then $A^{(d)}$, its quadratic twist by $\chi_d$, also has \RM \ by $F$.
\end{lemma}

\begin{proof}
	The character $\chi_d$ acts on $A^{(d)}$ (and thus on $A^{(d)} \times_L L(\sqrt{d}))$ as $-1 \in \Aut_L(A^{(d)})$	. But the action of $F$ on $\End_L(A)$ commutes with the automorphism $-1 \in \Aut_L(A)$, so the action of $F$ on $A \times_{L} L(\sqrt{d}) \simeq A^{(d)} \times_{L} L(\sqrt{d})$ is actually defined over $L$.
\end{proof}

\subsubsection{Decomposition of the Tate module}\label{sect: Decomp Tate module}

Fix a prime number $\ell \neq \charact(L)$ . Let  $\TlA \coloneqq \varprojlim_{n \geq 1} A \left[\ell^n \right]$ be the $\ell$-adic Tate module of $A$, and $\VlA \coloneqq \TlA \otimes \Ql$. It is well-known that $\VlA$ is a $\Ql$-vector space of dimension $2g$. The absolute Galois group $G_L$ acts continuously on $\VlA$, giving rise to a representation 
\begin{equation*}
	\rholA : G_{L} \longrightarrow \Aut_{\Ql} (\VlA) \simeq \GL_{2g}(\Ql).
\end{equation*}
Assume now that $A$ has \RM \ by a totally real field $F$, and define $\Fl \coloneqq F \otimes \Ql$. This is a semisimple $\Ql$-algebra that acts on $\VlA$, and the action of $G_L$ on $\VlA$ is $\Fl$-linear. For any place $\lambda$ of $F$ dividing $\ell$, define $\VlamA \coloneqq \VlA \, \otimes_{\Ql} \Flam$. In \cite[II \S 1]{Ribet76}, Ribet proves that $\VlA$ is a free $\Fl$-module of rank $2$, so $\VlamA$ is a $2$-dimensional $\Flam$-vector space. Again, the action of $G_L$ on $\VlamA$ is $\Flam$-linear, so we obtain a continuous representation
\begin{equation*}
	\rholamA : G_L \longrightarrow \Aut_{\Flam}(\VlamA) \simeq \GL_{2}(\Flam).
\end{equation*}
Shimura proved in \cite[\S11.10]{Shimura67} that the $\rholamA$'s form a weakly compatible system of Galois representations in the sense of \cite[Definition 3.12]{Bockle}. By the work of Fontaine \cite{Fontaine94}, this is moreover a strictly compatible system (see also \cite[Proposition 2.8.1]{BoxerCalegariGeePilloni} for further details).

Denote by $\Res_{\Flam / \Ql}(\VlamA)$ the restriction of scalars of $\Vlam$ to $\Ql$, and by $\Res_{\Flam / \Ql}(\rholamA)$ the associated representation of $G_L$. The isomorphism of $\Ql$-vector spaces $\Fl \simeq \prod_{\lambda \mid \ell} \Flam$ induces an isomorphism of $\Ql$-representations
\begin{equation}\label{eq: Decomp l adic lambda adic}
	\rholA \simeq \bigoplus_{\lambda \mid \ell} \Res_{\Flam / \Ql}(\rholamA).
\end{equation}
Finally, we define $\overline{\rho}_{A, \, \lambda} : G_L \rightarrow \GL_{2}(\Ff_{\lambda})$ to be the representation obtained by considering the action of $G_L$ on the semi-simplification of the module $\TlA \otimes_{\O_F} \Ff_{\lambda}$ (see \cite[\S II.2]{Ribet76}).

\subsubsection{Reduction types of the Néron model}

We now come back to the setting of \S \ref{sect: Notation}, by assuming that $L = K$ is a local field, and recover all the notation introduced above. In particular, $\frakm$ denotes the maximal ideal in the ring of integers of $K$.
\medskip

Let $A / K$ be an abelian variety (not necessarily having \RM). Let $\Amod \rightarrow \Spec(\Om)$ be the Néron model of $A/K$ (see \cite{BLR} for a detailed discussion). Denote by $\Amodm$ its special fiber, and by $\Amodm^{0}$ the connected component of the identity in it. By Chevalley's theorem (see \cite{Milne17}), there exists a torus $T / \Ffm$, a unipotent group $U / \Ffm$ and an abelian variety $B / \Ffm$ that fit in the short exact sequence
\begin{equation}\label{eq: ses Chevalley}
	1 \longrightarrow T \times U \longrightarrow \Amodm^{0} \longrightarrow B \longrightarrow 1.
\end{equation}
The datum of $T, U$ and $B$ is called the \textit{reduction type} of $\Amod$ at $\frakm$. We say that $\Amod$ is semistable if the unipotent group $U$ in \eqref{eq: ses Chevalley} is trivial. The semistable reduction theorem \cite[Exposé IX]{SGA7} states that any abelian variety over a local field has potential semistable reduction.
\medskip

Abelian varieties with real multiplication have restricted possibilities for their reduction types. 

\begin{proposition}\label{prop: reduction RM}
	Assume that $A$ has \RM. Then $\Amodm^{0}$ is either an abelian variety, a torus or a unipotent group, meaning that exactly one among $T, U$ and $B$ as in \eqref{eq: ses Chevalley} is non-trivial.
\end{proposition}

\begin{proof}
	This is \cite[Proposition 3.6.1]{Ribet76} for $\Amod$ semistable, and \cite[Proposition 2.6]{Melninkas} in general.
\end{proof}

\begin{remark}\label{rmk: Counterexample GL2}
	Ribet proved in \cite[\S 3.6]{Ribet76} that any abelian variety of $\GL_{2}$-type having semistable reduction has either good or totally toric reduction. However, Proposition~\ref{prop: reduction RM} does not hold if we only assume $A/K$ to be of $\GL_{2}$-type without having \RM. The following example, which was communicated to the author by Enric Florit, illustrates so.
	
	Consider the Jacobian of the genus $2$ modular curve $X_1(13) / \Q$ (labelled on the LMFDB \cite{lmfdb} as $169.a.169.1$). Its endomorphism algebra is $\End_{\Q}(\Jac(X_1(13))) \otimes \Q \simeq \Q(\sqrt{-3})$, so $\Jac(X_1(13)) / \Q$ is of $\GL_2$-type, but does not have \RM. Looking at the cluster picture of $X_1(13)$ at $p = 13$, one can check that the reduction type of the Néron model of $\Jac(X_1(13))$ is an extension of an elliptic curve by a unipotent group of dimension $1$ (\cf \ \cite[\S 5]{Hyperuser}).	Thus, Proposition~\ref{prop: reduction RM} does not hold for $\Jac(X_1(13)) / \Q$.
\end{remark}

\begin{definition}
	If $A$ has \RM, we say that $\Amod$ has good (resp. \textit{toric}, or \textit{unipotent}) reduction at $\frakm$ if $\Amodm^{0}$ is an abelian variety (resp. a torus, or unipotent).
\end{definition}

\begin{remark}\label{rmk: No mixed red types}
	In general, the Néron model of an abelian variety has ``mixed reduction" types, meaning that two or more among $T, U$ and $B$ as in \eqref{eq: ses Chevalley} are non-trivial. Whenever $\Amod_{\frakm}^{0}$ is a torus (resp. unipotent), we say that $\Amod$ has \textbf{totally} toric (resp. \textbf{totally} unipotent) reduction. In this paper we will only consider abelian varieties with \RM \ so, to simplify the terminology, we drop the adverb ``totally".
\end{remark}

As explained in \cite[Exposé IX, Théorème 3.6]{SGA7}, any abelian variety (having \RM \ is not necessary) over a local field attains semistable reduction over a finite extension of $K$. With this in mind, we introduce. 

\begin{definition}\label{def: semist defect}
	Let $A / K$ be an abelian variety, consider $M / \unr{K}$ be the minimal extension such that $A \times_K M$ has semistable reduction. We define the semistability defect of $A/K$ as $\sd_{A/K} \coloneqq [M : \unr{K}]$.
\end{definition}

\begin{remark}\label{rmk: Charact semist defect}
	Assume that $A/K$ has potential good reduction. We can characterise the field extension $M / \unr{K}$ from Definition~\ref{def: semist defect} (and thus the semistability defect of $A/K$) as follows:
	\begin{itemize}[leftmargin=*, itemsep=-1pt, topsep=5pt]
		\item $M$ is the subfield of $\sep{K}$ fixed by $\ker \left(\restr{\rholA}{I_{K}} \right)$.

		\item For any $m \geq 3$ coprime to the residue characteristic of $K$, we have $M = \unr{K}(A[m])$ (see \cite{SerreTate}).
	\end{itemize}
\end{remark}

\subsubsection{Conductors and associated \WD-representations}\label{sect: WD reps ab var RM}

Until the end of this subsection, we assume that $A/K$ has \RM \ by a totally real field $F$. Fix an embedding $i : \Ql \hookrightarrow \C$. In \eqref{eq: Decomp l adic lambda adic}, we decomposed $\rholA$ as the direct sum of $\Res_{\Flam / \Ql}(\rholamA)$, where $\lambda$ ranges through the places of $F$ dividing $\ell$, a prime different from the residue characteristic of $K$. Considering the associated complex \WD-representations, we get
\begin{equation}\label{eq: Decomp WD}
	\Wi(\rholA) \, \simeq \, \bigoplus_{\lambda \mid \ell}  \bigoplus_{\iota \, : \Flam \hookrightarrow \C} \W_{\iota}(\rholamA), 
\end{equation}
where the second direct sum on the RHS is indexed by embeddings of $\Flam$ into $\C$ that extend $i$. 

Knowing the reduction type of $A$, we can describe the inertial local type of $\W_{\iota}(\rholamA)$. The Néron-Ogg-Shafarevich criterion \cite[Theorem 1]{SerreTate} states that $A$ has good reduction at $\frakm$ if and only $\rholA(I_K)$ is trivial. Therefore, the main interesting case to study the inertial local type is when $A$ has bad reduction at $\frakm$. 

If $A/K$ has potential good reduction, let $M / \unr{K}$ be the extension cut out by $\restr{\rholA}{I_K}$, so that $\Gal(M / \unr{K}) \simeq \rholA(I_K)$  
and $[M : \unr{K}]$ is the semistability defect of $A/K$. 
Let also $I^{t}(M/K)$ be the tame part of $\Gal(M/\unr{K}) \simeq I(M / K)$.

\begin{proposition}\label{prop: local inert types RM}
	All the \WD-representations $\W_{\iota}(\rholamA)$ are $\Aut(\C)$-conjugate. 
	\begin{enumerate}[leftmargin=*, itemsep=0pt, topsep=5pt]
		\item If $A$ has potential toric reduction, then every $\W_{\iota}(\rholamA)$ is Steinberg.
		
		\item If $A$ has potential good reduction, then every $\W_{\iota}(\rholamA)$ is semisimple, and $\det(\W_{\iota}(\rholamA)) = \omega_K^{-1}$. Moreover, $\W_{\iota}(\rholamA)$ is a principal series $\Leftrightarrow M/ \unr{K}$ is abelian $ \Leftrightarrow |I^{t}(M/K)|$ divides $\abs{\Ffm} -1$.
	\end{enumerate}
\end{proposition}

\begin{proof}
	See Propositions 3.4 and 3.6, Theorem 4.1 and \S 5 in \cite{Melninkas}. 
\end{proof}

We now focus on the Artin conductor of $\rholamA$. The content of \cite[Appendix]{ACIKMM} shows that the conductor of $\W_{\iota}(\rholamA)$ does not depend on $\lambda$, and implies the equalities
\begin{equation}\label{eq: Conductor l adic lambda adic}
	\condtame{\rholA} = [F : \Q] \, \condtame{\rholamA}, \quad \text{ and } \quad \condwild{\rholA} = [F : \Q] \, \condwild{\rholamA}.
\end{equation}

We can compute the tame part of the conductor from the reduction type of the Néron model.

\begin{proposition}\label{prop: Tame cond lambda}
	The tame part of the conductor of $\rholamA$ is given by
	\begin{equation*}
		\condtame{\rholamA} = \begin{cases}
			0 & \text{ if } \Amod \text{ has good reduction at } \frakm, \\
			1 & \text{ if } \Amod \text{ has toric reduction at } \frakm, \\
			2 & \text{ if } \Amod \text{ has unipotent reduction at } \frakm.
		\end{cases}
	\end{equation*}
\end{proposition}

\begin{proof}
	Let $T$ and $U$ be the torus and the unipotent group as in \eqref{eq: ses Chevalley}. Grothendieck proved in \cite[Exposé IX, \S 4]{SGA7} that $\condtame{\rholA} = \dim T + 2 \dim U$. If $\Amod$ has toric (resp. unipotent) reduction at $\frakm$, then $\condtame{\rholA} = g$ (resp. $2g$). But we have $[F : \Q] = g$, so \eqref{eq: Conductor l adic lambda adic} gives the desired result.
\end{proof}

Concerning the wild part of the conductor, we have the following result.

\begin{lemma}\label{lem: Wild cond lambda}
	If $A$ attains semistable reduction over a tame extension of $K$, then $\condwild{\rholA} = 0$. If we further assume that $A/K$ has \RM, then $\condwild{\rholamA} = 0$.
\end{lemma}

\begin{proof}
	Let $L$ be a tame extension of $K$ where $A$ attains semistable reduction. By \cite[Exposé IX, \S 4]{SGA7}, we have $\condwild{\restr{\rholA}{G_L}} = 0$. The result follows by applying Lemma~\ref{lem: Behaviour nwild tame ext} and \eqref{eq: Conductor l adic lambda adic}.
\end{proof}

\subsection{Background on hyperelliptic curves}

We now recall some well-known facts about hyperelliptic curves and their models. Moreover, we briefly present the theory of cluster pictures, which allows to study the local arithmetic of hyperelliptic curves at odd places of bad reduction.

\subsubsection{Hyperelliptic curves and defining equations}\label{sect: Hyperell eqs}

For now, we let $L$ be any field. In this subsection we discuss the main properties of hyperelliptic curves. The reader can find more details on the topic in \cite{Liubook}, \cite{Lockhart}, \cite{Liu96}.

\begin{definition}
	A hyperelliptic curve $C/L$ is a smooth projective curve defined over $L$, such that there exists a finite separable morphism $C \rightarrow \Ppr_L$ of degree $2$.
\end{definition}

Let $C/L$ be a hyperelliptic curve, and write $L(C)$ for its function field. Let $\sigma$ be the generator of $\Gal(L(C) / L(\Ppr_L))$. It induces an automorphism of $C$ of order $2$, which we still denote by $\sigma$, and call the \textit{hyperelliptic involution} of $C$. If $C$ has genus $\geq 2$, then $\sigma$ is unique, and the degree $2$ morphism $C \rightarrow \Ppr_L$ is unique up to automorphism. From now on we assume that this is indeed the case, and denote by $g \geq 2$ the genus of $C$. We refer to \cite{Silverman86, Silverman94} for further details on the elliptic case $g=1$.

The curve $C$ can be described by a hyperelliptic equation. One can always choose two functions $x \in L(C)^{\langle \sigma \rangle}$ and $y \in L(C)$ such that $\lbrace 1, y \rbrace$ is a basis of the integral closure of $L[x]$ in $L(C)$ (see \cite[\S 1]{Liu96}). Then $C$ is described on an affine subset by a so-called \textit{Weierstrass equation}
\leqnomode
\begin{equation}\label{eq: Hhyp}\tag{$\Hhyp$}
	 y^2 + Q(x) y = P(x),
\end{equation}
where $Q(x), P(x) \in L[x]$, and $\deg Q \leq g+1, \deg P \leq 2g+2$. When $\charact(L) \neq 2$, one can choose the function $y$ in such a way that $Q=0$. When $Q = 0$ and the leading coefficient of $P$ is a non-zero square in $L$, then $C(L) \neq \emptyset$. Indeed $C$ has at least one or two $L$-points at infinity (\ie in the chart where $1/x \neq 0$), depending on whether $\deg P$ is odd or even. Following \cite{Liu96}, we introduce:

\begin{definition}
	Let $F(x) = 4P(x) + Q(x)^2$, and denote by $c_F$ the leading coefficient of $F$. When $\charact(L) \neq 2$, the discriminant $\Delta(\Hhyp)$ of equation \eqref{eq: Hhyp} is defined by 
	\begin{equation*}
		\Delta(\Hhyp) = \begin{cases}
			2^{-4(g+1)} \disc(F) \quad & \text{ if } \deg F = 2g+2, \\
			2^{-4(g+1)} c_F^2 \disc(F) \quad & \text{ if } \deg F = 2g+1.			
		\end{cases} 
	\end{equation*}
\end{definition}

In particular, when $Q=0$ and $P$ is monic, we have $\Delta(\Hhyp) = 2^{4g} \disc(P)$, whether $\deg P$ is odd or not. The curve described by \eqref{eq: Hhyp} is smooth if and only if $\Delta(\Hhyp) \neq 0$ (see \cite[Theorem 1.7]{Lockhart}).
\medskip

\begin{definition}
	Given a hyperelliptic equation $(\Hhyp) : y^2 + Q(x) y = P(x)$ describing a curve $C/L$, we define the Weierstrass model of $C$ associated to $(\Hhyp)$ as the glueing of the two open affine schemes
	\begin{equation*}
		\Spec \left(L[x, y] / (y^2 + Q(x) y = P(x))\right) \quad \text{ and } \quad \Spec \left(L[w, z] / (z^2 + w^{g+1}Q(w^{-1}) w = w^{2g+2}P(w^{-1}))\right)
	\end{equation*}
	along the identification $w = 1/x, z = y/x^{g+1}$. We define its discriminant as $\Delta(\Hhyp)$. If we let $\Wrst$ be the model associated to $(\Hhyp)$, we will write $\Wrst : y^2 + Q(x) y = P(x)$ to state that $\Wrst$ is described by $(\Hhyp)$.
\end{definition}
\medskip

A Weierstrass equation defining the hyperelliptic curve $C$ is not unique, so neither is a Weierstrass model of $C$. We now describe the admissible changes of variables to obtain a new hyperelliptic equation for $C$, and the behaviour of the discriminant under such transformations.

\begin{lemma}\label{lem: Legal change of vars}
	Let $(\Hhyp) \, : \, y^2 + Q(x) y = P(x)$ and $(\widetilde{\Hhyp}) \, : \, Y^2 + \widetilde{Q}(X)Y = \widetilde{P}(X)$ be two hyperelliptic equations describing the curve $C / L$. The coordinates $(x, y)$ and $(X,Y)$ are related by a change of variables
	\begin{equation*}
		x = \frac{a X + b}{cX + d}, \quad \text{ and } \quad y = \frac{e Y + R(X)}{(cX+ d)^{g+1}},
	\end{equation*}
	with $a, b, c, d, e \in L$, $R(X) \in L[X]$, and $ad - bc, e \neq 0$. We have the equality between discriminants
	\begin{equation*}
		\Delta(\widetilde{\Hhyp}) = e^{-4(2g+1)} (ad-bc)^{2(g+1)(2g+1)} \Delta(\Hhyp).
	\end{equation*}
\end{lemma}

\begin{proof}
	See \cite[\S 1]{Lockhart} or \cite[\S 2]{Liu96}.
\end{proof}

\begin{definition}\label{def: Twist hyper crv}
	Let $C$ be a hyperelliptic curve described by $(\Hhyp) : y^2 = P(x)$. Let $\delta \in L$. We define the quadratic twist of $C$ by $\delta$, which we denote by $C^{(\delta)}$, to be the hyperelliptic curve described by the affine equation $y^2 = \delta P(x)$. When $\deg P$ is odd, we can also describe $C^{(\delta)}$ by 
	\reqnomode
	\begin{equation}\label{eq: quadratic twist}
		Y^2 = \delta^{\deg P} P \left( X/ \delta \right), \qquad \text{ where } \qquad x = \frac{X}{\delta}, \ y = \frac{Y}{\delta^{(\deg P - 1)/2}}.
	\end{equation}
	The curves $C$ and $C^{(\delta)}$ are isomorphic over $L(\sqrt{\delta})$, so if $\delta$ is a square in $L$, then $C \simeq C^{(\delta)}$ over $L$. 
\end{definition}

\subsubsection{Local models and cluster pictures}\label{sect: Models and clusters}

We now discuss the behaviour of hyperelliptic curves defined over local fields. That is why, we recover the previous notation from \S \ref{sect: Notation} and assume that $L = K$ is local. We keep denoting by $C / K$ a hyperelliptic curve. 

\begin{definition}
	$1)$ A model of $C/K$ over $\Om$ is a flat proper $\Om$-scheme together with a $K$-isomorphism of its generic fiber with $C$.
	
	$2)$ We say that the Weierstrass model $\Wrst$ of $C$ corresponding to $(\Hhyp) : y^2 + Q(x) y = P(x)$ is $\Om$-integral if $P, Q \in \Om[x]$. In this case, $\Wrst \rightarrow \Spec(\Om)$ is a model in the sense above. We say that $\Wrst$ is a minimal Weierstrass model if $\vm(\Delta(\Hhyp))$ is minimal among all integral Weierstrass models of $C$.
\end{definition}

\begin{definition}
	Let $\Xmod \rightarrow \Spec(\Om)$ be a model of $C$ over $\Om$, and let $\Xmod_{\frakm}$ be its special fiber.
	\begin{enumerate}[leftmargin=*, topsep = 5pt, itemsep=0pt]
		\item We say that $\Xmod$ has good reduction at $\frakm$ if $\Xmod_{\frakm}$ is smooth.
		
		\item We say that $\Xmod$ has semistable reduction at $\frakm$ if $\Xmod_{\frakm}$ is geometrically reduced and has only ordinary double points as singularities.
	\end{enumerate}	
\end{definition}

We say that $C/K$ has good (resp. semistable) reduction at $\frakm$ if there is some model of $C$ over $\Om$ having good (resp. semistable) reduction at $\frakm$. Deligne and Mumford proved in \cite{DeligneMumford} that $C / K$ has semistable reduction if and only if its Jacobian $\Jac(C) / K$ has semistable reduction.

\begin{remark}\label{rmk: Bad red congr}
	Suppose that $C / K$ is described by the equation $(\Hhyp) : y^2 + Q(x) y = P(x)$. Applying a change of variables as in Lemma~\ref{lem: Legal change of vars} modifies the valuation of $\Delta(\Hhyp)$ by adding a multiple of $2(2g+1)$. Therefore, if $\vm(\Delta(\Hhyp)) \not \equiv 0 \mod 2(2g+1)$, then $C$ does not have good reduction over $\Om$.
\end{remark}

Denote by $\Cmod \rightarrow \Spec(\O_K)$ the \textit{minimal regular model} of $C$, and by $\Cmodm$ its special fiber (see \cite[\S 8.3]{Liubook} for details about its construction and main properties). If $C(K) \neq \emptyset$, then $\Cmod$ is the minimal desingularization of a minimal Weierstrass model of $C$. In practice, if $\Wrst$ is a minimal Weierstrass model, we construct $\Cmod$ by successively blowing-up $\Wrst$ along its singular points and normalising a minimal amount of times until we obtain a regular model that contains no exceptional divisors. Assume that the $\lcm$ of the multiplicities of the irreducible components of the special fiber is $1$ (\textit{e.g.} $C(K) \neq \emptyset$ or $\Cmodm$ is reduced). Then, the special fiber of the Néron model of the Jacobian $\Jac(C)$ is $\Pic^{0}(\Cmodm)$, so the reduction type of the Néron model of $\Jac(C)$ is encoded in the geometry of $\Cmodm$. 
\medskip

The following criterion is very helpful for checking if a scheme is regular or not.

\begin{lemma}\label{lem: Regularity criterion}
	Let $(A, \frakm)$ be a regular Noetherian local ring, and let $f \in \frakm \setminus \lbrace 0 \rbrace$. Then $A / fA$ is regular if and only if $f \notin \frakm^2$.
\end{lemma}

\begin{proof}
	This is Corollary $2.12$ in \cite[\S 4.2.2]{Liubook}.
\end{proof}
\medskip

From now on, we assume that the residue characteristic of $K$ is odd. The theory of cluster pictures provides a combinatorial way to describe the special fiber of the minimal regular model, and to deduce arithmetic information concerning $C$ and $\Jac(C)$. We refer to \cite{M2D2} for further details on the theory, or to \cite{Hyperuser} for a survey. Assume that $C/K$ is described by the Weierstrass equation
\begin{equation*}
	y^2 = P(x) = c_P \prod_{\gamma \in \mathcal{R}}(x-\gamma) \vspace{-0.5em}
\end{equation*}
where $c_P \in K\st$, $P(x) \in K[x]$ is separable, and $\mathcal{R}$ denotes the set of roots of $P$ in $\sep{K}$. 

\begin{definition}
	A \textit{cluster} is a non-empty subset $\s \subseteq \mathcal{R}$ of the form $\s = D \, \cap \, \mathcal{R}$ for some disc $D= \lbrace x \in  \sep{K} \mid \vm(x-z) \geq d \rbrace$, some $z\in \sep{K}$ and some $d\in \Q$. If $|\s|>1$, we say that $\s$ is a \textit{proper} cluster, and define its \textit{depth} $d_\s$ as the maximal $d$ for which $\s$ is cut out by such a disc, \ie $d_\s  = \min_{\gamma,\gamma' \in \mathfrak{s}} \vm(\gamma - \gamma')$. If $\s\neq \mathcal{R}$, the \textit{parent} cluster $P(\s)$ of $\s$ is the smallest cluster with $\s\subsetneq P(\s)$, and the \textit{relative depth} of $\s$ is $\delta_\s\! =\! d_\s\! -\!d_{P(\s)}$.
	
	We refer to this data as the \textit{cluster picture} of $C$.
\end{definition}

We introduce some terminology concerning cluster pictures (see \cite{Hyperuser} for further details).

\begin{definition}
	\begin{itemize}[leftmargin=*, itemsep=0pt]
		\item If $\s'\subsetneq \s$ is a maximal subcluster, we say that $\s'$ is a \textit{child} of $\s$. For two clusters $\s_1$, $\s_2$, we write $\s_1\wedge \s_2$ for the smallest cluster containing both of them.
		
		\item A cluster $\s$ such that $|\s|$ is odd (resp. even) is called an \textit{odd} (resp. \textit{even}) cluster. 
		
		\item A cluster consisting of precisely two roots is called a \textit{twin}. 
		
		\item  An even cluster whose children are all even is called \textit{\"ubereven}. 
		
		\item A non-\"ubereven cluster with a child of size $2g$ is called a \textit{cotwin}.
		
		\item A cluster $\s$ is called principal if one of the following holds \vspace{-0.5em}
		\begin{equation*}
			\begin{cases}
				|\s|\neq 2g+2 \text{ and } \s \text{ is proper, not a twin or a cotwin;}\\
				|\s|=2g+2 \text{ and } \s \text{ has more than }3 \text{ children.}
			\end{cases}\vspace{-0.5em}
		\end{equation*}
	
		\item Recall that $c_P$ denotes the leading coefficient of $P$. We define the quantity
		\begin{equation*}
			\nu_{\s}= \vm(c_P)+|\s |d_{\s}+ \sum_{\gamma \notin \s} \, 	d_{ \{\gamma\} \wedge \, \s}.
		\end{equation*}
	\end{itemize}
\end{definition}

Cluster pictures provide criteria to check if $C/K$ has good, or semistable reduction. The following two theorems are (part of) \cite[Theorem 10.8]{M2D2}.

\begin{theorem}\label{thm: Crit clusters good}
	The hyperelliptic curve $C/K$ has good reduction if and only if the following three conditions are satisfied: \vspace{-0.25em}
	\begin{enumerate}[itemsep=-1pt, topsep=5pt]
		\item the field extension $K(\Rroots) / K$ is unramified,
		\item every proper cluster has size at least $2g+1$,
		\item the (necessarily unique) principal cluster has $\nu_{\s} \in 2 \Z$.
	\end{enumerate}
\end{theorem}

\begin{theorem}\label{thm: Crit clusters semist}
	The hyperelliptic curve $C/K$ (equivalently $\Jac(C)$) has semistable reduction if and only if the following three conditions are satisfied: \vspace{-0.25em}
	\begin{enumerate}[itemsep=-1pt, topsep=5pt]
		\item the field extension $K(\Rroots) / K$ has ramification index at most $2$,
		\item every proper cluster is invariant under the action of the inertia group $I_K$,
		\item every principal cluster $\s$ has $d_{\s} \in \Z$ and $\nu_{\s} \in 2\Z$.
	\end{enumerate}
\end{theorem}

Moreover, \cite[Theorem 11.3]{M2D2}, which we recall below, gives the wild conductor of the $\ell$-adic representation arising from the action of $G_K$ on $\Vl (\Jac(C))$. Recall that $g$ denotes the genus of $C$.

\begin{theorem}\label{thm: Wild cond clusters}
	The wild conductor of $\rho_{\Jac(C), \, \ell} : G_K \rightarrow \GL_{2g}(\Ql)$ is given by 
	\begin{equation*}
		\condwild{\rho_{\Jac(C), \, \ell}} = \sum_{\gamma \in \Rroots / G_K} \left( \vm(\Delta(K(\gamma) / K)) - [K(\gamma) : K] + f_{K(\gamma) / K}\right).
	\end{equation*}
	The sum is taken over representatives of orbits in $\Rroots$ under the action of $G_K$, and $\Delta(K(\gamma) / K)$ (resp. $f_{K(\gamma) / K}$) denotes the discriminant (resp. residue degree) of the extension $K(\gamma) / K$.
\end{theorem}

\begin{remark}
	Recently, there has been ongoing work on the use of cluster pictures in the case of even residue characteristic. Some partial results have been obtained, namely in the study of the stable reduction of hyperelliptic curves. Let us mention the work of Dokchitser--Morgan \cite{DokchitserMorgan23}, and Gehrunger \cite{Gehrunger25}, which we will use later in section \ref{sect: Reduction types} to describe the reduction types of Jacobians of Frey hyperelliptic curves at even places. 
\end{remark}

% -- % -- % -- % -- % -- % -- % -- % -- % -- % -- % -- % -- % -- % -- % -- % -- % -- % -- % -- % -- % -- % -- % -- % -- % -- % -- % -- % -- % -- % -- % -- % -- % -- % -- % -- % -- % -- % -- % -- 

\section{Construction of Frey representations and curves}\label{sect: Constr Frey reps crvs}

We now recall the definition of Frey representations, and their construction for the specific cases of signatures $(p,p,r)$ and $(r,r,p)$. We will then obtain by specialisation Frey curves attached to putative solutions to the generalised Fermat equations of the mentioned signatures. We will conclude the section by presenting a common framework to deal with both signatures at once.

\subsection{Geometric construction of Frey representations}

In this subsection, we recall the construction of Frey representations of signature $(p, p, r)$ and $(r, r, p)$. Our main references are \cite{Darmon00} for the former signature, and \cite{BCDF23} for the latter.
\medskip 

Recall that we introduced above two primes $r \geq 5$, $p$, and that we denote by $\Kgl = \Q(\omega)$ the maximal totally real subfield of $\Q(\zr)$. Denote by $\Ppr$ the projective line over $\Kglbar$, with local coordinate $s$ and function field $\Kglbar(s)$. Let $x \in \Ppr(\C)$, $x \neq 0, 1, \infty$, and consider the topological fundamental group $\pi_1 \left( \Ppr(\C) \setminus \lbrace 0, 1, \infty \rbrace, x\right)$.  Its profinite completion equals the Galois group $\Gal(\Omega / \Kglbar(s))$, where $\Omega \subset \overline{\Kgl(s)}$ is the maximal extension of $\overline{\Kgl}(s)$ unramified outside the points $0, 1, \infty$ (see \cite[\S 6.3]{Serre08}). Thus, any continuous representation of $\pi_1 \left( \Ppr(\C) \setminus \lbrace 0, 1, \infty \rbrace, x\right)$, usually called a \textit{monodromy} representation, extends to $\Gal(\Omega / \overline{\Kgl(s)})$, and thus to $G_{\Kglbar(s)}$. Since the latter is a normal subgroup of $G_{\Kgl(s)}$, we obtain in this way continuous representations of $G_{\Kgl(s)}$. We refer to \cite{GolfieriPacetti} for broader use of monodromy representations in the context of Darmon's program.

\begin{definition}\label{def: Frey rep}
	A Frey representation of signature $(p, p, r)$ (resp. $(r, r, p)$) is a Galois representation $\varrho : G_{\Kgl(s)} \longrightarrow \GL_2(\Ff)$, where $\Ff$ is a finite field, satisfying the following:
	\begin{enumerate}[itemsep=0pt]
		\item The restriction of $\varrho$ to $G_{\Kglbar(s)}$ has trivial determinant and is irreducible. Denote its projectivization by $\overline{\varrho}^{\operatorname{geom}} : G_{\Kglbar(s)} \rightarrow \operatorname{PSL}_2(\Ff)$.
		\item  The homomorphism $\overline{\varrho}^{\operatorname{geom}}$ is unramified outside $\lbrace 0, 1, \infty \rbrace$.
		\item $\overline{\varrho}^{\operatorname{geom}}$ maps the inertia subgroups at $0, 1, \infty$ to subgroups of $\operatorname{PSL}_2(\Ff)$ of order $p, p, r$ (resp. $r, r, p$).
	\end{enumerate}
\end{definition}

\subsubsection{Frey representations of signature $(p,p,r)$}\label{sect: Frey reps ppr}

We now proceed to recall the geometric construction of Frey representations of signatures $(p, p, r)$. We refer to \cite[\S 1.3]{Darmon00} for details. Recall that $h_r \in \Z[x]$ denotes the minimal polynomial of $\omega = \zr + \zr^{-1}$.

\begin{definition}\label{def: Curve Cpm(s)}
	We define $\Cms(s), \Cpl(s) / \Q(s)$ to be the hyperelliptic curves given by the equations 
	\begin{equation*}
		(\Hhyp_{r, s}^{-}) \, : \, y^2 = \gminuss(x) \qquad \text{ and } \qquad (\Hhyp_{r, s}^{+}) \, : \, y^2 = \gpluss(x)
	\end{equation*}
	respectively, where $\gminuss(x), \, \gpluss(x) \in \Z[x]$ are the polynomials defined by \vspace{-0.25em}
	\reqnomode
	\begin{equation}\label{eq: Def gpms}
		\gminuss(x) \coloneqq (-1)^{\frac{r-1}{2}} x \, h_r(2 - x^2) + 2 - 4s, \qquad \text{ and } \qquad \gpluss(x) \coloneqq \gminuss(x)(x+2).
	\end{equation}
	We also define $\Jac(\Cpm(s)) / \Q(s)$ to be the Jacobian of $\Cpm(s)$, and let $J_{r}^{\pm}(s) \coloneqq \Jac(\Cpm(s)) \times_{\Q(s)} \Kgl(s)$.
\end{definition}

\begin{remark}
	We point out that the notation above slightly differs from the one used by Darmon. In \cite[\S 1.3]{Darmon00}, he defines $\Cms(s)$ using the minimal polynomial $- \omega$, which he denotes by $g$. The polynomials $h_r$ and $g$ are related by the equality $h_r(x) = (-1)^{\frac{r-1}{2}} g(-x)$.
\end{remark}

The polynomial $\gminuss$ has degree $r$, and $\gpluss$ has degree $r+1$, so both curves $\Cpm(s)$ have genus $\frac{r-1}{2}$. The polynomial $\gminuss$ can be expanded: using \cite[Lemma 2.16]{BCDF23}, we rewrite \eqref{eq: expanded Cms(s)} as \vspace{-0.5em}
\leqnomode
\begin{equation}\label{eq: expanded Cms(s)}\tag{$\Hhyp_{r, s}^-$}
	y^2 = \sum_{k=0}^{\frac{r-1}{2}} (-1)^k c_k \, x^{r-2k} + 2 - 4s, \quad \text{ where } \quad  c_k \coloneqq \frac{r}{r-k} \binom{r-k}{k} \in \Z.
\end{equation}
Given the equality $\gpluss(x) = (x+2) \gminuss(x)$, we easily deduce an expanded equation defining $\Cpl(s)$.

\begin{example}
	For small values of $r$, the hyperelliptic curve $\Cms(s)$ is described by
	\begin{align}
		y^2 & = x^3 - 3 x + 2 - 4s, \tag{$\Hhyp_{3, s}^-$}\\
		y^2 & = x^5 - 5 x^3 + 5  x + 2 - 4s, \tag{$\Hhyp_{5, s}^-$}\\
		y^2 & = x^7 - 7 x^5 + 14 x^3 - 7x + 2 - 4s, \tag{$\Hhyp_{7, s}^-$}\\
		y^2 & = x^{11} - 11 x^9 + 44 x^7 - 77 x^5 + 55 x^3 - 11 x + 2 - 4s. \tag{$\Hhyp_{11, s}^-$}
	\end{align}
\end{example}

\begin{theorem}\label{thm: Jpm have RM}
	The abelian varieties $J_{r}^{\pm}(s) / \Kgl(s)$ have \RM \ by $\Kgl$, and more precisely
	\begin{equation*}
		\End_{\Kgl(s)} \left(J_{r}^{\pm}(s)\right) \simeq \OK.
	\end{equation*}
\end{theorem}

\begin{proof}
	This follows from Theorem $1$, Corollary $6$ in \cite{TautzTopVerberkmoes}. See also \cite[Proposition 2.1]{DarmonMestre}.
\end{proof}

Let now $\frak p$ be a place of $\Kgl$ lying above $p$, $\Kp$ be the completion of $\Kgl$ at $\frak p$, and $\Ffp$ its residue field. The discussion in \S \ref{sect: Decomp Tate module} applies to $J_{r}^{\pm}(s)$, giving rise to $2$-dimensional representations
\begin{equation*}
	\rho_{J_{r}^{\pm}(s), \, \frak p} : G_{\Kgl(s)} \longrightarrow \GL_{2}(\Kp) \quad \text{ and } \quad \overline{\rho}_{J_{r}^{\pm}(s), \, \frak p} : G_{\Kgl(s)} \longrightarrow \GL_{2}(\Ffp)
\end{equation*}

\begin{theorem}\label{thm: Frey rep ppr}
	Both $\overline{\rho}_{J_{r}^{+}(s), \, \frak p}$ and $\overline{\rho}_{J_{r}^{-}(s), \, \frak p}$ are Frey representations of signature $(p, p, r)$. The former is even, and the latter is odd in the sense of \cite[Definition 1.4]{Darmon00}. 
\end{theorem} 

\begin{proof}
	See \cite[Theorem 1.10]{Darmon00}.
\end{proof}

\subsubsection{Frey representations of signature $(r, r, p)$}\label{sect: Frey reps rrp}

We now introduce a Frey representation of signature $(r, r, p)$, building upon the construction of the curve $\Cms(s)$ introduced above. In \cite{Pacetti25}, Pacetti explains how to recover this construction using the theory of hypergeometric motives. We will not discuss these topics here, and we refer the curious reader to the mentioned reference. The discussion below is inspired by the results in \cite[\S 2.4]{BCDF23}. 
\medskip

In order to simplify the notation, let $\stilde \coloneqq 2 - 4 s$ be a new indeterminate in $\Kglbar(s)$. With this notation, the constant term of \eqref{eq: expanded Cms(s)} is simply $\stilde$. Consider two new copies of the projective line over $\Kglbar$, which we denote by $\Ppr_{\stilde}$ and $\Ppr_t$, with local coordinates $\stilde$ and $t$. The definition of $\stilde$ implies an equality of function fields $\Kglbar(s) = \Kglbar(\stilde)$. Let $\B / \Kglbar$ be the smooth projective curve defined by the affine equation
\reqnomode
\begin{equation}\label{eq: Def curve B}
	t(t-1) = \frac{1}{\stilde^2 - 4}, \quad \text{ or equivalently,} \quad t(t-1) = \frac{1}{2^4s(s-1)}.
\end{equation}
Denote by $\Kglbar(\B)$ its function field. The projection onto the $s$ and $t$ coordinates give two morphisms $\B \rightarrow \Ppr_{s}$ and $\B \rightarrow \Ppr_t$ of degree two, having $0$ and $1$ as branching points. Therefore, $\Kglbar(\B)$ is a quadratic extension of both function fields $\Kglbar(\stilde) = \Kglbar(s)$ and $\Kglbar(t)$. Equality \eqref{eq: Def curve B} shows that the points $0, 1 \in \Ppr_s$ (or equivalently, $\pm 2 \in \Ppr_{\stilde}$) correspond to $\infty \in \Ppr_t$, and similarly $0, 1 \in \Ppr_t$ correspond to $\infty \in \Ppr_s$ (and thus to $\infty \in \Ppr_{\stilde}$).
\medskip

\begin{remark}
	The authors of \cite{BCDF23} already introduce the function field $\Kglbar(\B)$. However, their notation is different, as they introduce the curve $\B' / \Kglbar$ described by $t(1-t) = 1/(\sigma^2 + 4)$. This difference between the considered models is due to the fact that the authors of \cite{BCDF23} do not manipulate $\Cms(s)$, but its quadratic twist by $i = \sqrt{-1}$. If we were to consider $\B$ and $\B'$ as curves over $\Kgl$ (instead of $\Kglbar$) defined by the same equations, their function fields would not be isomorphic.
\end{remark}

Define $\deltatilde \coloneqq (2t-1)/ \stilde = (2t-1) / (2 - 4s) \in \Kglbar(\B)$. Using \eqref{eq: Def curve B}, one can check that $\stilde^2 = \frac{(2t-1)^2}{t(t-1)}$, so $\deltatilde^2 = t(t-1)$. Consider $\Cms(s)^{(\deltatilde)}$, the quadratic twist by $\deltatilde$ of the base change of $\Cms(s)$ to $\Kgl(\B)$. Thanks to \eqref{eq: expanded Cms(s)}, we describe $\Cms(s)^{(\deltatilde)}$ by the hyperelliptic equation\vspace{-0.5em}
\leqnomode
\begin{equation*}\label{eq: def Crt twist}\tag{$\Hhyp_{r, \, t}$}
	y^2 = \sum_{k = 0}^{\frac{r-1}{2}} (-1)^k c_k \, (t(t-1))^k x^{r-2k} + (t(t-1))^{\frac{r-1}{2}} (2t-1).
\end{equation*}

\textit{A priori}, $\Cms(s)^{(\deltatilde)}$ is defined over $\Kglbar(\B)$, but \eqref{eq: def Crt twist} defines a model of $\Cms(s)^{(\deltatilde)}$ over $\Q(t)$.

\begin{definition}
	We let $C_r(t) / \Q(t)$ be the hyperelliptic curve given by the Weierstrass equation \eqref{eq: def Crt twist}.	We also define $\Jac(C_r(t)) / \Q(t)$ to be the Jacobian of $C_r(t)$, and let $J_{r}(t) \coloneqq \Jac(C_r(t)) \times_{\Q(t)} \Kgl(t)$. 
\end{definition}

\begin{theorem}
	The abelian variety $J_{r}(t) / \Kgl(t)$ has \RM \ by $\Kgl$.
\end{theorem}

\begin{proof}
	This is \cite[Theorem 2.38]{BCDF23}: although it only states that $\Kgl \hookrightarrow \End_{\Kgl(t)}(J_r(t)) \otimes \Q$, the arguments in its proof show that the inclusion is actually an isomorphism. 
\end{proof}

Recall that $\frakp$ is a place of $\Kgl$ lying above $p$. Again, the discussion in section~\ref{sect: AV with RM} applies to $J_{r}(t)$, giving rise to $2$-dimensional representations
\begin{equation*}
	\rho_{J_{r}(t), \, \frak p} : G_{\Kgl(t)} \longrightarrow \GL_{2}(\Kp) \quad \text{ and } \quad \overline{\rho}_{J_{r}(t), \, \frak p} : G_{\Kgl(t)} \longrightarrow \GL_{2}(\Ffp)
\end{equation*}

\begin{theorem}
	The representation $\overline{\rho}_{J_{r}(t), \, \frak p}$ is a Frey representation of signature $(r, r, p)$. 
\end{theorem} 

\begin{proof}
	This is a consequence of Theorem~\ref{thm: Frey rep ppr}, and the fact that $0, 1 \in \Ppr_s$ (resp. $\infty \in \Ppr_s$) correspond to $\infty \in \Ppr_t$ (resp. $0, 1 \in \Ppr_t$). See Theorem 2.38 and Lemma 2.3 in \cite{BCDF23} for details.
\end{proof}

% -- % -- % -- % -- % -- % -- % -- % -- % -- % -- % -- % -- % -- % -- % -- % -- % -- % -- % -- % -- % -- % -- % -- % -- % -- % -- % -- % -- % -- % -- % -- % -- % -- % -- % -- % -- % -- % -- % -- 

\subsection{Assumptions on Diophantine equations and their solutions}

We now discuss the main assumptions that we will be doing on the parameters of the generalised Fermat equations under consideration, and on their putative solutions.

\begin{definition}\label{def: Parameters GFE}
	We let $r \geq 5$ be a fixed prime number, and $A, B, C \in \Z \setminus \lbrace 0 \rbrace$ be three pairwise coprime integers, that are free of $r$-th powers. We denote by $p$ any prime number greater than $r$.
\end{definition}

The generalised Fermat equations of signatures $(p, p, 3)$ and $(3, 3, p)$ have already been studied in the literature. We refer the curious reader to \cite{BennettVatsalYazdani}, where Bennett, Vatsal and Yazdani study the former equation; and to \cite{BennettDahmen}, where Bennett and Dahmen study the latter. Moreover, the case $r = p$ was also studied by Kraus in \cite{Kraus97}, and by Dieulefait and Soto in \cite{DieulefaitSoto21}.

Restricting to the case $p > r$ is not problematic, as this is the generic case. The assumptions on $A, B, C$ are natural to consider for proving that \eqref{eq: GFE ppr} and \eqref{eq: GFE rrp} do not have non-trivial primitive solutions. For the GFE \eqref{eq: GFE ppr} : $A x^p + B y^p = Cz^r$, one can perfectly assume that $C$ is free of $r$-th powers, and that $A, B$ are free of $p$-th powers. Similarly, for \eqref{eq: GFE rrp} : $Ax^r + By^r = Cz^p$, one can assume that $A$ and $B$ are free of $r$-th powers, and $C$ free of $p$-th powers. If, for any of the signatures, this was not satisfied, one could replace the coefficients by some others being free of $r$-th and $p$-th powers. Proving the non-existence of solutions for the newly obtained equation implies the non-existence of solutions for the original one. 

Assuming that all the coefficients are free of $r$-th powers implies a loss of generality. Nevertheless, it allows for a uniform treatment of the families of equations \eqref{eq: GFE ppr}$_p$ and \eqref{eq: GFE rrp}$_p$, as the coefficients are necessarily independent of $p$. Moreover, under this (restrictive) hypothesis we have:

\begin{lemma}\label{lem: Aa Bb Cc coprime}
	Let $r, p, A, B, C$ be as in Definition~\ref{def: Parameters GFE}. 
	\begin{enumerate}[leftmargin=*, topsep=5pt, itemsep=0pt]
		\item If $(a, b, c)$ is a primitive non-trivial solution to \eqref{eq: GFE ppr}, then $Aa^p, Bb^p$ and $Cc^r$ are pairwise coprime.
		
		\item If $(a, b, c)$ is a primitive non-trivial solution to \eqref{eq: GFE rrp}, then $Aa^r, Bb^r$ and $Cc^p$ are pairwise coprime.
	\end{enumerate}	
\end{lemma}

\begin{proof}
	We prove the first statement, as the second one is treated in the exact same way. Assume that $Aa^p + Bb^p = Cc^r$. If any prime number divides two terms among $Aa^p, Bb^p$ and $Cc^r$, then it necessarily divides the third one too. Assume by contradiction that such a prime $q$ exists. Two different cases appear.
	\begin{itemize}[topsep=5pt, itemsep=0pt] 
		\item If $q \mid A$, then $q \nmid B, C$, as $A, B, C$ are pairwise coprime, so $q \mid b, c$. We assumed that $A$ is free of $r$-th powers, so $v_q(Aa^p) = v_q(A) < r < p \leq v_q(Bb^p)$.	We have 
		\begin{equation*}
			r \leq r v_q(c) = v_q(Cc^r) = v_q(Aa^p + Bb^p) = \min(v_q(Aa^p), v_q(Bb^p)) = v_q(Aa^p) < r,
		\end{equation*}
		hence a contradiction. The case $q \mid B$ is treated in the same way.
		
		\item If $q \nmid AB$, then $q \mid a^p, b^p$, and since $a, b, c$ are coprime, then $q \nmid c$. In this case we have
		\begin{equation*}
			v_q(C) = v_q(Cc^r) = v_q(Aa^p + Bb^p) \geq \min (v_q(Aa^p), v_q(Bb^p)) = \min (p v_q(a), p v_q(b)) \geq p > r.
		\end{equation*}
		This contradicts the fact that $C$ is free of $r$-th powers. The desired result follows.
	\end{itemize}\vspace{-1.75em}
\end{proof}

For the signature $(r, r, p)$, the following statement gives $r$-adic constraints on the solutions. %it might happen that the considered GFE has no solutions for elementary arithmetic reasons. The following statement gives one of such conditions. 

\begin{lemma}\label{lem: r2 divides}
	Let $r, A, B, a, b$ be as in Definition~\ref{def: Parameters GFE}. If $Aa^r + Bb^r \equiv 0 \mod r$, then
	\begin{equation*}
		A^{r-1} \equiv B^{r-1} \pmod{r^2} \ \Leftrightarrow \  v_{r}(Aa^r + Bb^r) \geq 2.
	\end{equation*}
\end{lemma}

\begin{proof}
	Assume first $A^{r-1} \equiv B^{r-1} \pmod{r^2}$. Since $r \mid Aa^r +  Bb^r$, we know that $r \nmid ABab$ by coprimality. Moreover, $Aa^r \equiv -   Bb^r \pmod{r}$, so by Fermat's little theorem we have $Aa \equiv -  Bb \pmod{r}$. Hence, there exists some $u \in \Z$ such that $Aa = -   Bb + ru$, and applying Newton's binomial formula, we obtain
	\begin{equation*}
		(Aa)^r = (-   Bb + ru)^r = -   (Bb)^r + \sum_{j = 1}^{r} \binom{r}{j} (-Bb)^{j} (ru)^{r-j} \equiv -   (Bb)^r \pmod{r^2}.
	\end{equation*}
	The first equality uses the fact that $r$ is odd, and the congruence uses the fact that $r \mid \binom{r}{r-1}$. Therefore we have $A^{r-1} A a^r \equiv -   B^{r-1} B b^r \pmod{r^2}$, and since $A, B \not \equiv 0 \pmod{r^2}$, we can simplify by $A^{r-1} \equiv B^{r-1} \pmod{r^2}$ to obtain $Aa^r \equiv -   Bb^r \pmod{r^2}$.
	
	Assume now that $r^2 \mid Aa^r +   Bb^r$, then $Aa^r \equiv -   Bb^r \pmod{r^2}$. The multiplicative group $\left( \Z /r^2 \Z \right)\st$ has order $r(r-1)$. Since $r \nmid ABab$ by coprimality, then raising both sides of the previous congruence to the $(r-1)$-th power gives $A^{r-1} \equiv B^{r-1} \pmod{r^2}$.
\end{proof}

It follows from the Lemma above that certain equations do not have non-trivial primitive solutions for elementary arguments. When solving specific families of equations in section~\ref{sect: Solving instances GFEs}, we will restrict ourselves to choices of coefficients that avoid this straightforward arguments for proving the non-existence of solutions. %is not straightforward to prove. We will give further details about this in section~\ref{sect: Solving instances GFEs}.

% -- % -- % -- % -- % -- % -- % -- % -- % -- % -- % -- % -- % -- % -- % -- % -- % -- % -- % -- % -- % -- % -- % -- % -- % -- % -- % -- % -- % -- % -- % -- % -- % -- % -- % -- % -- % -- % -- % -- 

\subsection{From putative solutions to Frey curves}

We now discuss how to associate Frey curves to putative solutions to the generalised Fermat equations \eqref{eq: GFE ppr} and \eqref{eq: GFE rrp}. Those are obtained as specialisations of the curves giving rise to the corresponding Frey representations. We first construct such curves for the signature $(p,p,r)$, and then for the signature $(r,r,p)$. Both constructions fit into a single framework, which will allow for a uniform study of both signatures.
\medskip

Recall that we introduced $r, p, A, B, C$ in Definition~\ref{def: Parameters GFE}.

\subsubsection{Frey curves $\Cpmabc$ for the signature $(p, p, r)$}\label{sect: Frey curves ppr}

The generalised Fermat equation \eqref{eq: GFE ppr bis} of signature $(p,p,r)$ is 
\begin{equation}\label{eq: GFE ppr bis}
	\tag{$\Eppr$}	 Ax^p + By^p = C z^r.
\end{equation}
Assume that there exists a primitive non-trivial solution $(a, b ,c) \in \Z^3$ to \eqref{eq: GFE ppr}. Following \cite{Darmon00}, we construct two Frey curves associated to $(a, b, c)$ as follows. Fix $s_0 \coloneqq Aa^p/Cc^r$, and let $\Cpm(Aa^p/Cc^r)$ be the specialisation at $s = s_0$ of the curve $\Cpm(s) / \Q(s)$ introduced in \S \ref{sect: Frey reps ppr}.

\begin{definition}
	We define $\Cmsabc / \Q$ to be the quadratic twist by $Cc$ of $\Cms(Aa^p / Cc^r)$, and $\Cplabc \coloneqq \Cpl(Aa^p/Cc^r)$. We call these the \textit{Frey curves} of signature $(p,p,r)$ attached to $(a,b,c)$.
\end{definition}

\begin{remark}
	We note that, since $g_{r, s_0}^+$ has even degree, one can apply a change of variables as in Lemma~\ref{lem: Legal change of vars} to obtain an integral hyperelliptic equation for $\Cpl(s_0)$, whose RHS polynomial is monic. On the other hand, $\deg (g_{r, s_0}^-)$ is odd, so if one applies a change of variables as in Lemma~\ref{lem: Legal change of vars} to get rid of the denominator of $2-4s_0$, the RHS polynomial cannot be monic. In order to have a monic polynomial on the RHS, one needs to consider a twist of the curve. This motivates the definition of $\Cmsabc$ as a twist of the specialisation, whereas $\Cplabc$ is the specialisation itself.
\end{remark}

Using \eqref{eq: expanded Cms(s)}, we describe the curves $\Cpmabc$ by the hyperelliptic equations \vspace{-0.5em}
\begin{align}
	\hspace{2.5cm} y^2 & = \sum_{k=0}^{\frac{r-1}{2}} (-1)^k c_k \,(Cc)^{2k} \, x^{r-2k} + 2C^{r-1} \left(Bb^p - Aa^p \right), \tag{$\Hhyp_{r}^{-}(a, b, c)$}\label{eq: Hmsabc} \\
	y^2 & = \left( \sum_{k=0}^{\frac{r-1}{2}} (-1)^k c_k \,(Cc)^{2k} \, x^{r-2k} + 2C^{r-1} \left(Bb^p - Aa^p \right) \right)  \left( x + 2Cc\right). \tag{$\Hhyp_{r}^{+}(a, b, c)$}\label{eq: Hplabc}
\end{align}

The reader can find more details about these curves in \cite{Darmon00} and \cite{ChenKoutsianas} for the particular cases $A = B = C = 1$, with an emphasis on the case $r = 5$ in the second reference. When $r = 3$, the elliptic curve $C_3^{+}(a, b, c)$ was studied by Bennett, Vatsal, Yazdani in \cite[\S 2]{BennettVatsalYazdani}. However, the defining model therein is different, and the roles of $Aa^3$ and $Bb^3$ are switched compared to our presentation.

\subsubsection{Frey curve $\Crabc$ for the signature $(r, r, p)$}\label{sect: Frey curves rrp}

The generalised Fermat equation \eqref{eq: GFE rrp bis} of signature $(r, r, p)$ is 
\begin{equation}\label{eq: GFE rrp bis}
	\tag{$\Errp$} Ax^r + By^r = C z^p.
\end{equation}
Assume that there exists a primitive non-trivial solution $(a, b ,c) \in \Z^3$ to \eqref{eq: GFE rrp}. Following \cite{BCDF23}, we construct a Frey curve associated to $(a, b, c)$ as follows. Fix $t_0 \coloneqq Aa^r/Cc^p$, and let $C_r(Aa^r/Cc^p)$ be the specialisation at $t=t_0$ of the curve $C_r(t) / \Q(t)$ introduced in the \S \ref{sect: Frey reps rrp}. Thanks to \eqref{eq: def Crt twist}, we see that $C_r(Aa^r / Cc^p)$ is described by
\begin{equation*}
	y^2 = \sum_{k=0}^{\frac{r-1}{2}} c_k \frac{(ABa^r b^r)^k}{(Cc^p)^{2k}} \, x^{r-2k} + \frac{(-ABa^r b^r)^{\frac{r-1}{2}}(Aa^r - Bb^r)}{(Cc^p)^r}.
\end{equation*}

\begin{definition}
	We define $\Crabc / \Q$ to be the quadratic twist by $- Cc^p / (-ab)^{\frac{r-1}{2}}$ of the specialisation $C_r(Aa^r / Cc^p)$. We call it the \textit{Frey curve} of signature $(r,r,p)$ attached to $(a,b,c)$.
\end{definition}

In the equation above describing $C_r(Aa^r / Cc^p)$, the constant term has a non-trivial denominator. The polynomial on the RHS has odd degree, so if one applies a change of variables as in Lemma~\ref{lem: Legal change of vars} to get rid of this denominator, one obtains a RHS polynomial that is not monic. This motivates the definition of $\Crabc$ as a twist of the specialisation $C_r(Aa^r / Cc^p)$. We describe $\Crabc$ by the hyperelliptic equation
\begin{equation}\label{eq: Hhyp rabc}\tag{$\Hhyp_{r}(a, b)$}
	y^2 = \sum_{k=0}^{\frac{r-1}{2}} c_k \,(ABab)^{k} x^{r-2k} + (AB)^{\frac{r-1}{2}} \left(Bb^r - Aa^r \right).
\end{equation}

This curve generalises a construction by Kraus \cite{KrausC}. In the particular case of trivial coefficients, \ie $A=B=C=1$, the reader can find details on this curve in \cite{BCDF23}. When $r = 3$ and $A, B$ are odd, the elliptic curve $C_3(a, b, c)$ was studied by Bennett and Dahmen (\cite[\S 13]{BennettDahmen}) (we point out that our notation slightly differs from the one therein).
\medskip 

The following elementary remark will have huge Diophantine consequences:

\begin{remark}\label{rmk: Curves 2 paramts}
	Assume first that $(a, b, c)$ is a primitive non-trivial solution to the GFE \eqref{eq: GFE rrp}. The Frey curve attached to $(a, b, c)$ is described by the model \eqref{eq: Hhyp rabc}, which depends only on the coefficients $A, B$ and the members of the solution $a, b$. 
	
	Assume now that $(a, b, c)$ is a primitive non-trivial solution to \eqref{eq: GFE ppr}. We have $Bb^p = Cc^r - Aa^p$. Using this equality, one can rewrite the defining equations \eqref{eq: Hmsabc} and \eqref{eq: Hplabc}, using only the parameters $A, C$ and the members of the solution $a, c$. For instance, the former hyperelliptic equation can be rewritten as 
	\begin{equation}\label{eq: Hmsac}\tag{$\Hhyp_{r}^{-}(a, c)$}
		y^2 = \sum_{k=0}^{\frac{r-1}{2}} (-1)^k c_k \,(Cc)^{2k} \, x^{r-2k} + 2C^{r-1} \left(Cc^r - 2Aa^p \right).
	\end{equation}
	One can also give a model of $\Cplabc$ depending only on $A, C, a, c$. The fact that these Frey curves depend only on $2$ among the $3$ parameters will be used in \S \ref{sect: Solving instances GFEs} for solving many GFEs at the cost of only solving one.
\end{remark}

\subsection{Common framework to both signatures}\label{sect: Common framework}

By construction, the curve $\Cms(a, b, c)$ is a quadratic twist of the specialisation $\Cms(Aa^p / Cc^r)$. Recall from \S \ref{sect: Frey reps rrp} that the curve $C_r(t) / \Kglbar(\B)$ is a quadratic twist of $\Cms(s) / \Kglbar(\B)$. Therefore, $\Crabc$ can also be obtained as a twist of a specialisation of $\Cms(s) / \Q(s)$. In this subsection we introduce a common framework covering both signatures, and we explain how to recover each of the Frey curves associated to solutions by correctly choosing the parameters introduced below.
\medskip

\begin{definition}\label{def: s0 deltaQ}
	We fix $s_0 \in \Qbar \setminus \lbrace 0 , 1 \rbrace$, and $\deltaQ \in \Zbar \setminus \lbrace 0 \rbrace $ such that $\deltaQ^r(2 - 4s_0)$ and $\deltaQ^2$ belong to $\Z$. 
\end{definition}\vspace{-0.5em}

\begin{lemma}\label{lem: s0(s0-1) in Q}
	We have $2^4 \deltaQ^{2r} s_0(s_0-1) \in \Z$, and $s_0 (s_0-1) \in \Q$. If $\deltaQ \in \Z$, then $s_0 \in \Q$.
\end{lemma}

\begin{proof}
	By assumption, $\deltaQ^r (2 - 4s_0) \in \Z$, so its square $2^2 \deltaQ^{2r}(4 s_0(s_0-1) +1)$ also belongs to $\Z$. But $\deltaQ^{2} \in \Z$, so $2^4 \deltaQ^{2r} s_0(s_0-1) = (\deltaQ^r (2 - 4s_0))^2 - 4\deltaQ^{2r}$ also belongs to $\Z$. Finally, dividing the last term by $2^4 \deltaQ^{2r}$, we deduce that $s_0(s_0-1) \in \Q$. The last claim follows easily.
\end{proof}

\begin{definition}\label{def: Curves Cpm + pols}
	Denote by $\Cpm(s_0)$ the specialisation of $\Cpm(s) / \Q(s)$ at $s = s_0$, and by $\CmssQ$ the quadratic twist by $\deltaQ$ of $\Cms(s_0)$. After applying to $\CmssQ$ and $\Cpl(s_0)$ changes of variables as in Lemma~\ref{lem: Legal change of vars}, we describe these curves by the hyperelliptic equations \eqref{eq: expanded Cms} $:  y^2 = \gminus(x)$ and \eqref{eq: expanded Cpl} $: y^2 = \gplus(x)$, where
	\reqnomode
	\begin{equation}\label{eq: pols gminus gplus}
		\gminus(x) \coloneqq \deltaQ^r \, g_{r, s_0}^{-}(x/\deltaQ) \in \Q[x], \qquad \text{and} \qquad \gplus(x) \coloneqq \deltaQ^{r+1} \, g_{r, s_0}^{+}(x/\deltaQ) \in \Q[x].
	\end{equation}
\end{definition}
\medskip

Just as in \eqref{eq: expanded Cms(s)}, we can rewrite the defining equations above in terms of expanded sums
\begin{align}
	y^2 & = \sum_{k=0}^{\frac{r-1}{2}} (-1)^k \, c_k \, \deltaQ^{2k} \,  x^{r-2k} + \deltaQ^r (2 - 4s_0), \label{eq: expanded Cms} \tag{$\Hhyp_{r, s_0}^{- \ (\deltaQ)}$} \\
	y^2 & = \left( \sum_{k=0}^{\frac{r-1}{2}} (-1)^k \, c_k \, \deltaQ^{2k} \, x^{r-2k} + \deltaQ^r (2 - 4s_0) \right) (x+ 2 \deltaQ) \label{eq: expanded Cpl} \tag{$\Hhyp_{r, s_0}^{+ \ (\deltaQ)}$}.
\end{align}

\begin{remark}\label{rmk: Polynomials in Zx}
	As we will see in Proposition~\ref{prop: right choices s0 deltaQ}, in practice, $s_0$ has a non-trivial denominator. We introduce the parameter $\deltaQ$ to obtain integral models for $\CmssQ$ and $\Cpl(s_0)$. The assumptions $\deltaQ^2, \deltaQ^r (2 - 4s_0) \in \Z$ show that $\gminus \in \Z[x]$, and \eqref{eq: expanded Cms} defines a model of $\CmssQ$ over $\Z$. Moreover, if $\deltaQ \in \Z$, then $\gplus \in \Z[x]$, and \eqref{eq: expanded Cpl} also defines a model of $\Cpl(s_0)$ over $\Z$. From now on, we only consider the curve $\Cpl(s_0)$ when $\deltaQ$ is an element of $\Z$. 
\end{remark}

\begin{proposition}\label{prop: right choices s0 deltaQ}
	\begin{enumerate}[leftmargin=*]
		\item\label{enum: s0 deltaQ ppr} Assume that $(a, b, c)$ is a primitive non-trivial solution to \eqref{eq: GFE ppr}. Then the curve $\Cmsabc$ (resp. $\Cplabc$) is the curve $\CmssQ$ (resp. $\Cpl(s_0)$), described by \eqref{eq: expanded Cms} (resp. \eqref{eq: expanded Cpl}), with the specific choices 
		\begin{equation*}
			s_0 = \frac{Aa^p}{Cc^r}, \qquad  \text{ and } \qquad \deltaQ = Cc.
		\end{equation*}
		\item\label{enum: s0 deltaQ rrp} Assume that $(a, b, c)$ is a primitive non-trivial solution to \eqref{eq: GFE rrp}, and let $z_0 \coloneqq \sqrt{-ABab} \in \Qbar$ be a square root of $-ABab$. Then the curve $\Crabc / \Q$ is the curve $\CmssQ$, described by \eqref{eq: expanded Cms}, with the specific choices 
		\begin{equation*}
			s_0 = \frac{1}{2} + \frac{(AB)^{\frac{r-1}{2}} (Aa^r - Bb^r)}{4 z_0^r}, \qquad  \text{ and } \qquad \deltaQ = z_0.
		\end{equation*}		
	\end{enumerate}
\end{proposition}

\begin{proof}
	The first statement follows directly from the construction of $\Cpmabc$ done in \S \ref{sect: Frey curves ppr}. Assume now that $(a,b,c)$ is a primitive non-trivial solution to \eqref{eq: GFE rrp}. Let $t_0 = Aa^r / Cc^p$, consider $s_0$ as in the Proposition, and let $\stilde_0 \coloneqq 2 - 4 s_0 = (AB)^{\frac{r-1}{2}}(Bb^r - Aa^r)/z_0^r$. Recall from \S \ref{sect: Frey curves rrp} that $t_0 = Aa^r / Cc^p$. Since $Aa^r + Bb^r = Cc^p$, we have the equalities
	\begin{equation*}
	 	t_0(t_0-1) = \frac{-ABa^r b^r}{(Aa^r + Bb^r)^2} \qquad \text{ and } \qquad \stilde_0^2 = \frac{(Bb^r - Aa^r)^2}{- ABa^r b^r} = \frac{1}{t_0(t_0-1)} + 4,
	\end{equation*}
	so $(\stilde_0, t_0)$ defines a point in $\B(\Kglbar)$. Recall from \S\ref{sect: Frey curves rrp} that $C_r(t_0)$ is the quadratic twist of $\Cms(s_0)$ by
	\begin{equation*}
		\frac{2 t_0 -1}{\stilde_0} = \frac{(Aa^r - Bb^r) \, z_0^r}{Cc^p \,  (AB)^{\frac{r-1}{2}} (Bb^r - Aa^r)} = \frac{-z_0^r}{Cc^p \, (AB)^{\frac{r-1}{2}}}.
	\end{equation*}
	Moreover, $\Crabc$ is the quadratic twist of $C_r(t_0)$ by $-Cc^p / (-ab)^{\frac{r-1}{2}}$, so composing twists, we conclude that $\Crabc$ is the quadratic twist of $\Cms(s_0)$ by $z_0$.
\end{proof}

Whenever the element $\deltaQ$ does not belong to $\Z$, the curve $\Cpl(s_0)$ does not necessarily admit a model over $\Z$. For instance, given a putative non-trivial solution to \eqref{eq: GFE rrp}, the twisting parameter $\deltaQ = z_0$ is not an element of $\Z$ in general. This is why, when solving the GFE of signature $(r,r,p)$, one cannot employ $\Cpl(s_0)$, and only the curve $\CmssQ$ is available.

\begin{remark}
	We note that, whenever $(a, b, c)$ is a solution to \eqref{eq: GFE rrp}, our choice of $z_0$ as in Proposition~\ref{prop: right choices s0 deltaQ} differs from the one done in \cite{BCDF23}. In \textit{loc. cit.}, it is suggested to choose $z_0 = \sqrt{ABab}$, as the authors of \cite{BCDF23} manipulate the quadratic twist by $i$ of $\Cms(s_0)$. Choosing $z_0$ to be a square root of $-ABab$ allows one to work directly with $\Cms(s_0)$, without having to twist by $i$.
\end{remark}

\begin{remark}\label{rmk: Signature (2, r, p)}
	In \cite{ChenKoutsianas25}, the authors explain that one can attach a hyperelliptic curve to a putative solution to the generalised Fermat equation of signature $(2, r, p)$. Its construction relies on that of Kraus' hyperelliptic curve, introduced in \S \ref{sect: Frey curves rrp}. It would be interesting to study if the hyperelliptic curve associated to a solution to the equation of signature $(2, r, p)$ fits in the framework we introduce here, for some specific choice of $s_0$ and $\deltaQ$.
\end{remark}

In order to solve the generalised Fermat equations \eqref{eq: GFE ppr} and \eqref{eq: GFE rrp}, we will need to consider the base change to $\Kgl$ of the curves introduced above.

\begin{definition}\label{def: Curve Cpm}
	\begin{enumerate}[leftmargin=*]
		\item We let $\Cms$ be the base change of $\CmssQ$ to $\Kgl$, and $\Jms\coloneqq \Jac(\Cms) / \Kgl$ be its Jacobian. We let $\Wms \rightarrow \Spec(\OK)$ be the Weierstrass model of $\Cms / \Kgl$ defined by \eqref{eq: expanded Cms}, and $\Jmodms \rightarrow \Spec(\OK)$ be the Néron model of $\Jms / \Kgl$.
		
		\item If $\deltaQ \in \Z$, we let $\Cpl$ be the base change $\Cpl(s_0)$ to $\Kgl$, and $\Jpl \coloneqq \Jac(\Cpl) / \Kgl$ be its Jacobian. We let $\Wpl \rightarrow \Spec(\OK)$ be the Weierstrass model of $\Cpl / \Kgl$ defined by \eqref{eq: expanded Cpl}, and $\Jmodpl \rightarrow \Spec(\OK)$ be the Néron model of $\Jpl / \Kgl$.
	\end{enumerate}
\end{definition}
\medskip

Although the objects introduced in Definition~\ref{def: Curve Cpm} depend on the parameters $s_0$ and $\deltaQ$, we write them without any reference to the latter in order to simplify the notation. Whenever we make a specific choice for $s_0$ and $\deltaQ$, we will state it in an explicit way. The choices for $s_0$ and $\deltaQ$ done in Proposition~\ref{prop: right choices s0 deltaQ} satisfy indeed the assumptions in Definition~\ref{def: s0 deltaQ}. By considering generic choices of the parameters $s_0, \deltaQ$, we can manipulate the Frey objects for both signatures $(r,r,p)$ and $(p,p,r)$ in a single framework.

We now introduce a new parameter $\deltaK$.

\begin{definition}\label{def: Twisted Cpm}
	Let $\deltaK$ be a non-zero square-free element of $\OK$. We let $\Cpmtw$ be the quadratic twist of $\Cpm$ by $\deltaK$. We describe these curves by the hyperelliptic equations
	\begin{equation*}
		(\Hhyp_r^{-})^{(\deltaK)} \, : \, y^2 = \deltaK \, \gminus(x) \qquad \text{ and } \qquad (\Hhyp_r^{+})^{(\deltaK)} \, : \, y^2 = \deltaK \, \gplus(x)
	\end{equation*}
	We also let $\Wpmtw \rightarrow \Spec(\Oq)$ be the Weierstrass model associated to each of the equations above, $\Jpmtw$ be the Jacobian of $(\Cpm)^{(\deltaK)}$, and $(\Jmodpm)^{(\deltaK)} \rightarrow \Spec(\OK)$ be its Néron model.  
\end{definition}

\begin{remark}
	Lemma~\ref{lem: quad twist has RM} combined with Theorem~\ref{thm: Jpm have RM} shows that the quadratic twist $\Jpmtw$ also has real multiplication by $\Kgl$. Therefore, all the discussion from \S \ref{sect: AV with RM} applies to $\Jpmtw$ too. 
\end{remark}

The interest of introducing this extra twist by $\deltaK$ is to obtain new curves $(\Cpm)^{(\deltaK)}$ having better reduction properties. When solving a Diophantine equation like \eqref{eq: GFE ppr} or \eqref{eq: GFE rrp}, we view the parameters $s_0$ and $\deltaQ$ as fixed, associated to the equation and a putative solution. On the other hand, we treat $\deltaK$ as a variable parameter, whose value we will choose depending on the behaviour of $s_0, \deltaQ$. We search for values of $\deltaK$ such that, for any finite place $\frakq$ of $\Kgl$, the conductor exponent of $\rho_{\Jpmtw , \, \lambda}$ at $\frakq$ is smaller or equal than the one of $\rho_{\Jpm, \, \lambda}$. In practice, $\deltaK$ will be divisible only by primes lying above $2$ and $r$.

\begin{remark}\label{rmk: deltaQr s0(s0-1)}
	To conclude this section, we summarise in Table~\ref{table: Main values} all the important quantities that we will manipulate in the rest of the article for each of the considered signatures. As we will see in Corollary~\ref{cor: Discriminants Wpm}, the discriminants of the models $\Wpm$ depend on $2^4 \deltaQ^{2r} s_0 (s_0-1)$. In the column corresponding to the signature $(p, p, r)$ (resp. $(r, r, p)$), the triple $(a, b, c)$ denotes a primitive non-trivial solution to \eqref{eq: GFE ppr} (resp. \eqref{eq: GFE rrp}).

	\begin{table}[h!]
		\renewcommand{\arraystretch}{1.75}
		\captionsetup{justification=centering,margin=1.5cm}
		\centering
		\begin{tabular}{!{\vrule width 1.25pt}P{3cm}!{\vrule width 1.25pt}P{4.5cm}!{\vrule width 1.25pt}P{4.5cm}!{\vrule width 1.25pt}}
			\noalign{\hrule height 1.25pt}
			Signature  &  $(p, p, r)$  &  $(r, r, p)$  \\  \noalign{\hrule height 1.25pt}
			$s_0$  & $\frac{Aa^p}{Cc^r}$  & $\frac{1}{2} + \frac{(AB)^{\frac{r-1}{2}} (Aa^r - Bb^r)}{4 z_0^r}$ \\  \hline
			$\deltaQ$  &  $Cc$  & $z_0 \coloneqq \sqrt{-ABab}$ \\  \hline
			$\deltaQ^r(2-4s_0)$  & $2 C^{r-1}(Bb^p - Aa^p)$  & $(AB)^{\frac{r-1}{2}}(Bb^r - Aa^r)$ \\  \hline
			$s_0 (s_0-1)$  &  $\frac{-A a^p B b^p}{(Cc^r)^2}$ & $\frac{- (Cc^p)^2}{16 Aa^r Bb^r}$  \\ \hline
			$2^4 \deltaQ^{2r} s_0 (s_0-1)$ & $-2^4 A a^p B b^p C^{2(r-1)}$ & $(AB)^{r-1}(Cc^p)^2$  \\  \noalign{\hrule height 1.25pt}
		\end{tabular}
		\caption{Table describing the explicit quantities introduced above for each of the signatures $(p,p,r)$ and $(r,r,p)$}
		\label{table: Main values}
	\end{table}
\end{remark}

Let $q$ be a rational prime. From now on, we assume that $s_0$ and $\deltaQ$ satisfy the hypotheses below.

\begin{hypothesis}\label{hyp: v(s0(s0-1)) = 0}
	If $v_q(s_0 (s_0-1)) \geq 0$, then $v_q(\deltaQ) = 0$, and $v_q(2 -4s_0) \in \Z$.
\end{hypothesis}

\begin{hypothesis}\label{hyp: val equiv 0 mod r}
	If $q \geq 3$, $v_q(s_0(s_0-1)) \leq 0$ and $v_q(s_0(s_0-1)) \equiv 0 \mod r$, then $v_q(\deltaQ^{2r} s_0 (s_0-1)) = 0$.
\end{hypothesis}

\begin{hypothesis}\label{hyp: 2adic val not -1 -3}
	If $q =2$ and $v_2^{ }(s_0 (s_0-1)) \leq 0$, then $v_2^{ }(s_0(s_0-1)) \notin \lbrace -3, -1 \rbrace$.
\end{hypothesis}

\begin{lemma}\label{lem: s0 deltaQ fill hypoth}
	For each signature, the choice of $s_0, \deltaQ$ done in Proposition~\ref{prop: right choices s0 deltaQ} satisfies Hypotheses~\ref{hyp: v(s0(s0-1)) = 0}, \ref{hyp: val equiv 0 mod r} and \ref{hyp: 2adic val not -1 -3}.
\end{lemma}

\begin{proof}
	By Lemma~\ref{lem: Aa Bb Cc coprime}, $Aa, Bb$ and $Cc$ are pairwise coprime, so Table~\ref{table: Main values} shows that Hypothesis~\ref{hyp: v(s0(s0-1)) = 0} holds. Moreover, we assume that $A, B, C$ are free of $r$-th powers (\cf \ Definition~\ref{def: Parameters GFE}).
	\begin{enumerate}[leftmargin=*]
		\item Assume that $(a, b, c)$ is a primitive non-trivial solution to \eqref{eq: GFE ppr}. If $q \geq 3$, $v_q(s_0(s_0-1)) \leq 0$ and $v_q(s_0(s_0-1)) \equiv 0 \mod r$, then $q \mid c$. Table~\ref{table: Main values} shows that $q \nmid \deltaQ^{2r} s_0(s_0-1)$, so Hypothesis~\ref{hyp: val equiv 0 mod r} is satisfied. If $q = 2$ and $v_{2}^{ }(s_0 (s_0-1)) \leq 0$, then $v_{2}^{ }(s_0 (s_0-1)) = 2 v_{2}^{ }(Cc^r) \in 2 \Z$, so $v_{2}^{ }(s_0(s_0 -1)) \neq -3, -1$.
		
		\item Assume now that $(a, b, c)$ is a primitive non-trivial solution to \eqref{eq: GFE rrp}. If $q \geq 3$, $v_q(s_0(s_0-1)) \leq 0$ and $v_q(s_0(s_0-1)) \equiv 0 \mod r$, then $q \mid ab$. Again, Table~\ref{table: Main values} shows that $q \nmid \deltaQ^{2r} s_0(s_0-1)$. If $q = 2$ and $v_{2}^{ }(s_0 (s_0-1)) \leq 0$, two cases appear. If $v_{2}^{ }(Cc^p) \geq 0$, then $v_{2}^{ }(s_0(s_0-1)) = 2 v_{2}^{ }(Cc^p) -4$, which is even. If $v_{2}^{ }(Cc^p) = 0$, then $v_{2}^{ }(s_0(s_0-1)) \leq -4$. In both cases, we have $v_{2}^{ }(s_0(s_0-1)) \neq -3, -1$.
	\end{enumerate}\vspace{-1em}
\end{proof}

% -- % -- % -- % -- % -- % -- % -- % -- % -- % -- % -- % -- % -- % -- % -- % -- % -- % -- % -- % -- % -- % -- % -- % -- % -- % -- % -- % -- % -- % -- % -- % -- % -- % -- % -- % -- % -- % -- % -- 

\section{Reduction types of the Néron models $\Jmodpmtw$}\label{sect: Reduction types}

In this section we study the geometry of the Néron models $\Jmodpmtw$ as introduced in Definition~\ref{def: Twisted Cpm}, and their reduction types. The main results in here are Theorems~\ref{thm: reduction types Jms} and \ref{thm: reduction types Jpl}, which describe how $\Jmodpmtw$ reduce at any place $\frakq$, depending on the behaviour of $s_0$, $\deltaQ$ and $\deltaK$. Since the Jacobians $\Jpmtw$ have \RM \ by $\Kgl$, they have either good, (totally) toric or (totally) unipotent reduction at $\frakq$, and there is no phenomenon of mixed reduction (see Proposition~\ref{prop: reduction RM} and Remark~\ref{rmk: No mixed red types}). Understanding these reduction types will help us compute Artin conductors of $\rhopmK$ in \S \ref{sect: Cond & local types}, but also to prove modularity of $\rhopmK$ (see \S \ref{sect: Modularity}), and finiteness of the residue representations (\cf \ \S \ref{sect: Level lowering}).
\medskip

As explained in \S \ref{sect: Common framework}, the parameters $s_0$ and $\deltaQ$ are imposed by the considered Diophantine equation and the putative solution. However, the twisting parameter $\deltaK$ is not fixed \textit{a priori}, and we search for values of $\deltaK$ that minimise the Artin conductor of the representations attached to $\Jpmtw$ restricted to $D_{\frakq} \simeq G_{\Kq}$. As explained in the literature \cite{ACIKMM, BCDF23, ChenKoutsianas}, a low semistability defect for $\Jpmtw / \Kq$ tends to lead to a small value of the Artin conductor of $\rhopmK$ at $\frakq$. With this in mind, we aim at minimising the semistability defect of $\Jpmtw$ among all twists of $\Jpm$.
\medskip 

At even places, under some $2$-adic conditions, we study the geometry of different models of the curves $\Cpmtw$, to deduce the reduction types of $\Jmodpmtw$. We use the recent work of Gehrunger \cite{Gehrunger25} and Dokchitser--Morgan \cite{DokchitserMorgan23}, which help understanding the stable reduction of a hyperelliptic curve in even residue characteristic. At odd places, we use the machinery of cluster pictures and the corresponding criteria for having good and semistable reduction. Part of the content of this section is inspired by \cite{ACIKMM}, and the results presented below are a generalisation of the ones established therein. The study of the reduction type of $\Jpmtw$ at even places also makes use of the cluster pictures in some situations. That is why, in \S\ref{sect: Clusters Cpm}, we compute the cluster pictures of the curves $\Cpm$ at all finite places of $\Kgl$, included the even ones. The author would like to thank again Mar Curc\'{o}-Iranzo, Maleeha Khawaja, C\'{e}line Maistret, Diana Mocanu and Tim Gehrunger for fruitful conversations.

\begin{definition}\label{def: nuq}
	Let $\frakq$ be a finite place of $\OK$ and $q \in \Z$ the prime below $\frakq$. We define \vspace{-0.25em}
	\begin{equation*}
		\nu_q \coloneqq v_q(s_0 (s_0-1)). \vspace{-0.25em}
	\end{equation*} 
	Given an element $x \in \Oq \st$, we say that the property $\SQ{x}$ holds if $\vq(x)$ is even, and $x$ is a square $\mod \frakq^2$, \ie there is some $\tau \in \Oq$ such that $\vq(x - \tau^2) \geq 2$.
\end{definition}

\begin{theorem}\label{thm: reduction types Jms}
	The reduction type of $\Jmodmstw$ at $\frakq$ is described in Figure~\ref{fig: Reduction Jmodmstw}.
	\begin{figure}[h!]
		\captionsetup{justification=centering,margin=0.5cm}
		\begin{center}
			\forestset{
		circ/.style={
			circle,
			draw
		},
		rdsq/.style={
			rounded rectangle,
			draw,
			inner sep=4pt,
			font=\footnotesize
		},
		square/.style={
			rectangle,
			draw,
			inner sep=4pt,
			font=\footnotesize
		}
		}

    \scalebox{0.75}{
	\begin{forest}
		for tree={
			grow=east,                    % Tree grows to the right
			edge={->},                    % Edges with arrows
			l sep=6cm,                  % Level separation
			s sep=0.5cm,                    % Sibling separation
			circle,                       % Node shape: circle
			draw,                         % Draw a border around the circle
			inner sep=1.5pt,              % Adjust the size of the circle
			fill=white,                   % Background color
			font=\scriptsize              % Font size for text inside nodes
		}
		[ , circ
			[ , circ, edge label={node[midway,below, sloped]{$\frakq \nmid 2r$}}
				[ , circ, edge label={node[midway,below, sloped]{ $\nu_q \leq 0$, $\nu_q \not \equiv 0 \mod r$}}
				[unipotent, rdsq, edge label={node[midway,above, sloped]{ }}] ]
				[ , circ, edge label={node[midway,above, sloped]{ \quad $\nu_q \leq 0$, $\nu_q \equiv 0 \mod r$}}
					[unipotent, rdsq, edge label={node[midway,below, sloped]{$\vq(\deltaK) =1$}}]
					[good, rdsq, edge label={node[midway,above, sloped]{$\vq(\deltaK) =0$}}] ]
				[ , circ, edge label={node[midway,above, sloped]{$\nu_q > 0$}}
					[unipotent, rdsq, edge label={node[midway,below, sloped]{$\vq(\deltaK) = 1$}}]
					[toric, rdsq, edge label={node[midway,above, sloped]{$\vq(\deltaK) = 0$}}] ] ]
			[ , circ, edge label={node[midway,above, sloped]{$\frakq = \frakr$}}
				[ , circ, edge label={node[midway,below, sloped]{$\nu_r \leq 2$}}
					[unipotent, rdsq, edge label={node[midway,below, sloped]{ }}] ]
				[ , circ, edge label={node[midway,above, sloped]{$\nu_r > 2$}}
					[unipotent, rdsq, edge label={node[midway,below, sloped]{$\vr(\deltaK) = 0$}}]
					[toric, rdsq, edge label={node[midway,above, sloped]{$\vr(\deltaK) = 1$}}] ] ]
			[ , circ, edge label={node[midway,above, sloped]{$\frakq \mid 2$}}
				[ , circ, edge label={node[midway,below, sloped]{$\nu_2 \leq -8, \ \nu_2 \not \equiv -8 \mod r$}}
				[unipotent, rdsq, edge label={node[midway,above, sloped]{ }}]]
				[ , circ, edge label={node[midway,above, sloped]{\quad $\nu_2 \leq -8, \ \nu_2 \equiv -8 \mod r$}}
					[unipotent, rdsq, edge label={node[midway,below, sloped]{not $\SQ{\deltaK \deltaQ^r (2-4s_0)}$}}]
					[good, rdsq, edge label={node[midway,above, sloped]{$\SQ{\deltaK \deltaQ^r (2-4s_0)}$}}] ]
				[ , circ, edge label={node[midway,above, sloped]{$\nu_2 > 0$}}
					[unipotent, rdsq, edge label={node[midway,below, sloped]{$\vq(\deltaK) = 1$}}]
					[toric, rdsq, edge label={node[midway,above, sloped]{$\vq(\deltaK) = 0$}}] ] ]
		]
	\end{forest}}
			\caption{\centering Decision tree describing the reduction type of the Néron model $\Jmodmstw$ at $\frakq$. At each node, follow the branch whose condition is satisfied by the considered parameters.}
			\label{fig: Reduction Jmodmstw}
		\end{center}
	\end{figure}
\end{theorem}

\begin{theorem}\label{thm: reduction types Jpl}
	
	The reduction type of $\Jmodpltw$ at $\frakq$ is described in Figure~\ref{fig: Reduction Jmodpltw}. 
	\begin{figure}
		\captionsetup{justification=centering,margin=0.5cm}
        \forestset{
		circ/.style={
			circle,
			draw
		},
		rdsq/.style={
			rounded rectangle,
			draw,
			inner sep=4pt,
			font=\footnotesize
		},
		square/.style={
			rectangle,
			draw,
			inner sep=4pt,
			font=\footnotesize
		}
	}
		\begin{center}
			\scalebox{0.6}{
            \begin{forest}
		for tree={
			grow=east,                    % Tree grows to the right
			edge={->},                    % Edges with arrows
			l sep=8.5cm,                  % Level separation
			s sep=0.65cm,                    % Sibling separation
			circle,                       % Node shape: circle
			draw,                         % Draw a border around the circle
			inner sep=1.5pt,              % Adjust the size of the circle
			fill=white,                   % Background color
			font=\scriptsize              % Font size for text inside nodes
		}
		[ , circ
			[ , circ, edge label={node[midway,below, sloped]{$\frakq \nmid 2r$}}
				[ , circ, edge label={node[midway,below, sloped]{ $\nu_q \leq 0$, $\nu_q \not \equiv 0 \mod r$}}
					[unipotent, rdsq, edge label={node[midway, above, sloped]{ }}] ]
				[ , circ, edge label={node[midway,above, sloped]{ \quad $\nu_q \leq 0$, $\nu_q \equiv 0 \mod r$}}
					[unipotent, rdsq, edge label={node[midway,below, sloped]{$\vq(\deltaK) = 1$}}]
					[good, rdsq, edge label={node[midway,above, sloped]{$\vq(\deltaK) = 0$}}] ]
				[ , circ, edge label={node[midway,above, sloped]{$\nu_q > 0$}}
					[unipotent, rdsq, edge label={node[midway,below, sloped]{$\vq(\deltaK) = 1 $}}]
					[toric, rdsq, edge label={node[midway,above, sloped]{$\vq(\deltaK) = 0$}}] ] ]
			[ , circ, edge label={node[midway,above, sloped]{$\frakq = \frakr$}}
				[ , circ, edge label={node[midway,below, sloped]{$\nu_r \leq 2, \ \nu_r \not \equiv 0 \mod r$}}
					[unipotent, rdsq, edge label={node[midway,below, sloped]{ }}] ]
				[ , circ, edge label={node[midway,above, sloped]{$\vr(s_0-1) > 2$}}
					[unipotent, rdsq, edge label={node[midway,below, sloped]{$\vr(\deltaK) = 0$}}]
					[toric, rdsq, edge label={node[midway,above, sloped]{$\vr(\deltaK) = 1$}}] ] 
				[ , circ, edge label={node[midway,above, sloped]{$\vr(s_0) > 2$}}
					[unipotent, rdsq, edge label={node[midway,below, sloped]{$\vr(\deltaK) = 1$}}]
					[toric, rdsq, edge label={node[midway,above, sloped]{$\vr(\deltaK) = 0$}}] ] ]
			[ , circ, edge label={node[midway,above, sloped]{$\frakq \mid 2$}}
				[ , circ, edge label={node[midway,below, sloped]{$\nu_2 \leq -4, \ \nu_2 \not \equiv 0 \mod r$}}
					[unipotent, rdsq, edge label={node[midway,above, sloped]{ }}]]
				[ , circ, edge label={node[midway,above, sloped]{\quad $\nu_2 \leq -4, \ \nu_2 \equiv 0 \mod r$}}
					[unipotent, rdsq, edge label={node[midway,below, sloped]{not $\SQ{\deltaK}$}}]
					[good, rdsq, edge label={node[midway,above, sloped]{$\SQ{\deltaK}$}}] ] 
				[ , circ, edge label={node[midway,above, sloped]{$\nu_2 > 0$}}
					[unipotent, rdsq, edge label={node[midway,below, sloped]{$\vq(\deltaK) = 1$}}]
					[toric, rdsq, edge label={node[midway,above, sloped]{$\vq(\deltaK) = 0$}}] ]] 
		]
	\end{forest}
            }
			\caption{\centering Decision tree describing the reduction type of the Néron model $\Jmodpltw$ at $\frakq$. At each node, follow the branch whose condition is satisfied by the considered parameters.}
			\label{fig: Reduction Jmodpltw}
		\end{center}
	\end{figure}
\end{theorem}

\begin{example}\label{examp: Reduction types Fermat}
	We can specialise the content of Theorems~\ref{thm: reduction types Jms} and \ref{thm: reduction types Jpl} for the specific choices of $s_0, \deltaQ$ done in \S \ref{sect: Common framework} for each of the signatures $(p,p,r)$ and $(r,r,p)$. For instance, we describe the reduction types of the Néron models of the Jacobians considered in \S \ref{sect: Frey curves ppr} and \S \ref{sect: Frey curves rrp} at $\frakq \nmid 2r$. Describing this reduction types for $\frakq \mid 2$, or $\frakq= \frakr$ can be done in a similar way.
	\begin{enumerate}[leftmargin=*, itemsep=0pt]
		\item[] \boxppr Assume that $(a, b, c)$ is a primitive non-trivial solution to \eqref{eq: GFE ppr}. The choice of $s_0, \deltaQ$ done in Proposition~\ref{prop: right choices s0 deltaQ} gives the equality $2^4 s_0 (s_0-1) =  \frac{-2^4 Aa^p Bb^p}{(Cc^r)^2}$. If $\frakq \nmid 2r$, $q$ lies below $\frakq$, and $\deltaK$ is such that $\vq(\deltaK) = 0$, then the Néron model of $\Jac(\Cpmabc \times \Kgl)^{(\deltaK)}$ has 
		\begin{equation*}
			\begin{cases}
				\text{toric reduction at } \frakq & \text{ if } q \mid ABab, \\
				\text{good reduction at } \frakq & \text{ if } q \mid c \text{ and } q \nmid C, \\
				\text{unipotent reduction at } \frakq & \text{ if } q \mid C .
			\end{cases}
		\end{equation*}
	
		\item[] \boxrrp  Assume that $(a, b, c)$ is a primitive non-trivial solution to \eqref{eq: GFE rrp}. The choice of $s_0, \deltaQ$ done in Proposition~\ref{prop: right choices s0 deltaQ} gives the equality $2^4 s_0 (s_0-1) =  \frac{(Cc^p)^2}{Aa^r Bb^r}$. If $\frakq \nmid 2r$, $q$ lies below $\frakq$, and  $\deltaK$ is such that $\vq(\deltaK) = 0$, then the Néron model of $\Jac(\Crabc \times \Kgl)^{(\deltaK)}$ has
		\begin{equation*}
			\begin{cases}
				\text{toric reduction at } \frakq & \text{ if } q \mid Cc, \\
				\text{good reduction at } \frakq & \text{ if } q \mid ab \text{ and } q \nmid AB, \\
				\text{unipotent reduction at } \frakq & \text{ if } q \mid AB.
			\end{cases}
		\end{equation*}
	\end{enumerate}
\end{example}

\subsection{Roots of the defining polynomials $\gpm$}\label{sect: Computing roots}

In this subsection we exhibit algebraic expressions for the roots of the polynomials $\gpm(x) \in \Z[x]$ introduced in Definition~\ref{def: Curves Cpm + pols}, and we deduce the discriminants of the Weierstrass models $\Wpmtw$. Recall that we only consider the curve $\Cpl$ (and thus $\gplus$) when $\deltaQ \in \Z$. In this case, the extra root of $\gplus$ is $-2 \deltaQ$, which belongs to $\Z$ (see Remark~\ref{rmk: Polynomials in Zx} ). Therefore, we focus on describing the roots of the polynomial $\gminus(x)$.
\medskip 

All the content below is local, and concerns the curves $\Cpm$, which we consider as defined over $\Kq$. From now on, we treat $s_0, \deltaQ, \deltaK$ as belonging to $\Qqbar$, and we view $\gminus$ as an element of $\Qq[x]$, through the choice of an embedding $\Qbar \hookrightarrow \Qqbar$.

\begin{definition}\label{def: alpha0 beta0}
	Let $\sqts, \sqtsm \in \Qqbar$ be square roots of $s_0$ and $s_0 -1$ respectively, and let $\sqt \coloneqq \sqts \, \sqtsm$. Define $\alpha_0 \in \Qqbar$ to be an $r$-th root of \vspace{-0.5em}
	\begin{equation*}
		\left( \sqts + \sqtsm \right)^2 = 2 s_0 -1 + 2 \sqt. \vspace{-0.5em}
	\end{equation*}
	For any $j \in \Iintv{0, r-1}$, define $\alpha_j \coloneqq \zr^j \, \alpha_0$, $\beta_j \coloneqq 1/\alpha_j$ and $\gamma_j\coloneqq \deltaQ (\alpha_j + \beta_j) \in \Qqbar$. 
	For simplicity, we write $\gamma_{-j} \coloneqq \gamma_{r-j}$. Moreover, we set $\gamma_r \coloneqq -2 \deltaQ$.
\end{definition}

\begin{remark}\label{rmk: alpha0 rational}
	By definition of $\deltaQ$ and $s_0$, $\deltaQ \alpha_0$ is an $r$-th root of $\deltaQ^r (\sqts + \sqtsm)^2$. By \cite[Section VI, Theorem 9.1]{Lang}, the polynomial $x^{r} - \deltaQ^r(\sqts + \sqtsm)^2 \in \Qq(\deltaQ^r \sqt)[x]$ is reducible if and only if it has a root in $\Qq(\deltaQ^r	 \sqt)$. We choose the embedding $\Qbar \hookrightarrow \Qqbar$ in such a way that, whenever the polynomial above is reducible, then $\deltaQ \alpha_0 \in \Qq(\deltaQ^r \sqt)$ (see Proposition~\ref{prop: gamma0 rational} for further details). Therefore, the extension $\Qq(\deltaQ \alpha_0) / \Qq(\deltaQ^r \sqt)$ has degree either $1$ or $r$. 
\end{remark}

\begin{lemma}\label{lem: useful eqs}
	For any $j \in \Iintv{0, r-1}$, we have the following properties:
	\begin{enumerate}[itemsep=0pt]
		\itemsep0em 
		\item\label{enum: useful eqs1} The element $\beta_j$ is an $r$-th root of $\left( \sqts - \sqtsm \right)^2 = 2 s_0 -1 - 2 \sqt$. 
		\item\label{enum: useful eqs2} $\alpha_j^r + \beta_j^r = 4s_0 -2$.
		\item\label{enum: useful eqs3} $\alpha_j^r - \beta_j^r = 4 \sqt$.
		\item\label{enum: useful eqs4} $\left( \alpha_j^r + 1 \right)^2 / \alpha_j^r = 4s_0$.
	\end{enumerate}
\end{lemma}

\begin{proof}
	These are algebraic manipulations that follow from Definition~\ref{def: alpha0 beta0}.
\end{proof}

\begin{proposition}
	For any $j \in \Iintv{0, r-1}$, $\gamma_j$ is a root of $\gminus(x) \in \Qq[x]$. 
\end{proposition}

\begin{proof}
	Fix $j \in \Iintv{0, r-1}$. Since $\gminus(x) = \deltaQ^r \, g_{r, s_0}^{-}(x/\deltaQ)$, it suffices to check that $\alpha_j  + \beta_j$ is a root of $g_{r, s_0}^{-}$. Recall that $\alpha_j \beta_j = 1$ and $\alpha_j^r + \beta_j^r = 4s_0-2$. Let us recall that \cite[Lemma 2.14]{ACIKMM} gives the equality of polynomials
	\begin{equation*}
		(-X Y)^{\frac{r-1}{2}} (X + Y) h_{r}\left(2 - \frac{(X + Y)^2}{X Y}\right) = X^r + Y^r. 
	\end{equation*}
	Combining the latter with Lemma~\ref{lem: useful eqs}\eqref{enum: useful eqs2}, we compute
	\begin{equation*}
		g_{r, s_0}^-(\alpha_j + \beta_j)  = (-\alpha_j \beta_j)^{\frac{r-1}{2}} \left(\alpha_j + \beta_j \right) \, h_r \left(2 - \frac{(\alpha_j + \beta_j)^2}{\alpha_j \beta_j} \right) + 2 -4 s_0  = \alpha_j^r + \beta_j^r + 2-4 s_0 = 0. \vspace{-1.5em}
	\end{equation*}
\end{proof}

\begin{remark}
	The notation used above is slightly different from the one employed in \cite{ACIKMM}. When $A=B=C=1$ and $(a,b,c)$ is a primitive non-trivial solution to \eqref{eq: GFE ppr}, we have $s_0 = a^p/c^r$, $\deltaQ = c$ with the choice of Proposition~\ref{prop: right choices s0 deltaQ}. We recover the elements that are denoted by $\alpha_j, \beta_j$ in \cite[\S 4.1]{ACIKMM} as $\deltaQ \alpha_j, \deltaQ \beta_j$ in our notation.
\end{remark}

\begin{proposition}\label{prop: dif 2 roots}
	For any $0 \leq j, k \leq r-1$, we have the equalities \vspace{-0.5em}
	\begin{equation*}
		\gamma_k - \gamma_j = \deltaQ \, \zeta_{r}^{k}(1-\zeta_{r}^{j-k})(\alpha_{0} - \zeta_{r}^{-j-k}\beta_{0}) \qquad \text{ and } \qquad \gamma_j - \gamma_r = \frac{\deltaQ \, \zr^j}{\alpha_0} \left(\alpha_0 + \zr^{-j} \right)^2.
	\end{equation*}
	For any $k \in \Iintv{0, r-1}$, we have
	\begin{equation*}
		\prod_{\substack{0 \leq j \leq r-1 \\ j \neq k}} (\gamma_k - \gamma_j) = \frac{4 r \deltaQ^{r-1} \sqt}{\zr^k (\alpha_0 - \zr^{-2k} \beta_0)} \qquad \text{ and } \qquad \prod_{j=0}^{r-1} (\gamma_j - \gamma_r) = 4 \deltaQ^r s_0.
	\end{equation*}
	In particular, the $\gamma_j$'s are pairwise distinct.
\end{proposition}

\begin{proof}
	The first identity is obtained by developing the right-hand side. Using such equality, we compute, for $k \in \Iintv{0, r-1}$
	\begin{equation*}
		\prod_{\substack{0 \leq j \leq r-1 \\ j \neq k}} (\gamma_k - \gamma_j) = \prod_{\substack{0 \leq j \leq r-1 \\ j \neq k}} \deltaQ \, \zr^k (1 - \zr^{j-k}) (\alpha_0 - \zr^{-j-k} \beta_0) = \frac{\deltaQ^{r-1} \zr^{k(r-1)} r \, (\alpha_0^r - \beta_0^r)}{\alpha_0 - \zr^{-2k} \beta_0},
	\end{equation*}
	and the result follows from Lemma~\ref{lem: useful eqs} \eqref{enum: useful eqs3}.  On the other hand, using $\beta_0 = 1/ \alpha_0$, we compute
	\begin{equation*}
		\frac{\deltaQ \, \zr^j}{\alpha_0} \left(\alpha_0 + \zr^{-j} \right)^2 = \frac{\deltaQ \, \zr^j}{\alpha_0} \left( \alpha_0^2 + 2 \, \zr^{-j} \alpha_0 + \zr^{-2j} \right) = \deltaQ  \left( \zr^j \alpha_0 + 2 +\zr^{-j} \beta_0  \right) =  \gamma_j - \gamma_r.
	\end{equation*}
	Now Lemma~\ref{lem: useful eqs} \eqref{enum: useful eqs4} implies $\prod_{j=0}^{r-1} (\gamma_j - \gamma_r) = \deltaQ^r \, \alpha_0^{-r}(\alpha_0^r + 1)^2 = 4 \deltaQ^r s_0$. Finally, the $\gamma_j$'s are pairwise distinct because we have $s_0 \neq 0 ,1$, and $\deltaQ \neq 0$, so the products of the differences are non-zero.
\end{proof}

\begin{definition}
	We denote by $\Rroots \coloneqq \lbrace \gamma_0, \ldots, \gamma_{r-1} \rbrace$ the set of roots of $\gminus$, and by $\Rroots^+ \coloneqq \Rroots \cup \lbrace \gamma_r \rbrace$ the set of roots of $\gplus$.
\end{definition}

\begin{corollary}\label{cor: Discriminants Wpm}
	The discriminants of the Weierstrass models $\Wpmtw$ are given by 
	\begin{align*}
		\Delta \left( \Wmstw \right) & = (-1)^{\frac{r-1}{2}} \, 2^{2(r-1)} \, r^r \, \left(2^4 \deltaQ^{2r} s_0(s_0-1)\right)^{\frac{r-1}{2}} \, \deltaK^{2r}, \\
		\Delta \left( \Wpltw \right) & = (-1)^{\frac{r-1}{2}} \, 2^{2(r-1)} \, r^r \, \left(2^4 \deltaQ^{2r} s_0(s_0-1)\right)^{\frac{r-1}{2}} (4 \deltaQ^r s_0)^2 \, \deltaK^{2r}. \vspace{-0.5em}
	\end{align*}
	In particular, the places of bad reduction for $\Wpmtw$ divide $2, r$, $\deltaQ^{2r} s_0(s_0 -1)$ and $\deltaK$.
\end{corollary}

\begin{proof}
	We begin by computing the discriminants of $\Wpm$, which satisfy $\Delta(\Wpm) = 2^{4g} \disc(\gpm)$, where $g = \frac{r-1}{2}$ is the genus of $\Cpm$, and $\disc(\gpm)$ is the discriminant of $\gpm$. First we have
	\begin{align*}
		\disc(\gminus) & = (-1)^{r \frac{r-1}{2}} \prod_{\substack{0 \leq j, k \leq r-1 \\ j \neq k}} (\gamma_k - \gamma_j) \\
		& = (-1)^{\frac{r-1}{2}} \prod_{k = 0}^{r-1} \frac{4 \deltaQ^{r-1} r \sqt}{\zr^k (\alpha_0 - \zr^{-2k} \beta_0)} = (-1)^{\frac{r-1}{2}} \frac{\left( 4 \deltaQ^{r-1} r \sqt\right)^r }{4 \sqt} 
	\end{align*}
	so regrouping terms yields the displayed result. On the other hand, the discriminants of $\gminus$ and $\gplus$ are related by $\disc(\gplus) = \disc(\gminus) \prod_{j= 0}^{r-1}(\gamma_j - \gamma_r)^2 = \disc(\gminus) \, 4^2 \deltaQ^{2r} s_0^2$. The discriminants of $\Wpmtw$ are obtained from those of $\Wpm$ by applying Lemma~\ref{lem: Legal change of vars}.
\end{proof}

\begin{remark}
	Recall that we assume that $s_0, \deltaQ$ satisfy Hypothesis~\ref{hyp: val equiv 0 mod r}. The latter implies that, when $\frakq \nmid 2r$, if $\nu_q \leq 0$ and $\nu_q \equiv 0 \mod r$, then $\frakq \nmid \deltaQ^{2r} s_0(s_0-1)$. Therefore, the curves $\Cpm$ then have good reduction at $\frakq$, and the same holds for $\Cpmtw$ as $\deltaK$ is supported at places where $\Cpm$ have bad reduction.
\end{remark}

\begin{example}\label{examp: Discriminants Fermat}
	We now specialise the content of Definition~\ref{def: alpha0 beta0} and Corollary~\ref{cor: Discriminants Wpm} to the specific choices of $s_0, \deltaQ$ done in Proposition~\ref{prop: right choices s0 deltaQ} (see also Table~\ref{table: Main values}).
	\begin{enumerate}[leftmargin=*]
		\item[] \boxppr Assume that $(a, b, c)$ is a non-trivial primitive solution to the GFE \eqref{eq: GFE ppr}. Let us fix square roots $\sqrt{Aa^p}$ and $\sqrt{-Bb^p}$ for $Aa^p$ and $-Bb^p$ respectively. We have the equalities
		\begin{equation*}
			\alpha_0^r =  \frac{1}{Cc^r} \left( \sqrt{Aa^p} + \sqrt{-Bb^p} \right)^2 \qquad \text{ and } \qquad \beta_0^r = \frac{1}{Cc^r} \left( \sqrt{Aa^p} - \sqrt{-Bb^p} \right)^2. 
		\end{equation*}
		Combining Corollary~\ref{cor: Discriminants Wpm} and Remark~\ref{rmk: deltaQr s0(s0-1)}, we describe the discriminants of \eqref{eq: Hmsabc} and \eqref{eq: Hplabc}, the equations describing the curves $\Cpmabc$ (\cf \ \S \ref{sect: Frey curves ppr}):		
		\begin{align*}
			\Delta(\Hhyp_{r}^{-}(a, b, c)) & = 2^{4(r-1)} \, r^r (Aa^p)^{\frac{r-1}{2}} (Bb^p)^{\frac{r-1}{2}} C^{(r-1)^2}, \\ 
			\Delta(\Hhyp_{r}^{+}(a, b, c)) & = 2^{4r} \, r^r (Aa^p)^{\frac{r+3}{2}} (Bb^p)^{\frac{r-1}{2}} C^{r^2-1}.
		\end{align*}
		
		\item[] \boxrrp Assume that $(a, b, c)$ is a non-trivial primitive solution to the GFE \eqref{eq: GFE rrp}. We have
		\begin{equation*}
			\alpha_0^r =  \frac{(AB)^{\frac{r-1}{2}} \, Aa^r}{z_0^r} \qquad \text{ and } \qquad \beta_0^r = - \frac{ (AB)^{\frac{r-1}{2}} \, Bb^r}{z_0^r}. 
		\end{equation*}
		Just as above, Corollary~\ref{cor: Discriminants Wpm} and Remark~\ref{rmk: deltaQr s0(s0-1)} yield the discriminant of \eqref{eq: Hhyp rabc}, the equation defining $\Crabc$ (see \S \ref{sect: Frey curves rrp}):
		\begin{equation*}
			\Delta(\Hhyp_{r}(a, b)) = (-1)^{\frac{r-1}{2}} 2^{2(r-1)} \, r^r (Cc^p)^{r-1} (AB)^{\frac{(r-1)^2}{2}}.
		\end{equation*}
	\end{enumerate}
	In the particular case of trivial coefficients $A = B = C = 1$, we recover the results of \cite{Darmon00, BCDF23, ChenKoutsianas, ACIKMM} for both signatures.
\end{example}

\subsection{Cluster pictures of $\Cpm$}\label{sect: Clusters Cpm}

We now describe the cluster pictures of $\Cpm$ at all places of $\Kgl$. The curves $\Cpl(s_0)$ and $\CmssQ$ are defined over $\Q$, and one could compute their cluster pictures over $\Qq$, in the style of \cite{ACIKMM}. However, for the Diophantine applications we address, we only need to know the cluster pictures of the base-changed curves $\Cpm / \Kq$. In particular, in \S \ref{sect: Conclusion red types} we will exploit the fact that the base changed Jacobians $\Jpm / \Kq$ have \RM \ by $\Kgl$, whereas those defined over $\Qq$ do not. We note that the cluster pictures of the twisted curves $\Cpmtw$ are the same ones as those of $\Cpm$. The only difference appears in the valuation of the leading coefficient of the defining polynomial.
\medskip

The theory of cluster pictures concerns hyperelliptic curves over local fields of odd residue characteristic. Nevertheless, Dokchitser and Morgan \cite{DokchitserMorgan23} have recently showed that, for some very specific cluster pictures, one can describe the stable reduction of hyperelliptic curves in even residue characteristic. Similarly, the recent work of Gehrunger~\cite{Gehrunger25} provides combinatorial tools to understand this stable reduction, and some of his results may be phrased in terms of cluster pictures. This is why we compute and draw these cluster pictures at all places of $\Kgl$, included the even ones. However, we will only be able to extract arithmetic information from the cluster pictures at even places in some particular cases (those of (totally) toric reduction). For this reason, we only compute the cluster pictures of $\Cpm$ at even places for certain specific choices of parameters $s_0, \deltaQ$.
\medskip

Recall that we set in Definition~\ref{def: nuq} $\nu_q \coloneqq v_q(s_0(s_0-1))$. We keep denoting by $\frakq$ a finite place of $\Kgl$ and $q$ the rational prime lying below $\frakq$. We do not assume that $\frakq$ is an odd place of $\Kgl$, as some of the results below will apply to the particular case where $\frakq \mid 2$. We begin by describing the $\frakq$-adic valuation of $\alpha_0$ and $\beta_0$ as introduced in Definition~\ref{def: alpha0 beta0}.

\begin{lemma}\label{lem: valuations alpha0 beta0}
	We have the following properties:
	\begin{enumerate}[itemsep=0pt, topsep=5pt]
		\item\label{enum: val alpha0 beta 0 i} If $\nu_q \geq 0$, then $\vq(\alpha_0) = \vq(\beta_0) = 0$.
		\item\label{enum: val alpha0 beta 0 ii} If $q$ is odd and $\nu_q < 0$ then $\vq(s_0 - 1) = \vq(s_0) < 0$, and $\left \lbrace \vq(\alpha_0) , \vq(\beta_0)  \right \rbrace = \left \lbrace \frac{1}{r} \vq(s_0) , - \frac{1}{r} \vq(s_0)  \right \rbrace$.
	\end{enumerate}
\end{lemma}

\begin{proof}
	\begin{enumerate}[leftmargin=*]
		\item If $\nu_q > 0$, then exactly one among $\vq(s_0)$ and $\vq(s_0-1)$ is zero and the other is positive. The statement follows from the equalities $\alpha_0^r = \left(\sqts + \sqtsm \right)^2$ and $\beta_0^r = \left(\sqts - \sqtsm \right)^2$. If $\nu_q = 0$, assume by contradiction that $\vq(\alpha_0) \neq 0$. We have $\beta_0 = 1/ \alpha_0$, so up to switching $\alpha_0, \beta_0$, we may assume that $\vq(\alpha_0) > 0$ and $\vq(\beta_0) < 0$. Using Lemma~\ref{lem: useful eqs} \eqref{enum: useful eqs3}, we obtain
		\begin{equation*}
			0 = \frac{1}{2} \, \vq(s_0(s_0 -1)) = \vq(\alpha_0^r - \beta_0^r) = \min(r \, \vq(\alpha_0), r \, \vq(\beta_0)) = r \, \vq(\beta_0) < 0,
		\end{equation*}
		hence a contradiction. We conclude that $\vq(\alpha_0) = \vq(\beta_0) = 0$.
		
		\item If $\nu_q < 0$, then at least one among $s_0$ and $s_0 - 1$ has negative valuation, and thus both of them have the same valuation. By definition of $\alpha_0, \beta_0$ and Lemma~\ref{lem: useful eqs} \eqref{enum: useful eqs1}, we have
		\begin{equation*}
			\vq(\alpha_0^r) = \vq(s_0) + 2 \vq \left(1 +  \frac{\sqtsm}{\sqts}  \right) \quad \text{ and } \quad \vq(\beta_0^r) = \vq(s_0) + 2 \vq \left(1 -  \frac{\sqtsm}{\sqts} \right)
		\end{equation*}
		But $\sqtsm / \sqts$ is a square root of $1 - \frac{1}{s_0}$. Since $\frac{1}{s_0}$ belongs to $\frakq$ and $\frakq$ is an odd place of $\Kgl$, we can compute the Taylor expansion of the square root of $1 - \frac{1}{s_0}$ in the ring of integers $\Oq$. There is some $\eps \in \lbrace \pm 1 \rbrace$ such that $\frac{\sqtsm}{\sqts} = \eps (1 - \frac{1}{2s_0} + \ldots)$ and further powers of $\frac{1}{s_0}$. We deduce that exactly one among 
		\begin{equation*}
			1 +  \frac{\sqtsm}{\sqts} \quad \text{ and } \quad 1 -  \frac{\sqtsm}{\sqts} 
		\end{equation*}
		has valuation $0$ and the other has valuation $- \vq(s_0)$. Therefore, we conclude that exactly one among $\alpha_0^r$ and $\beta_0^r$ has $\frakq$-adic valuation $\vq(s_0)$, and the other has valuation $- \vq(s_0)$.
	\end{enumerate}\vspace{-1em}
\end{proof}

Recall that we denote by $\frakr = (2 - \omega)$ the unique prime of $\OK$ dividing $r$. The valuation $\vr$ being normalised with respect to $\Kr$, we have $\vr(1-\zr)=\frac{1}{2}$. As explained in Proposition~\ref{prop: dif 2 roots}, the difference of roots $\gamma_k - \gamma_j$ is divisible by $\alpha_0 - \zr^{-j-k} \beta_0$. We now describe the valuation of the latter.

\begin{lemma}\label{lem: val alpha - zr beta}
	We have the following properties.
	\begin{enumerate}[leftmargin=*, itemsep=0pt]
		\item\label{enum: val alpha - zr beta i} If $\, \frac{1}{2} \vq(s_0(s_0-1)) > r \, \vq(1-\zr)$, there is a unique $j_0 \in \Iintv{0,r-1}$ such that $\vq(\alpha_0 - \zr^{j_0} \beta_0) > \vq(1 - \zr)$, and $\vq(\alpha_0 - \zr^j \beta_0) = \vq(1-\zr)$ for $j \in \Iintv{0, r-1} \setminus \lbrace j_0 \rbrace $. Moreover, if $\zr \notin \Kq$, then $j_0 = 0$.
		
		\item\label{enum: val alpha - zr beta ii} If $\frakq$ is odd and $\, \frac{1}{2} \vq(s_0(s_0-1)) \leq r \, \vq(1-\zr)$, then $\vq(\alpha_0 - \zr^j \beta_0) = \vq(s_0(s_0-1))/2r$ for all $j \in \Iintv{0, r-1}$.	
	\end{enumerate}
\end{lemma}

\begin{proof}
	Recall from Lemma~\ref{lem: useful eqs} \eqref{enum: useful eqs3} that $\alpha_0^r - \beta_0^r  = 4\sqt$, so $\vq(\alpha_0^r - \beta_0^r) = \frac{1}{2} \vq(2^4 s_0(s_0-1))$.
	\begin{enumerate}[leftmargin=*]
		\item If $\frakq$ is odd and $\frac{1}{2}\vq(s_0(s_0-1)) > r \vq(1 - \zr) \geq 0$, then Lemma~\ref{lem: valuations alpha0 beta0} \eqref{enum: val alpha0 beta 0 i} yields $\vq(\alpha_0) = \vq(\beta_0) = 0$. The statement is a particular case of \cite[Lemma 2.15]{ACIKMM}, applied with $K = \Kq$, $v = \vq$, $\alpha = \alpha_0$ and $\beta = -\beta_0$. Although the mentioned Lemma is stated in the context of an odd prime $q$, it is easy to check that the proof is exactly the same when $q=2$, so we can indeed use it as above.
		
		\item Assume that $\frakq$ is odd and that $\frac{1}{2}\vq(s_0(s_0-1)) \leq r \vq(1-\zr)$. We claim that the valuation of $\alpha_0 - \zr^j \beta_0$ does not depend on $j$. Different possibilities arise. \vspace{-0.5em}
		\begin{enumerate}[leftmargin=*]
			\item If $\vq(\alpha_0 - \zr^j \beta_0) \geq \vq(1 - \zr)$ for every $j \in \Iintv{0, r-1}$, then we have \vspace{-0.5em}
			\begin{equation*}
				r \vq(1 - \zr) \leq \sum_{j=0}^{r-1} \vq(\alpha_0 - \zr^j \beta_0) = \vq(\alpha_0^r - \beta_0^r) = \frac{1}{2} \vq(s_0(s_0-1)) \leq r \vq(1 - \zr). \vspace{-0.75em}
			\end{equation*}
			Thus, $\frac{1}{2}\vq(s_0(s_0-1)) = r \vq(1-\zr)$, and $\vq(\alpha_0 - \zr^j \beta_0) = \frac{1}{r} \vq(\alpha_0^r - \beta_0^r)$ for all $j \in \Iintv{0, r-1}$.
			
			\item Suppose that there is some $j_1 \in \Iintv{0, r-1}$ such that $\vq(\alpha_0 - \zr^{j_1} \beta_0) < \vq(1 - \zr)$.
			\begin{enumerate}
				\item If $\vq(\beta_0) \geq 0$ then equality $\alpha_0 - \zr^j \beta_0 = (\alpha_0 - \zr^{j_1} \beta_0 )  + \zr^{j_1}(1 - \zr^{j - j_1}) \beta_0$ shows that $\vq(\alpha_0 - \zr^j \beta_0) = \vq(\alpha_0 - \zr^{j_1} \beta_0)$.
				
				\item  If $\vq(\beta_0) < 0$, equality $\alpha_0 \beta_0 = 1$ implies that $\vq(\alpha_0)>0$, and thus $\vq(\alpha_0 - \zr^j \beta_0) = \vq(\beta_0)$.
			\end{enumerate}
		\end{enumerate} \vspace{-0.25em}
		In all cases, $\vq(\alpha_0 - \zr^j \beta_0)$ is independent of $j$, and therefore equals $ \frac{1}{r} \vq(\alpha_0^r - \beta_0^r) = \frac{1}{2r}\vq(s_0(s_0-1))$.
	\end{enumerate}\vspace{-0.75em}
\end{proof}

The second statement in the Lemma above is presented only in the case where the place $\frakq$ is odd. One could consider an analogue statement when $\frakq$ is even, but here we would not make use of such a statement. Indeed, in such a situation, we do not know how to extract arithmetic information at an even place from the corresponding cluster picture.

\begin{remark}
	Lemma~\ref{lem: val alpha - zr beta} describes the behaviour of $\vq(\alpha_0 - \zr^j \beta_0)$ under certain conditions on $\vq(s_0 (s_0-1))$. If $\frakq \neq \frakr$, the criterion is $\vq(s_0(s_0-1))> 0$, and if $\frakq = \frakr$, the criterion is $\vq(s_0(s_0-1))> r$. Since $\vq(1-\zr)=0$ if $\frakq \neq \frakr$, and $\vr(1-\zr)=\frac{1}{2}$, the condition $\frac{1}{2} \vq(s_0(s_0-1)) > r \vq(1-\zr)$ recovers both criteria above, whether $\frakq = \frakr$ or not.
\end{remark}

\begin{definition}\label{def: i0 j0}
	If $\frakq \neq \frakr$ and $\vq(s_0(s_0-1)) > 0$, we let $j_0 \in \Iintv{0, r-1}$ be as in Lemma~\ref{lem: val alpha - zr beta} \eqref{enum: val alpha - zr beta i}. If $\frakq = \frakr$ and $\vr(s_0(s_0-1)) > r$, we let $j_0 = 0$. Since $\zr \notin \Kr$, we have $\vq(\alpha_0 - \zr^{j_0} \beta_0) > \vq(1 - \zr)$ and $\vq(\alpha_0 - \zr^{j} \beta_0) = \vq(1 - \zr)$ for $j \neq j_0$, whether $\frakq = \frakr$ or not. We define $i_0 \in \Iintv{0, r-1}$ to be such that $-2i_0 \equiv j_0 \mod r$. 
\end{definition}

\begin{theorem}\label{thm: Cluster Cms}
	Let $\frakq$ be a place of $\Kgl$ dividing $2r \deltaQ^{2r} s_0 (s_0-1)$. Let $n \coloneqq \frac{\vq(2^4 s_0(s_0-1))}{2}$ and $m \coloneqq n - \frac{r}{2}$. 
	\begin{enumerate}
		\item If $\frakq \neq \frakr$ and $\vq(s_0(s_0-1)) > 0$, then the cluster picture of $\Cms / \Kq$ is
		\begin{center}
\scalebox{1.05}{
	\clusterpicture            
	\Root[D] {3} {first} {r2};
	\Root[D] {6} {r2} {r3};
	\ClusterLDName c1[][n][\gamma_0, \ \gamma_{2i_0 }] = (r2)(r3);
	\Root[D] {4} {c1} {r4};
	\Root[D] {6} {r4} {r5};
	\ClusterLDName c2[][n][\gamma_1, \ \gamma_{2i_0 - 1}] = (r4)(r5);
	\Root[E] {8} {c2} {r6};
	\Root[E] {2} {r6} {r7};
	\Root[E] {2} {r7} {r8};
	\Root[D] {8} {r8} {r9};
	\Root[D] {6} {r9} {r10};
	\ClusterLDName c3[][n][\gamma_{\frac{r-1}{2}}, \ \gamma_{2 i_0 - \frac{r-1}{2}}] = (r9)(r10);
	\Root[D] {8} {c3} {r11};
	\ClustercLDName c5[][][\gamma_{i_0}] = (r11);
	\ClusterLD c4[][\, 0] = (c1)(c1n)(c2)(c2n)(r6)(r7)(r8)(c3)(c3n)(c5)(c5n);
	\endclusterpicture 
}

		\end{center}

		Note that, when $k = i_0$, the roots $\gamma_k$ and $\gamma_{2i_0 - k}$ are equal, so the twin $\lbrace \gamma_k, \gamma_{2i_0 - k} \rbrace$ is simply the singleton $\lbrace \gamma_{i_0} \rbrace$.
		
		\item If $\frakq \nmid 2r$ and $\vq(s_0(s_0-1)) \leq 0$, then the cluster picture of $\Cms / \Kq$ is
		\begin{center}
\scalebox{1.25}{
		\clusterpicture            
		\Root[D] {2} {first} {r2};
		\Root[D] {4} {r2} {r3};
		\Root[D] {4} {r3} {r4};
		\Root[E] {4} {r4} {r5};
		\Root[E] {1} {r5} {r6};
		\Root[E] {1} {r6} {r7};
		\Root[D] {4} {r7} {r8};
		\Root[D] {4} {r8} {r9};
		\Root[D] {4} {r9} {r10};
		\ClusterLD c1[][\, v_{\frakq}(\delta_{\Q}) +  \frac{n}{r}] = (r2)(r3)(r4)(r5)(r6)(r7)(r8)(r9)(r10);
		\endclusterpicture
}
		\end{center}

		\item If $\frakq = \frakr$ and $\vr(s_0(s_0-1)) > r$, then the cluster picture of $\Cms / \Kr$ is 
		\begin{center}
        \scalebox{1.1}{
	\clusterpicture         
	\Root[D] {3} {first} {r2};
	\Root[D] {5} {r2} {r3};
	\ClusterLDName c1[][m][\gamma_1, \ \gamma_{-1}] = (r2)(r3);
	\Root[D] {3} {c1} {r4};
	\Root[D] {5} {r4} {r5};
	\ClusterLDName c2[][m][\gamma_2, \ \gamma_{-2}] = (r4)(r5);
	\Root[E] {8} {c2} {r6};
	\Root[E] {2} {r6} {r7};
	\Root[E] {2} {r7} {r8};
	\Root[D] {8} {r8} {r9};
	\Root[D] {5} {r9} {r10};
	\ClusterLDName c3[][m][\gamma_{\frac{r-1}{2}}, \ \gamma_{\frac{r+1}{2}}] = (r9)(r10);
	\Root[D] {5} {c3} {r11};
	\ClustercLDName c5[][][\gamma_0] = (r11);
	\ClusterLD c4[][\, 1 ] = (c1)(c1n)(c2)(c2n)(r6)(r7)(r8)(c3)(c3n)(c5)(c5n);
	\endclusterpicture
}
		\end{center}

		\item If $\frakq = \frakr$ and $\vr(s_0(s_0-1)) \leq r$, then the cluster picture of $\Cms / \Kr$ is 
		\begin{center}
		\scalebox{1.3}{
		\clusterpicture            
		\Root[D] {2} {first} {r2};
		\Root[D] {4} {r2} {r3};
		\Root[D] {4} {r3} {r4};
		\Root[E] {4} {r4} {r5};
		\Root[E] {1} {r5} {r6};
		\Root[E] {1} {r6} {r7};
		\Root[D] {4} {r7} {r8};
		\Root[D] {4} {r8} {r9};
		\Root[D] {4} {r9} {r10};
		\ClusterLD c1[][\, \frac{1}{2} + v_{\frakr}(\delta_{\Q}) + \frac{n}{r}] = (r2)(r3)(r4)(r5)(r6)(r7)(r8)(r9)(r10);
		\endclusterpicture
	}
		\end{center}
	\end{enumerate}
\end{theorem}

\begin{proof}
	We use the equality $\gamma_k - \gamma_j = \deltaQ \zr^{k} (1 - \zr^{j-k}) (\alpha_0 - \zr^{-j-k} \beta_0)$ from Proposition~\ref{prop: dif 2 roots} to compute the valuation of the pairwise differences of the roots.
	\begin{enumerate}[leftmargin=*]
		\item Suppose that $\frakq \neq \frakr$, and $\vq (s_0(s_0-1)) > 0$. By Hypothesis~\ref{hyp: v(s0(s0-1)) = 0}, we have $\vq(\deltaQ) = 0$. Let $i_0, j_0$ be as in Definition~\ref{def: i0 j0}. For any $j, k \in \Iintv{0, r-1}$, we have $\vq(\alpha_0 - \zr^{-j-k} \beta_0) > 0$ if and only if $-j -k \equiv -2i_0 \mod r$, \ie $j \equiv 2 i_0 - k \mod r$. If $j \not \equiv 2 i_0 - k \mod r$ then $\vq(\gamma_k - \gamma_j) = 0$, as $\vq(\deltaQ) = 0 = \vq(1- \zr)$. On the other hand, Lemma~\ref{lem: val alpha - zr beta} gives
		\begin{equation*}
			\vq(\gamma_{k} - \gamma_{2i_0 - k}) = \vq(\alpha_0 - \zr^{j_0} \beta_0) = \vq(\alpha_0^r - \beta_0^r) = \frac{1}{2} \vq(2^4 s_0(s_0-1)) =n.
		\end{equation*}
		The outer depth of the cluster picture is therefore $0$. Moreover, any choice of $k \neq i_0$ satisfies $- k - i_0 \not \equiv -2 i_0 \mod r$, implying that $\gamma_{i_0}$ is an isolated root in the cluster picture. Finally, there are $\frac{r-1}{2}$ twins of (relative) depth $n$, each of them consisting of the roots $\lbrace \gamma_{k}, \gamma_{2i_0 - k} \rbrace$, where $k\neq i_0$.
		
		\item Suppose that $\frakq \nmid 2r$ and $\vq(s_0(s_0-1)) \leq 0$. Lemma~\ref{lem: val alpha - zr beta} gives $\vq(\alpha_0 - \zr^{-j-k} \beta_0) = \frac{\vq(s_0(s_0-1))}{2r} = \frac{n}{r}$ for any $j, k \in \Iintv{0, r-1}$ with $j \neq k$, and thus $\vq(\gamma_k - \gamma_j) = \vq(\deltaQ) + n/r$. We conclude that all the roots lie in a single cluster of depth $\vq(\deltaQ) + n/r$.

		\item Suppose that $\frakq = \frakr$ and $\vr(s_0(s_0-1)) > r$. Again, by Hypothesis~\ref{hyp: v(s0(s0-1)) = 0}, we have $\vr(\deltaQ) = 0$.	Since $\zr \notin \Kr$, Lemma~\ref{lem: valuations alpha0 beta0} \eqref{enum: val alpha0 beta 0 i} implies that, for any $j, k \in \Iintv{0, r-1}$ such that $j \not \equiv -k \mod r$, we have $\vr(\alpha_0 - \zr^{-j-k} \beta_0) = \vr(1 - \zr)$, so $\vr(\gamma_k - \gamma_j) = 2 \vr(1 - \zr) = 1$. Just as in the first case, we deduce that $\gamma_0$ is an isolated root in the cluster picture. Using Lemma~\ref{lem: useful eqs} \eqref{enum: useful eqs3}, we compute 
		\begin{equation*}
			\vr(\alpha_0 - \beta_0) = \vr(\alpha_0^r - \beta_0^r) - (r-1) \vr(1 - \zr) = \frac{\vr(s_0 (s_0-1))}{2}  - \frac{r-1}{2} = n - \frac{r}{2} + \frac{1}{2}.
		\end{equation*}
		Therefore, for any $k \in \Iintv{0, \frac{r-1}{2}}$ we have $\vr(\gamma_k - \gamma_{-k}) = \vr(1 - \zr) + \vr(\alpha_0 - \beta_0)= 1 + n - \frac{r}{2}$. Hence, the outer depth of the cluster picture is $1$, and there are $\frac{r-1}{2}$ twins of relative depth $m = n - \frac{r}{2}$.

		\item Suppose that $\frakq = \frakr$ and $\vr(s_0(s_0-1)) \leq  r$. By Lemma~\ref{lem: val alpha - zr beta}, we have $\vr(\alpha_0 - \zr^{-j-k} \beta_0) = \frac{\vr(s_0(s_0-1))}{2r}$ for any $j, k \in \Iintv{0, r-1}$ with $j \neq k$, and thus $\vr(\gamma_k - \gamma_j) = \vr(1 - \zr) + \vr(\deltaQ) + n/r$. Again, we conclude that all the roots lie in a single cluster with depth $\frac{1}{2} + \vq(\deltaQ) + n/r$.
	\end{enumerate}\vspace{-1em}
\end{proof}

The second statement in the theorem above is only presented for odd places whereas the first one also includes the case $\frakq \mid 2$. Again, this is because in the case of the first statement we will be able to deduce the reduction type of $\Cpm$ from the cluster picture at an even place. On the other hand, in the case of the second statement we lack of theoretical arguments to deduce the reduction type from the cluster picture. 
\medskip

Proposition~\ref{prop: dif 2 roots} shows that the difference of roots $\gamma_j - \gamma_r$ is divisible by $\frac{1}{\alpha_0}(\alpha_0 + \zr^{-j})^2$.
We now describe the valuation of the latter in terms of $s_0$, in order to draw the cluster picture of $\Cpl / \Kq$. 
Recall from Definition~\ref{def: i0 j0} that $i_0$ is the unique element in $\Iintv{0, r-1}$ satisfying $\vq(\alpha_0 - \zr^{-2i_0} \beta_0) > \vq(1 - \zr)$.

\begin{lemma}\label{lem: val alpha0 + zr}
	We have the following properties.
	\begin{enumerate}[leftmargin=*, itemsep=0pt]
		\item\label{enum: val alpha + zr i} If $\ \frac{1}{2} \vq(s_0) > r \vq(1- \zr)$, then $\vq \left(\frac{1}{\alpha_0}(\alpha_0 + \zr^{-i_0})^2 \right) > 2 \vq(1-\zr)$ and $\vq \left(\frac{1}{\alpha_0}(\alpha_0 + \zr^{-j})^2 \right) = 2 \vq(1-\zr)$ for any $j \in \Iintv{0, r-1} \setminus \lbrace i_0 \rbrace$.
		
		\item\label{enum: val alpha + zr ii} If $\frakq$ is odd and $\, \frac{1}{2} \vq(s_0) \leq r \vq(1- \zr)$, then $\vq \left(\frac{1}{\alpha_0}(\alpha_0 + \zr^{-j})^2 \right) = \frac{1}{r} \vq(s_0)$ for all $j \in \Iintv{0, r-1}$.
	\end{enumerate}
\end{lemma}

\begin{proof}
	Recall from Lemma~\ref{lem: useful eqs} \eqref{enum: useful eqs4} that $(\alpha_0^r +1)^2 / \alpha_0^r = 4s_0$, so $\vq \left( (\alpha_0^r +1)^2 / \alpha_0^r \right) = \vq(s_0)$.
	\begin{enumerate}[leftmargin=*]
		\item Assume that $\frac{1}{2} \vq(s_0) > r \vq(1 - \zr)$. Lemma~\ref{lem: valuations alpha0 beta0} \eqref{enum: val alpha0 beta 0 i} yields $\vq(\alpha_0) = 0$, so $\vq(\alpha_0^r + 1) > r \vq(1-\zr)$. Lemma 2.15 from \cite{ACIKMM} applied with $K = \Kq$, $v = \vq$, $\alpha = \alpha_0$, $\beta = 1$ states that there is a unique $k_0 \in \Iintv{0, r-1}$ such that $\vq(\alpha_0 + \zr^{k_0}) > \vq(1-\zr)$ (again, we use the mentioned lemma even when $q$ is even). We claim that $k_0 \equiv -i_0 \mod r$. Indeed \vspace{-0.25em}
		\begin{equation*}
			\left(\alpha_0 + \zr^{k_0} \right) - \left(\alpha_0 - \zr^{-2 i_0} \beta_0 \right) = \frac{\zr^{k_0}}{\alpha_0} \left( \alpha_0 + \zr^{-2 i_0 - k_0} \right)
		\end{equation*}
		gives $\vq(\alpha_0 + \zr^{-2 i_0 - k_0}) > \vq(1-\zr)$, since
		\begin{equation*}
			\vq(\alpha_0 - \zr^{- 2i_0} \beta_0) = \vq(\alpha_0 - \zr^{j_0} \beta_0) > \vq(1 - \zr).
		\end{equation*}
		It follows that $-2i_0 - k_0 \equiv k_0 \mod r$, and therefore $k_0 \equiv - i_0 \mod r$.
		
		\item Assume now that $\frakq$ is odd and that $\frac{1}{2} \vq(s_0) \leq r \vq(1 - \zr)$, so that $\vq \left( (\alpha_0^r +1)^2 / \alpha_0^r \right) \leq 2r \vq(1-\zr)$. We claim that the valuation of $(\alpha_0 + \zr^{-j})^2 / \alpha_0$ does not depend on $j$. 
		Different cases arise.
		\begin{enumerate}[leftmargin=*]
			\item If $\vq ((\alpha_0 + \zr^{-j})^2 / \alpha_0 ) \geq 2 \vq(1 - \zr)$ for every $j \in \Iintv{0, r-1}$, then we have \vspace{-0.5em}
			\begin{equation*}
				2 r \vq(1 - \zr) \leq \sum_{j=0}^{r-1} \vq \left( \frac{1}{\alpha_0} \left( \alpha_0 + \zr^{-j} \right)^2 \right) = \vq \left( \frac{(\alpha_0^r + 1)^2}{\alpha_0^r}\right) = \vq(s_0) \leq 2r \vq(1 - \zr). \vspace{-0.75em}
			\end{equation*}
			Thus, $\vq(s_0) = 2r \vq(1 - \zr)$ and $\vq ((\alpha_0 + \zr^{-j})^2 / \alpha_0) = \frac{1}{r}\vq((\alpha_0^r +1)^2 / \alpha_0^r)$ for all $j \in \Iintv{0, r-1}$. \vspace{-0.5em}
			
			\item Suppose that there is some $k_1 \in \Iintv{0, r-1}$ such that $\vq ((\alpha_0 + \zr^{-k_1})^2 / \alpha_0 ) < 2 \vq(1 - \zr)$. Fix $j \in \Iintv{0, r-1}$. Again, we treat different cases separately. 
			\begin{enumerate}
				\item If $\vq(\alpha_0) = 0$, then $\vq(\alpha_0 + \zr^{-k_1}) < \vq(1-\zr)$. Equality	
				\begin{equation*}
					\alpha_0 + \zr^{-j} = ( \alpha_0 + \zr^{-k_1} ) + \zr^{k_1} ( 1 - \zr^{-j - k_1})
				\end{equation*}
				implies that $\vq(\alpha_0 + \zr^{-j}) = \vq(\alpha_0 + \zr^{-k_1})$, so $\vq((\alpha_0 + \zr^{-j})^2 / \alpha_0) = 2 \vq(\alpha_0 + \zr^{-k_1})$.
				
				\item If $\vq(\alpha_0) > 0$, then $\vq(\alpha_0 + \zr^{-j}) = 0$, so $\vq((\alpha_0 + \zr^{-j})^2 / \alpha_0) = - \vq(\alpha_0)$.
				
				\item If $\vq(\alpha_0) < 0$, then $\vq(\alpha_0 + \zr^{-j}) = \vq(\alpha_0)$, and thus $\vq((\alpha_0 + \zr^{-j})^2 / \alpha_0) = \vq(\alpha_0)$.
			\end{enumerate}
		\end{enumerate}
		In all cases, we see that $\vq((\alpha_0 + \zr^{-j})^2 / \alpha_0)$ is independent of $j$, and it is therefore equal to $\frac{1}{r} \vq((\alpha_0^r +1)^2 / \alpha_0^r) = \frac{1}{r} \vq(s_0)$. 
	\end{enumerate} \vspace{-1.5em}
\end{proof}

We now proceed to describe the cluster picture of the curve $\Cpl$ at $\frakq$.

\begin{theorem}\label{thm: Cluster Cpl}
	Let $\frakq$ be a place of $\Kgl$ dividing $2r \deltaQ^{2r} s_0 (s_0-1)$. Let us define $n \coloneqq \frac{\vq(2^4 s_0 (s_0-1))}{2}$, $\widetilde{n} \coloneqq \vq(4s_0)$,  and $m \coloneqq n - \frac{r}{2}$. 
	\begin{enumerate}
		\item If $\frakq \neq \frakr$ and $\vq(s_0(s_0 -1)) > 0$, then the cluster picture of $\Cpl / \Kq$ is
		\begin{center}
			\scalebox{1.1}{
	\clusterpicture            
	\Root[D] {3} {first} {r2};
	\Root[D] {6} {r2} {r3};
	\ClusterLDName c1[][n][\gamma_1, \ \gamma_{2i_0-1}] = (r2)(r3);
	\Root[D] {4} {c1} {r4};
	\Root[D] {6} {r4} {r5};
	\ClusterLDName c2[][n][\gamma_2, \ \gamma_{2i_0- 2}] = (r4)(r5);
	\Root[E] {8} {c2} {r6};
	\Root[E] {2} {r6} {r7};
	\Root[E] {2} {r7} {r8};
	\Root[D] {8} {r8} {r9};
	\Root[D] {6} {r9} {r10};
	\ClusterLDName c3[][n][\gamma_{\frac{r-1}{2}}, \ \gamma_{2i_0- \frac{r-1}{2}}] = (r9)(r10);
	\Root[D] {6} {c3} {r11};
	\Root[D] {6} {r11} {r12};
	\ClusterLDName c4[][\widetilde{n}][\gamma_{i_0}, \ \  \gamma_r] = (r11)(r12);
	\ClusterLD c4[][\, 0] = (c1)(c1n)(c2)(c2n)(r6)(r7)(r8)(c3)(c3n)(c4)(c4n);
	\endclusterpicture 
}
		\end{center}

		In particular, when $\frakq$ is odd and $\vq(s_0-1) > 0$, then $\widetilde{n} = 0$, so $\gamma_{i_0}$ and $\gamma_r$ are both isolated roots.

		\item If $\frakq \nmid 2 r$ and $\vq(s_0(s_0-1)) \leq 0$, then the cluster picture of $\Cpl / \Kq$ is 
		\begin{center}
	\scalebox{1.25}{
	\clusterpicture            
	\Root[D] {2} {first} {r2};
	\Root[D] {4} {r2} {r3};
	\Root[D] {4} {r3} {r4};
	\Root[E] {4} {r4} {r5};
	\Root[E] {1} {r5} {r6};
	\Root[E] {1} {r6} {r7};
	\Root[D] {4} {r7} {r8};
	\Root[D] {4} {r8} {r9};
	\Root[D] {4} {r9} {r10};
	\Root[D] {4} {r10} {r11}
	\ClusterLD c1[][\, v_{\frakq}(\delta_{\Q}) + \frac{n}{r}] = (r2)(r3)(r4)(r5)(r6)(r7)(r8)(r9)(r10)(r11);
	\endclusterpicture
}
		\end{center}

		\item If $\frakq = \frakr$ and $\vq(s_0) > r$, then the cluster picture of $\Cpl / \Kr$ is 
		\begin{center}
	\scalebox{1.05}{
	\clusterpicture            
	\Root[D] {3} {first} {r2};
	\Root[D] {6} {r2} {r3};
	\ClusterLDName c1[][m][\gamma_1, \ \gamma_{-1}] = (r2)(r3);
	\Root[D] {4} {c1} {r4};
	\Root[D] {6} {r4} {r5};
	\ClusterLDName c2[][m][\gamma_2, \ \gamma_{- 2}] = (r4)(r5);
	\Root[E] {8} {c2} {r6};
	\Root[E] {2} {r6} {r7};
	\Root[E] {2} {r7} {r8};
	\Root[D] {8} {r8} {r9};
	\Root[D] {6} {r9} {r10};
	\ClusterLDName c3[][m][\gamma_{\frac{r-1}{2}}, \ \gamma_{- \frac{r-1}{2}}] = (r9)(r10);
	\Root[D] {6} {c3} {r11};
	\Root[D] {6} {r11} {r12};
	\ClusterLDName c4[][2m][\gamma_{0}, \ \gamma_r] = (r11)(r12);
	\ClusterLD c4[][\, 1] = (c1)(c1n)(c2)(c2n)(r6)(r7)(r8)(c3)(c3n)(c4)(c4n);
	\endclusterpicture 
}
		\end{center}

		\item If $\frakq = \frakr$ and $\vr(s_0 -1)>r$, then the cluster picture of $\Cpl / \Kr$ is 
		\begin{center}
			\scalebox{1.05}{
	\clusterpicture            
	\Root[D] {3} {first} {r2};
	\Root[D] {6} {r2} {r3};
	\ClusterLDName c1[][m][\gamma_1, \ \gamma_{-1}] = (r2)(r3);
	\Root[D] {4} {c1} {r4};
	\Root[D] {6} {r4} {r5};
	\ClusterLDName c2[][m][\gamma_2, \ \gamma_{-2}] = (r4)(r5);
	\Root[E] {8} {c2} {r6};
	\Root[E] {2} {r6} {r7};
	\Root[E] {2} {r7} {r8};
	\Root[D] {8} {r8} {r9};
	\Root[D] {6} {r9} {r10};
	\ClusterLDName c3[][m][\gamma_{\frac{r-1}{2}}, \ \gamma_{- \frac{r-1}{2}}] = (r9)(r10);
	\Root[D] {6} {c3} {r11};
	\ClustercLDName c7[][][\gamma_{0}] = (r11);
	\ClusterLD c4[][1] = (c1)(c1n)(c2)(c2n)(r6)(r7)(r8)(c3)(c3n)(c7)(c7n);
	\Root[D] {4} {c4} {r12};
	\ClustercLDName c6[][][\gamma_r] = (r12);
	\ClusterLD c5[][\, 0] = (c1)(c1n)(c2)(c2n)(r6)(r7)(r8)(c3)(c3n)(c4)(c6)(c6n);
	\endclusterpicture
}
		\end{center}

		\item If $\frakq = \frakr$ and $\vr(s_0(s_0-1)) \leq r$, then the cluster picture of $\Cpl / \Kr$ is
		\begin{center}
			\scalebox{1.25}{
	\clusterpicture            
	\Root[D] {2} {first} {r2};
	\Root[D] {4} {r2} {r3};
	\Root[D] {4} {r3} {r4};
	\Root[E] {4} {r4} {r5};
	\Root[E] {1} {r5} {r6};
	\Root[E] {1} {r6} {r7};
	\Root[D] {4} {r7} {r8};
	\Root[D] {4} {r8} {r9};
	\Root[D] {4} {r9} {r10};
	\ClusterLD c1[][\frac{1}{2} + \widetilde{m}] = (r2)(r3)(r4)(r5)(r6)(r7)(r8)(r9)(r10);
	\Root[D] {3} {c1} {r12};
	\ClustercLDName c5[][][\gamma_r] = (r12);
	\ClusterLD c4[][\,  v_{\frakr}(\delta_{\Q}) + \frac{n}{r}] = (c1)(r6)(r7)(r8)(c5)(c5n);
	\endclusterpicture
}
		\end{center}
		where $\widetilde{m} \coloneqq \frac{1}{2r} \left( \vr(s_0 -1) - \vr(s_0) \right)$.
	\end{enumerate}
\end{theorem}

\begin{proof}
	By Theorem~\ref{thm: Cluster Cms}, it remains to compute the $\frakq$-adic valuations of $\gamma_j - \gamma_r$ for $j \in \Iintv{0,r-1}$. We recall from Proposition~\ref{prop: dif 2 roots} that $\gamma_j - \gamma_r = \frac{\deltaQ \zr^j}{\alpha_0} (\alpha_0 + \zr^{-j})^2$.
	\begin{enumerate}[leftmargin=*]
		\item Suppose that $\frakq \neq \frakr$ and $\vq(s_0(s_0-1)) > 0$. Hypothesis~\ref{hyp: v(s0(s0-1)) = 0} yields $\vq(\deltaQ) = 0$. For $j \in \Iintv{0, r-1} \setminus \lbrace i_0 \rbrace$, we have $\vq((\alpha_0 + \zr^{-j})^2 / \alpha_0) = 0$ by Lemma~\ref{lem: val alpha0 + zr} \eqref{enum: val alpha + zr i}, so $\vq(\gamma_j - \gamma_r) = 0$. On the other hand,
		\begin{equation*}
			\vq \left(\frac{(\alpha_0 + \zr^{-i_0})^2}{\alpha_0} \right) = \vq \left( \frac{(\alpha_0^r + 1)^2}{\alpha_0^r} \right) = \vq(4s_0)
		\end{equation*}
		by Lemma~\ref{lem: useful eqs} \eqref{enum: useful eqs4}. It follows that $\vq(\gamma_{i_0} - \gamma_r) = \vq(4 s_0) = \widetilde{n}$, so $\gamma_{i_0}$ and $\gamma_r$ lie in a twin of depth $\vq(s_0) = \widetilde{n}$. If $\frakq$ is odd and $\vq(s_0) = 0$, then this twin has depth $0$, and therefore $\gamma_{i_0}, \gamma_r$ are isolated roots.

		\item Suppose that $\frakq \nmid 2r$ and $\vq(s_0(s_0-1)) \leq 0$. Lemma~\ref{lem: val alpha0 + zr} \eqref{enum: val alpha + zr ii} yields $\vq((\alpha_0 + \zr^{-j})^2 / \alpha_0) = \vq(s_0)/r$ for any $j \in \Iintv{0, r-1}$, hence $\vq(\gamma_j - \gamma_r) = \vq(\deltaQ) + \vq(s_0)/r$. We conclude that all the roots of $\gplus$ lie in a common cluster of depth $\vq(\deltaQ) + \vq(s_0)/r = \vq(\deltaQ) + n/r$.
		
		\item Suppose that $\frakq = \frakr$ and $\vq(s_0) > r$, so that $\vr(s_0 -1) = 0$ and $n = \vr(s_0)/2$. Hypothesis~\ref{hyp: v(s0(s0-1)) = 0} gives $\vr(\deltaQ) = 0$, and Lemma~\ref{lem: valuations alpha0 beta0} \eqref{enum: val alpha0 beta 0 i} implies that $\vr(\alpha_0) = 0$. Moreover, $i_0 = 0$ as $\zr \notin \Kr$. Thus, we have $\vr((\alpha_0 + \zr^{-j})^2 / \alpha_0) = 2\vr(1 - \zr)$ for any $j \in \Iintv{1, r-1}$ (\cf \ Lemma~\ref{lem: val alpha0 + zr} \eqref{enum: val alpha + zr i}), and $\vr((\alpha_0 + 1)^2 / \alpha_0) > 2 \vr(1-\zr)$. Thanks to Lemma~\ref{lem: useful eqs} \eqref{enum: useful eqs4}, we compute 
		\begin{equation*}
			\vr(\alpha_0 + 1) = \vr(\alpha_0^r + 1) - (r-1) \vr(1 - \zr) = \frac{1}{2} \, \vr(s_0) + \frac{1 - r}{2}.
		\end{equation*} 
		It follows that $\vr(\gamma_j - \gamma_r) = 1$ for any $j \neq 0$, and $\vr(\gamma_0 - \gamma_r) = 1 + \vr(s_0) -r$. Therefore, the outer depth remains $1$, and $\gamma_0$ and $\gamma_r$ lie in a twin of relative depth $\vr(s_0) -r = 2m$.
		
		\item Suppose that $\frakq = \frakr$ and $\vr(s_0 - 1) > r$, which implies $\vr(s_0) = 0$, and $\vr(\deltaQ) = 0$ (\cf \ Hypothesis~\ref{hyp: v(s0(s0-1)) = 0}). For any $j \in \Iintv{0,r-1}$ we have $\vr((\alpha_0 +\zr^{-j})^2 / \alpha_0) = 0$ (Lemma~\ref{lem: val alpha0 + zr} $ii)$), so $\vr(\gamma_j - \gamma_r) = 0$. Therefore, $\gamma_r$ is an isolated root lying outside of the cluster picture of $\gminus$, and the outer depth is $0$.
		
		\item Suppose that $\frakq = \frakr$ and $\vr(s_0(s_0-1)) \leq r$. Lemma~\ref{lem: val alpha0 + zr} \eqref{enum: val alpha + zr ii} yields $\vr((\alpha_0 + \zr^{-j})^2 / \alpha_0) = \vr(s_0)/r$ for all $j \in \Iintv{0, r-1}$, so $\vr(\gamma_j - \gamma_r) = \vr(\deltaQ) + \vr(s_0)/r$. Fix $j, k \in \Iintv{0, r-1}$ with $j \neq k$. Recall from Theorem~\ref{thm: Cluster Cms} that $\vr(\gamma_j - \gamma_k) = \frac{1}{2} + \vr(\deltaQ) + \frac{1}{2r} \vr(s_0(s_0-1))$. We claim that \vspace{-0.25em}
		\begin{equation*}
			\vr(\gamma_j - \gamma_k) \geq \vr(\gamma_j - \gamma_r), \qquad \text{ or, equivalently, } \qquad r + \vr(s_0-1) \geq \vr(s_0). \vspace{-0.25em}
		\end{equation*}
		Let us prove that the latter is true by treating the different cases separately.
		\begin{enumerate}
			\item If $ 0 < \vr(s_0) \leq r$, then $\vr(s_0-1) = 0$ and the inequality becomes $\vr(s_0) \leq r$. 
			
			\item If $ 0 < \vr(s_0 -1) \leq r$, then $\vr(s_0) = 0$, and the inequality becomes $r + \vr(s_0-1) \geq 0$.
			
			\item If $\vr(s_0(s_0-1)) \leq 0$, then we have $\vr(s_0) = \vr(s_0-1)$, and the inequality is simply $r \geq 0$.
		\end{enumerate}
		We conclude that the outer depth of the cluster picture is $\vr(\deltaQ) + \vr(s_0)/r$, and $\gamma_r$ is an isolated root in it. The cluster containing the roots of $\gminus$ has relative depth \vspace{-0.25em}
		\begin{equation*}
			\vr(\gamma_j - \gamma_k) - \vr(\gamma_j - \gamma_r) = \frac{1}{2} + \frac{1}{2r} \left( \vr(s_0-1) - \vr(s_0)\right) \geq 0.
		\end{equation*}
	\end{enumerate}\vspace{-1em}
\end{proof}

\subsection{Local behaviour at even places}\label{sect: Red even places}

In this subsection we assume that the finite place $\frakq$ is even and we describe the reduction type of the Néron models $\Jmodpmtw$. We know that $\Cpmtw$ have potential semistable reduction, and since their Jacobians $\Jpmtw$ have real multiplication by $\Kgl$, they have either potential good or potential (totally) toric reduction. We treat separately these two cases depending on the parameter $\nu_2 = v_2(s_0(s_0-1))$. In the former case, we use the recent work of Gehrunger \cite{Gehrunger25} and Dokchitser--Morgan \cite{DokchitserMorgan23}, who study the stable reduction of hyperelliptic curves in even residue characteristic. Their work allows us to recover the reduction type of $\Jmodpmtw$ at $\frakq$ from the cluster picture of $\Cpm$ at $\frakq$. For the case of potential good reduction, we analyse the geometry of different models of $\Cpmtw$, in order to deduce the reduction type of $\Jmodpmtw$ at $\frakq$, and the extension of $\Kq$ where $\Jpmtw$ attains good reduction.

\subsubsection{The case of (potential) toric reduction}

We first treat the case of potential toric reduction for $\Jmstw$.

\begin{proposition}\label{prop: Toric red Jmodmstw}
	Assume that $\nu_2 > 0$. If $\vq(\deltaK) = 0$, then $\Jmodmstw$ has toric reduction at $\frakq$. Otherwise, it has unipotent reduction at $\frakq$, and it attains toric reduction over $\Kq(\sqrt{\deltaK})$.
\end{proposition}

\begin{proof}
	Assume first that $\vq(\deltaK) = 0$. We use the work of Gehrunger \cite{Gehrunger25} to show that $\Jmodmstw$ has totally toric reduction at $\frakq$. Let us recall from Theorem~\ref{thm: Cluster Cms} the description of the cluster picture of $\Cms$ at $\frakq$:  the outer depth of the cluster picture is $0$, and there are $\frac{r-1}{2}$ twins. Each of these twins has depth $\frac{1}{2} \vq(2^4 s_0(s_0-1))$, and the assumption $\nu_2 = v_2(s_0 (s_0 - 1)) > 0$ implies that this common depth is $>4$. Moreover, assuming that $\vq(\deltaK) = 0$ implies that the defining polynomial is monic. In the terminology of \cite[\S 5.2]{Gehrunger25}, such a description of the cluster picture means that $\Cmstw$ has $\frac{r-1}{2}$ grounded double points with thickness strictly larger than $4$. Then \cite[Corollary 5.2.7]{Gehrunger25} implies that $\Jmodmstw$ has toric rank equal to $\frac{r-1}{2}$ (which is the maximal value), so $\Jmstw$ has totally toric reduction at $\frakq$. 
	
	Assume now that $\vq(\deltaK) = 1$, so that $\Kq(\sqrt{\deltaK}) / \Kq$ is ramified. The curve $\Cmstw$ is the quadratic twist of $\Cms$ by $\deltaK$, so these curves are isomorphic over $\Kq(\sqrt{\deltaK})$. Now $\Cms$ is obtained by setting $\deltaK = 1$, so the discussion above shows that $\Jmodms$ has toric reduction at $\frakq$, and thus $\Jmodmstw$ attains toric reduction over $\Kq(\sqrt{\deltaK})$. Now fix an odd prime number $\ell$ and write $\chi_{\deltaK}^{ }$ for the character of $G_{\Kq}$ cut out by $\Kq(\sqrt{\deltaK})$. By Grothendieck's inertial criterion (see \cite[Exposé IX]{SGA7}), the action of $I_{\frakq}$ on $\Vl(\Jms)$ is unipotent. But the action of $I_{\frakq}$ on $\Vl(\Jmstw)$ is given by the action above twisted by $\chi_{\deltaK}^{ }$. The character $\chi_{\deltaK}$ being ramified, $I_{\frakq}$ does not act unipotently on $\Vl(\Jmstw)$. Applying again Grothendieck's inertial criterion (this time in the reverse direction), we conclude that $\Jmodmstw / \Oq$ is not semistable.
\end{proof}

We now turn to the study of the Jacobian $\Jpltw$. 

\begin{proposition}\label{prop: Toric red Jmodpltw}
	Assume that $\nu_2 > 0$. If $\vq(\deltaK) = 0$, then $\Jmodpltw$ has toric reduction at $\frakq$. Otherwise, it has unipotent reduction at $\frakq$, and attains toric reduction over $\Kq(\sqrt{\deltaK})$.
\end{proposition}

\begin{proof}
	Assume first that $\vq(\deltaK) = 0$. Let us recall from Theorem~\ref{thm: Cluster Cpl} the description of the cluster picture of $\Cpl$ at $\frakq$. If $\vq(s_0-1) > 0$, then the cluster picture contains $\frac{r+1}{2}$ twins, among whose $\frac{r-1}{2}$ have depth $> 4$ and the one containing $\gamma_r, \gamma_{i_0}$ has depth $2$. Since $\vq(\deltaK) = 0$, the defining polynomial is unitary and the outer depth is $0$. Therefore, $\Cpl$ has $\frac{r-1}{2}$ grounded double points with thickness $> 4$, and so Corollary 5.2.7 in \cite{Gehrunger25} implies that $\Jpltw$ has (totally) toric reduction at $\frakq$. On the other hand, if $\vq(s_0) > 0$, then there are $\frac{r+1}{2}$ twins with depth $> 2$. In this case, we use Proposition 1.5 in \cite{DokchitserMorgan23}, which states that the special fiber of the stable model of $\Cpl$ is a union of two rational curves intersecting at $\frac{r+1}{2}$ points. It follows that the toric rank is maximal, and therefore $\Jpltw$ has totally toric reduction. 
	
	If $\vq(\deltaK) = 1$, the same discussion as in the previous proof shows that $\Jpltw$ attains toric reduction over $\Kq(\sqrt{\deltaK})$, and that it has unipotent reduction over $\Kq$. 
\end{proof}

\begin{remark}
	Above, we made use of the theory of cluster pictures in even residue characteristic to describe the reduction types of $\Jmodpmtw$. This is a truly special situation, as the cluster pictures of $\Cpm$ have many symmetries: Dokchitser--Morgan \cite{DokchitserMorgan23} and Gehrunger \cite{Gehrunger25} then explain how to extract arithmetic information from them. It does not seem likely that these arguments may adapted for other cluster pictures at even places. The study of the reduction of hyperelliptic curves in even residue characteristic is very rich, and one needs more complex and advanced combinatorial objects. We refer the curious reader to \cite{GehrungerPink} and \cite{Gehrunger25} for further details.
\end{remark}

\subsubsection{The case of (potential) good reduction}

We focus now on the case where $\Jpmtw$ has potential good reduction. We are going to manipulate different models of $\Cpmtw$ to describe the reduction type of $\Jmodpmtw$. We begin by introducing some notation to simplify the discussion below.

\begin{definition}\label{def: delta mu}
	Let $\pi$ be a uniformizer of $\Oq$. We define $\mu \coloneqq \vq(\deltaQ^r (2-4s_0)) \in \Z$, $\kappa \coloneqq \frac{\deltaQ^r (2-4s_0)}{\pi^{\mu}}$, and $\delta_2 \coloneqq \frac{\deltaK}{\pi^{\vq(\deltaK)}}$. 
\end{definition}

Recall that, given an element $x \in \Oq \st$, we say that the property $\SQ{x}$ holds if $\vq(x)$ is even, and $x$ is a square $\mod \frakq^2$, \ie there is some $\tau \in \Oq$ such that $\vq(x - \tau^2) \geq 2$.

\begin{remark}\label{rmk: Valuation Delta Wpm}
	If $\nu_2 < -2$, then the equality $(\deltaQ^{r} (2-4s_0))^2 = 4 \deltaQ^{2r} (4s_0 (s_0-1) +1 )$ implies that $\mu = \frac{1}{2} \vq(2^4 \deltaQ^{2r} s_0 (s_0-1)) \leq 0$. With the notation above, we can rewrite the valuations of the discriminants given in Corollary~\ref{cor: Discriminants Wpm} as
	\begin{align*}
		\vq(\Delta(\Wmstw)) & = (r-1) \left(2 + \mu \right) + 2r \vq(\deltaK), \\ 
		\vq(\Delta(\Wpltw)) & = 2(r-1) + \mu (r+1) + 2r \vq(\deltaK).
	\end{align*}
	We note that $\nu_2 < -2$ is equivalent to $\nu_2 \leq -4$, thanks to Hypothesis~\ref{hyp: 2adic val not -1 -3}.	
\end{remark}

We focus from now on the case $\nu_2 \leq -4$, as in this case the discriminants above can be expressed in terms of the valuation $\mu$ of the constant term of $\gminus$. Our analysis of the reduction of the Néron models $\Jmodpmtw$ at $\frakq$ is not complete, but it suffices for our Diophantine purposes. 

\begin{definition}
	Until the end of this subsection, we let $L$ be a finite extension of $\Kq$ with ramification index $r$. We let $\frakP$ be the maximal ideal of its ring of integers $\OP$, $\vP$ a valuation on $\overline{L}$ normalised with respect to $L$. We let $\pi_L \in \OP$ be such that $\pi_L^r = \pi$, so that $\pi_L$ is a uniformizer of $\OP$.
\end{definition}

\begin{proposition}\label{prop: Pot good Cmstw}
	Assume that $\nu_2 \leq -8$. If $\SQ{\deltaK \deltaQ^r (2-4s_0)}$ holds, then $\Cmstw$ attains good reduction over $L$. Otherwise, it only attains good reduction over a quadratic ramified extension of $L$.
\end{proposition}

\begin{proof}
	Assume first that $\SQ{\deltaK \deltaQ^r (2-4s_0)}$ holds. Then $\mu + \vq(\deltaK) = \vq(\deltaK \deltaQ^r (2 - 4s_0))$ is even. Moreover, $\delta_2 \kappa$ is a square modulo $\frakq^2$, so there is some $\tau \in \Oq \st$ such that $\vq(\delta_2 \kappa - \tau^2 ) \geq 2$. Consider the model $\Gmodms$ of the base change $\Cmstw \times_{\Kq} L$ obtained from $\Wmstw \times_{\Kq} L$ by setting the change of variables
	\reqnomode
	\begin{equation}\label{eq: Change vars Pmodms}
		x = \pi_{L}^{\mu + 2} X, \qquad \text{ and } \qquad y = \pi_{L}^{r (\mu + \vq(\deltaK))/2}  ( \pi_L^{r} Y + \tau).
	\end{equation}
	Recall that we describe $\Wms$ by \eqref{eq: expanded Cms}, and that $\Wmstw$ is given by the latter equation with the RHS multiplied by $\deltaK$. With the notation introduced in Definition~\ref{def: delta mu} and after rearranging terms, we describe $\Gmodms$ by
	\leqnomode
	\begin{equation}\label{eq: Pmodms}\tag{$\Gmodms$}
		Y^2 + \frac{2 \tau}{\pi_L^r} \, Y = \delta_2 \sum_{k = 0}^{\frac{r-1}{2}} (-1)^k \, c_k \, \left( \frac{\deltaQ}{\pi_{L}^{\mu + 2} }\right)^{2k} \, X^{r-2k} + \frac{\delta_2 \kappa - \tau^2}{\pi_{L}^{2r}}.
	\end{equation}
	We have $2 \tau / \pi_{L}^{r} \in \OP \st$. Since $\SQ{\deltaK \deltaQ^r (2-4s_0)}$ holds, the constant term belongs to $\Oq$. Moreover,
	\begin{equation*}
		\vP \left( \deltaQ / \pi_{L}^{\mu +2}\right) = r \vq(\deltaQ) -  (\vq(\deltaQ^r(2-4s_0)) + 2) = - (\vq(2-4s_0) + 2), 
	\end{equation*}
	and the latter is equal to $- \frac{\nu_2}{2} - 4$, which is non-negative as we assume that $\nu_2 \leq -8$. We deduce that $\Gmodms / \OP$ describes an integral model of $\Cmstw \times_{\Kq} L$. Recall that the curve $\Cmstw$ has genus $\frac{r-1}{2}$, so combining Lemma~\ref{lem: Legal change of vars} and Remark~\ref{rmk: Valuation Delta Wpm}, we compute
	\begin{align*}
		\vP(\Delta(\Gmodms)) & = r \vq(\Delta(\Wmstw)) - 4r \left( \frac{r (\mu + \vq(\deltaK))}{2} + r \right) + r(r+1)(\mu + 2) \\
		& = r (r-1)(2+\mu) + 2r^2 \vq(\deltaK) - r (r-1)(2 + \mu) - 2r^2 \vq(\deltaK) = 0.
	\end{align*}
	We conclude that $\Gmodms$ has good reduction over $\OP$, so $\Cmstw$ attains good reduction over $L$.
	
	On the other hand, if $\SQ{\deltaK \deltaQ^r (2-4s_0)}$ does not hold, then $\Jmodmstw$, the curve attains good reduction over a quadratic extension of $L$. A similar argument as in the proof of Proposition~\ref{prop: Toric red Jmodmstw}, using the action of inertia, shows that $\Jmodmstw \times \O_L$ has unipotent reduction.
\end{proof}

\begin{proposition}\label{prop: Pot good Cpltw}
	Assume that $\nu_2 \leq -4$. If $\SQ{\deltaK}$ holds, then $\Cpltw$ attains good reduction over $L$. Otherwise, it only attains good reduction over a quadratic ramified extension of $L$.
\end{proposition}

\begin{proof}
	To assume that $\SQ{\deltaK}$ holds means that $\deltaK = \delta_2 \in \Oq\st$, and that there is some $\vartheta \in \Oq\st$ such that $\vq(\deltaK - \vartheta^2) \geq 2$. Consider the model $\Gmodpl$ of the base change $\Cpltw \times_{\Kq} L$ obtained from $\Wpltw \times_{\Kq} L$ by setting the change of variables
	\reqnomode
	\begin{equation}\label{eq: Change vars Pmodpl}
		x = \frac{\pi_{L}^{\mu - 2}}{X}, \qquad \text{ and } \qquad y =  \pi_{L}^{(\mu -2) \left(\frac{r+1}{2} \right)} \left( \frac{\pi_{L}^r Y + \vartheta}{X^{\frac{r+1}{2}}} \right).
	\end{equation}
	The model $\Wpltw$ is described by \eqref{eq: expanded Cpl}, with the RHS multiplied by $\deltaK$. Using the notation introduced in Definition~\ref{def: delta mu} and rearranging terms, we describe $\Gmodpl$ by
	\leqnomode
	\begin{equation}\label{eq: Gmodpl}\tag{$\Gmodpl$}
		(\pi_{L}^{r} Y + \vartheta)^2 = \delta_2 \left( \sum_{k = 0}^{\frac{r-1}{2}} (-1)^k \, c_k \, \left( \frac{\deltaQ}{\pi_{L}^{\mu - 2}} \right)^{2k} \, X^{2k} + \pi_{L}^{2r} \kappa X^r \right) \left( 1 + \frac{2 \deltaQ}{\pi_{L}^{\mu - 2}} \, X \right).
	\end{equation}
	But $\gplus$ is monic, so the constant term on the RHS is $\delta_2$. The equation above can be rewritten as
	\begin{equation*}\tag{$\Gmodpl$}
		Y^2 + \frac{2 \vartheta}{\pi_{L}^{2r}} \, Y = \frac{\delta_2 - \vartheta^2}{\pi_{L}^{2r}} + \frac{2 \deltaQ}{\pi_L^{\mu -2 + 2r}} +  \delta_2 \left( \sum_{k = 1}^{\frac{r-1}{2}} (-1)^k \, \frac{c_k}{\pi_{L}^{2r}} \, \left( \frac{\deltaQ}{\pi_{L}^{\mu-2}} \right)^{2k} \, X^{2k} + \kappa X^r \right) \left( 1 + \frac{2 \deltaQ}{\pi_{L}^{\mu-2}} \, X\right).
	\end{equation*}
	We have $2\vartheta / \pi_{L}^{2r} \in \OP\st$, and $\vP(\delta_2 - \vartheta^2) \geq 2r$ by definition of $\vartheta$. On top of that,
	\begin{equation*}
		\vP \left(\deltaQ / \pi_{L}^{\mu - 2} \right) = r \vq(\deltaQ) - (\vq(\deltaQ^r (2-4s_0)) - 2) = 2 - \vq(2-4s_0),
	\end{equation*}
	and the latter equals $- \nu_2/2$, which is positive as $\nu_2 \leq -4$. It follows that $\vP(2\deltaQ / \pi_{L}^{\mu - 2 +2r}) \geq 0$, and the terms inside the sum also belong to $\OP$, as the index $k$ satisfies $k \geq 1$. We deduce that $\Gmodpl$ is integral over $\OP$. Combining Lemma~\ref{lem: Legal change of vars} with Remark~\ref{rmk: Valuation Delta Wpm}, we compute
	\begin{align*}
		\vP(\Delta(\Gmodpl)) & = r \, \vq(\Delta(\Wpltw)) - 4r \left((\mu -2) \left( \frac{r+1}{2} \right) + r \right) + r(r+1)(\mu-2)  \\
		& = 2r(r-1) + \mu r(r+1) -\mu r(r+1) + 2 r(r+1) - 4r^2 = 0. 
	\end{align*}
	We conclude that $\Gmodpl$ has good reduction over $\OP$, so $\Cpltw$ attains good reduction over $L$.
	
	Again, if $\SQ{\deltaK}$ does not hold, then $\Cpltw$ attains good reduction over a quadratic ramified extension of $L$. The action of the inertia group $I_{\frakq}$ on the $\ell$-adic Tate module of $\Jpltw$ shows that $\Jmodpltw \times \O_L$ has unipotent reduction.
\end{proof}

When some extra congruences are satisfied, the curves $\Cpmtw$ have good reduction over $\Kq$. 

\begin{corollary}\label{cor: Good red Jmstw}
	Assume that $\nu_2 \leq -8$ and $\SQ{\deltaK \deltaQ^r (2 - 4s_0)}$ holds. If $\nu_2 \equiv -8 \mod r$, then $\Cmstw$ has good reduction over $\Kq$, and otherwise it has bad reduction.
\end{corollary}

\begin{proof}
	Since $\nu_2 \leq -8 < 0$, the assumption $\nu_2 \equiv -8 \mod r$ is equivalent to $\mu \equiv -2 \mod r$.  
	If this is the case, the model \eqref{eq: Pmodms} is defined over $\Kq$, and 	the change of variables given in \eqref{eq: Change vars Pmodms} too, so $\Cmstw$ has good reduction at $\frakq$. 
	If $\nu_2 \not \equiv -8 \mod r$, then $\vq(\Delta(\Wmstw)) \not \equiv 0 \mod 2r$ (\cf \ Remark~\ref{rmk: Valuation Delta Wpm}). Since $\Cmstw$ has genus $\frac{r-1}{2}$, Remark~\ref{rmk: Bad red congr} implies that $\Cmstw$ has bad reduction at $\frakq$.
\end{proof}

\begin{corollary}\label{cor: Good red Jpltw}
	Assume that $\nu_2 \leq -4$ and $\SQ{\deltaK}$ holds. If $\nu_2 \equiv 0 \mod r$, then $\Cpltw$ has good reduction over $\Kq$, and otherwise it has bad reduction.
\end{corollary}

\begin{proof}
	Since $\nu_2 \leq -4$, we have $\nu_2 \equiv 0 \mod r$ if and only if $\mu \equiv 2 \mod r$. If the latter holds, the model $\Gmodpl$ and the change of variables given in \eqref{eq: Change vars Pmodpl} are defined over $\Kq$. Thus, $\Cpltw$ has good reduction at $\frakq$. If $\nu_2 \not \equiv 0 \mod r$, then $\vq(\Delta(\Wpltw)) \not \equiv 0 \mod 2r$ by Remark~\ref{rmk: Valuation Delta Wpm}. The curve $\Cpltw$ also has genus $\frac{r-1}{2}$, so Remark~\ref{rmk: Bad red congr} shows that $\Cpltw$ has bad reduction at $\frakq$.
\end{proof}

\begin{corollary}\label{cor: Semist defect q even}
	If $\nu_2 \leq -8$, $\nu_2 \not \equiv -8 \mod r$ and $\SQ{\deltaK \deltaQ^r (2-4s_0)}$ holds, then the semistability defect of $\Jmstw / \Kq$ equals $r$. Similarly, if $\nu_2 \leq -4$, $\nu_2 \not \equiv 0 \mod r$ and $\SQ{\deltaK}$ holds, then the semistability defect of $\Jpltw / \Kq$ equals $r$.
\end{corollary}

\begin{remark}
	The equality $\pi_L^r = \pi$ shows that the models $\Gmodpm$ are both defined over $\Kq$. However, the change of variables relating the Weierstrass models $\Wpmtw$ to $\Gmodpm$ is only defined over $\Kq$ when the congruences displayed in the two corollaries above are satisfied.
\end{remark}

\begin{remark}\label{rmk: Remaining cases Red q even}
	Our analysis of the reduction type of $\Jmodpmtw$ at an even place $\frakq$ is not complete. Numerical computations lead the author to believe that, when $-8 < \nu_2 \leq 0$ (resp. $-4 < \nu_2 \leq 0$), then $\Jmodmstw$ (resp. $\Jmodpltw$) seems to have unipotent reduction at $\frakq$, and to attain good reduction over a wildly ramified extension of $\Kq$. 
\end{remark}

\subsection{Ramification indices and discriminants}\label{sect: Ramification}

Now that we understand the reduction types of $\Jmodpmtw$ at even places, we study these Néron models at odd places of $\Kgl$. For this, we make use of the theory of cluster pictures \cite{M2D2} in the style of \cite{ACIKMM}. in order to prove Theorems~\ref{thm: reduction types Jms} and \ref{thm: reduction types Jpl}, we need to understand the ramification indices of the local extensions obtained by adjoining the roots of the defining polynomials $\gpm$. 
\medskip 

We devote this subsection to the computation of such ramification indices, and we study the splitting field of the polynomial $\gminus$, which we view as an element of $\Qq[x]$. As explained in Definition~\ref{def: Curve Cpm}, we only consider the curve $\Cpl$ when the parameter $\deltaQ$ belongs to $\Z$ (and thus to $\Zq$). In that case, $\gamma_r = -2\deltaQ$ is a $\Qq$ rational root, so $\Qqroots$ is also the splitting field of $\gplus$. 
\medskip 

We aim at understanding the base-changed curve $\Cpm / \Kq$, so it could be reasonable to manipulate the polynomial $\gminus$ as an element of $\Kq[x]$. However, for computational reasons, we prefer to consider that $\gminus$ belongs to $\Qq[x]$, and then deduce the desired results about $\Cpm / \Kq$.

\subsubsection{Description of the splitting field}

Recall that, for any $j \in \Iintv{0, r-1}$, we defined $\omega_j = \zr^j + \zr^{-j}$, $\tau_j \coloneqq \zr^j - \zr^{-j}$, $\omega = \omega_1$, and $\tau = \tau_1$. We view all these as elements of $\Qqbar$, and we identify $\Qq(\omega) \simeq \Kq$. Recall that we denote by $\phi_r(X, Y)$ the polynomial given by
\begin{equation*}
	\phi_r(X, Y) \coloneqq (X^r + Y^r) / (X + Y) = \sum_{j = 0}^{r-1} (-1)^j X^{r-1-j} Y^j.
\end{equation*}

\begin{definition}
	We let $\Qqroots$ be the splitting field of the polynomial $\gminus(x) \in \Qq[x]$. We let $\Kqroots$ be the compositum of $\Qqroots$ and $\Kq$. 
\end{definition}

\begin{theorem}\label{thm: Splitting field}
	The splitting field of $\gminus$ is $\Qqroots = \Kq(\gamma_0, \deltaQ \tau \sqt)$.
\end{theorem}

\begin{proof}
	Let us begin with the reverse inclusion. We clearly have $\gamma_0 \in \Qqroots$, and the equality $\omega = (\gamma_1 + \gamma_{-1})/ \gamma_0$ shows that $\Kq \subset \Qqroots$. It remains to prove that $\deltaQ \tau \sqt \in \Qqroots$. Recall from Lemma~\ref{lem: useful eqs} \eqref{enum: useful eqs3} that $\alpha_0^r - \beta_0^r = 4 \sqt$. One can check that, for any $j \in \Iintv{0, r-1}$, we have $\deltaQ \tau_j (\alpha_0 - \beta_0) = \gamma_j - \gamma_{-j}  \in \Qqroots$. A direct application of \cite[Lemma 2.14]{ACIKMM} gives the equality 
	\begin{equation*}
		\phi_{r}(\alpha_0, - \beta_0) = h_r \left( \left(\frac{\gamma_0}{\deltaQ}\right)^2 - 2 \right).
	\end{equation*}
	But $h_r$ belongs to $\Z[x]$ and $\deltaQ^2 \in \Z$, so the latter quantity belongs to $\Qq(\gamma_0) \subset \Qqroots$. Finally, the equality 
	\begin{equation*}
		\deltaQ \tau \sqt = \frac{\deltaQ \tau}{4}  \, (\alpha_0^r - \beta_0^r) = \frac{\deltaQ \tau}{4} ( \alpha_0 - \beta_0) \, \phi_r(\alpha_0, - \beta_0)
	\end{equation*}
	gives the desired inclusion.

	Let us show now that $\Qqroots \subset \Kq(\gamma_0, \deltaQ \tau \sqt)$. Fix $j \in \Iintv{0, r-1}$. One can check that $\gamma_j = (\omega_j \gamma_0 + \deltaQ \tau_j(\alpha_0 - \beta_0))/2$. Moreover, $\deltaQ \tau_j (\alpha_0 - \beta_0) =  \frac{4 \deltaQ \tau_j}{\tau} \, \tau \sqt \ \phi_r(\alpha_0, - \beta_0)^{-1}$. But 
	\begin{equation*}
		\frac{\tau_j}{\tau} = \frac{\omega_{j + 2} - \omega_j }{\omega^2 - 4} \in \Qq(\omega) \simeq \Kq,
	\end{equation*}
	so the discussion above shows that $\phi_r(\alpha_0, - \beta_0) \in \Qq(\gamma_0)$, hence the result. 
\end{proof}

\begin{definition}
	We let $\Qone \coloneqq \Qq(\deltaQ^r \sqt), \Qquad \coloneqq \Kq \cdot \Qone$, and $\Qtau \coloneqq \Kq(\deltaQ^r \tau \sqt)$. We let $\frakQ$ be the maximal ideal in $\O_{\Qone}$, and $\vQ$ be a valuation on $\Qqbar$ normalised with respect to $\Qone$.
\end{definition}

\begin{remark}\label{rmk: Def Qquad Qtau}
	Theorem~\ref{thm: Splitting field} implies that $\Kq \subset \Qqroots$. Recall that $\deltaQ^{2r} s_0(s_0-1) \in \Q$, and we also have $\tau^2 = (\omega_2 - 2) \in \OK$. We deduce that $\Qone / \Qq$ is at most quadratic, and that $\Qquad$, $\Qtau$ are at most quadratic extensions of $\Kq$. With this notation, Theorem~\ref{thm: Splitting field} simply states that $\Qqroots = \Qtau(\gamma_0)$.
\end{remark}

\begin{proposition}\label{prop: gamma0 rational}
	If $\gminus$ is reducible over $\Qq$ then $\gamma_{0} \in \Qone$. If we further assume $\zr \notin \Qq$, then $\gamma_0 \in \Qq$.
\end{proposition}

\begin{proof}
	Consider the diagram of field extensions depicted in Figure~\ref{fig: Diamond diagram}.
	\begin{figure}[h!]	
		\captionsetup{justification=centering,margin=1.5cm}
		\begin{center}
			\begin{tikzcd}
			& \Qq(\deltaQ \alpha_0)                                      &                                            \\
			\Qone \arrow[ru, "n_0", no head] &                                                             & \Qq(\gamma_0) \arrow[lu, "d_0'"', no head] \\
			& \Qq \arrow[ru, "n_0'"', no head] \arrow[lu, "d_0", no head] &                                           
			\end{tikzcd}
			\caption{\centering Diagram of field extensions relating $\Qone$ and $\Qq(\gamma_0)$. The label on each line denotes the degree of the respective extension.}
			\label{fig: Diamond diagram}
		\end{center}
	\end{figure} 

	As we saw just above, $d_0 \leq 2$. Similarly, we have $\beta_0 = \alpha_0^{-1}$, so $\gamma_0 = \deltaQ\alpha_0 + \deltaQ^2 / (\deltaQ \alpha_0)$. Therefore, the polynomial $X^2 - \gamma_0 X + \deltaQ^2 \in \Qq(\gamma_0)[X]$ vanishes at $\deltaQ \alpha_0$, showing that $d_0' \leq 2$. Finally, the extension $\Qq(\deltaQ \alpha_0) / \Qone$ is defined by the polynomial $x^r - \deltaQ^r (\sqts + \sqtsm)^2$, so Remark~\ref{rmk: alpha0 rational} implies that $n_0 \in \lbrace 1, r \rbrace$. If $\gminus$ is reducible over $\Qq$, then $n_0' < r$, and thus $n_0 d_0 = n_0' d_0' < 2r$. It follows that $n_0 < r$, and therefore $\Qq(\deltaQ \alpha_0) = \Qone$, hence $\gamma_0 \in \Qone$, so the first statement is satisfied. 
	
	If $\Qone = \Qq$, the second statement holds too. Assume therefore that $\zr \notin \Qq$ and that $\Qq \neq \Qone$. Let $\sigma$ (resp. $\Nrm{\Qone}{\Qq}$, $\Tr{\Qone}{\Qq}$) be the non-trivial automorphism (resp. the norm and trace maps) of the quadratic extension $\Qone / \Qq$. From Definition~\ref{def: alpha0 beta0} we obtain the identities
	\begin{equation*}
		(\deltaQ \alpha_0 )^r = \deltaQ^r \left(2 s_0 -1 \right) + 2 \deltaQ^r \sqt \quad \text{ and, } \quad (\deltaQ \beta_0 )^r = \deltaQ^r \left(2 s_0 -1 \right) - 2 \deltaQ^r \sqt.
	\end{equation*}
	But we know from Definition \ref{def: s0 deltaQ} that $\deltaQ^r (2s_0-1) \in \Qq$, so $\sigma(\deltaQ \alpha_0)^r = (\deltaQ \beta_0)^r = \deltaQ^{2r} / (\deltaQ \alpha_0)^r$. Consequently $\Nrm{\Qone}{\Qq}(\deltaQ \alpha_0)^r = \deltaQ^{2r}$, and we deduce that $\Nrm{\Qone}{\Qq}(\deltaQ \alpha_0) = \deltaQ^{2}$, as $\zr \notin \Qq$. Therefore $\sigma(\deltaQ \alpha_0) = \deltaQ \beta_0$, and we conclude that $\gamma_0 = \Tr{\Qone}{\Qq}(\deltaQ \alpha_0)$ indeed belongs to $\Qq$.
\end{proof}

\begin{remark}\label{rmk: Ram ind quadratics}
	Proposition~\ref{prop: gamma0 rational} and Figure~\ref{fig: Diamond diagram} show that the degrees $n_0$, $n_0'$ are equal and belong to $\lbrace 1, r \rbrace$. Similarly, $d_0 = d_0' \in \lbrace 1, 2 \rbrace$. By coprimality of degrees, we also obtain the equalities of ramification indices $e_{\Qq(\gamma_0) / \Qq}^{ } = e_{\Qq(\deltaQ \alpha_0) / \Qone}^{ }$. Similarly, one can show that the extensions $\Kq(\gamma_0) /  \Kq$ and $\Kq(\deltaQ \alpha_0) / \Qquad$ have same degree and same ramification index.
\end{remark}

\subsubsection{Reducibility criteria for the polynomial $\gminus(x)$}
	
We now give reducibility criteria concerning the polynomial $\gminus(x) \in \Qq[x]$.

\begin{remark}\label{rmk: gminus irred Qq Kq}
	Since $\gminus$ has degree $r$ and $[\Kq : \Q_q]$ has degree at most $\frac{r-1}{2}$, it is easy to check that $\gminus$ is irreducible over $\Qq$ if and only if it is irreducible over $\Kq$. Moreover, we have $e_{\Qr(\gamma_0) / \Qr} = e_{\Kr(\gamma_0) / \Kr}$.
\end{remark}

\begin{proposition}\label{prop: gminus reducible}
	If $\, \frac{1}{2} \vq(s_0(s_0-1)) > r \vq(1 - \zr)$, then $\gminus$ is reducible over $\Qq$ and $\gamma_{i_0} \in \Qq$.
\end{proposition}

\begin{proof}
	As discussed at the beginning of \S \ref{sect: Clusters Cpm}, the cluster picture of $\CmssQ / \Qq$ is the one of $\Cms / \Kq$, except that the depths are multiplied by the ramification index of $\Kq / \Qq$. We know that the absolute Galois group $G_{\Qq}$ acts on clusters, preserving depths and containments (\cite[Remark 3.2]{Hyperuser}). When $\vq(s_0(s_0-1)) > r \vq(1 - \zr)$ (whether $\frakq = \frakr$ or not), Theorem~\ref{thm: Cluster Cms} implies that $\gamma_{i_0}$ is an isolated root in the cluster picture of $\CmssQ / \Qq$, hence fixed by all the elements of $G_{\Qq}$.
\end{proof}

\begin{proposition}\label{prop: gminus irreducible}
	If $v_q(s_0) \leq 0$ and $v_q(s_0) \not \equiv 0 \mod r$, then $\gminus$ is irreducible over $\Qq$.
\end{proposition}

\begin{proof}
	By definition we have $\gminus(x) = \deltaQ^r g_{r, s_0}^{-}(x / \deltaQ)$ (see \eqref{eq: pols gminus gplus}), so $\gminus$ is irreducible over $\Qq$ if and only if $g_{r, s_0}^{-}$ is so. The latter is given by $g_{r, s_0}^{-}(x) = \sum_{k = 0}^{(r-1)/2} (-1)^k c_k x^{r-2k} + 2 - 4s_0$. Since $c_k \in \Z$, assuming that $\vq(s_0) \leq 0$ implies that the Newton polygon of $g_{r, s_0}^{-}$ is made of a single line joining the points of coordinates $(0, \vq(s_0))$ and $(r, 0)$. The slope of such line is $\vq(s_0) / r$, so if $v_q(s_0) \not \equiv 0 \mod r$, then $g_{r, s_0}^{-}$ is irreducible over $\Qq$ (see \cite[Chapter II \S 6]{Neukirch99}). 
\end{proof}

\begin{proposition}\label{prop: irred when nur 1 or 2}
	If $q = r$, and $\nu_r = v_r(s_0(s_0-1)) \in \lbrace 1, 2 \rbrace$, then $\gminus$ is irreducible over $\Qr$.
\end{proposition}

\begin{proof}
	Recall that we denote by $\frakQ$ the maximal ideal in $\O_{\Qone}$. Until the end of the proof, write $f \coloneqq f_{\Qone / \Qr}$ for the residue degree of $\Qone / \Qr$, and $u \coloneqq (\deltaQ \alpha_0)^r = \deltaQ^r \left(2 s_0 -1 \right) + 2 \deltaQ^r \sqt \in \Qone$. Remark~\ref{rmk: Ram ind quadratics} shows that $\gminus$ is irreducible over $\Qr$ if and only if $x^r - u$ is irreducible over $\Qone$, which occurs if and only if \vspace{-0.5em}
	\begin{equation*}
		E(x) \coloneqq \left(x + u^{r^{f -1}} \right)^r - u \, = \, x^r + \sum_{j = 1}^{r-1} \binom{r}{j} (u^j)^{r^{f-1}} x^{r- j} + u^{r^f} - u \in \Qone[x] \vspace{-0.25em}
	\end{equation*}
	is irreducible over $\Qone$. We claim that $\nu_r \in \lbrace 1, 2 \rbrace$ implies that $E(x)$ is Eisenstein. The coefficients in the sum are divisible by $r$, so their $\frakQ$-adic valuation is $\geq 1$. By Fermat's little theorem, $\vQ(u^{r^f} - u) \geq 1$, and we are going to prove that the latter is an equality.
	
	Note first that the quadratic extension $\Qone / \Qr$ is ramified if and only if $v_r(\deltaQ^{2r} s_0 (s_0-1))$ is odd (see Remark~\ref{rmk: ramif Qone Qr} below). The assumption $\nu_r \in \lbrace 1, 2 \rbrace$ implies then $\vQ(\sqt) = 1$, and Hypothesis~\ref{hyp: v(s0(s0-1)) = 0} yields $v_r(\deltaQ) = 0$, so $\vQ(\deltaQ^r \sqt) = 1$ too. By Newton's binomial formula, we have \vspace{-0.5em} 
	\begin{equation*}
		u^{r^f} = \left( \deltaQ^r \left(2 s_0 -1 \right) + 2 \deltaQ^r \sqt \right)^{r^f} \equiv \left( \deltaQ^r \left(2 s_0 -1 \right) \right)^{r^f} \mod \frakQ^2.
	\end{equation*}
	Recall from Lemma~\ref{lem: s0(s0-1) in Q} that $\deltaQ^r (2s_0 - 1) \in \Zr$, so when $\nu_r = v_r(s_0 (s_0 -1))$ equals $1$, we have
	\begin{equation*}
		\vQ \left( (\deltaQ^r (2s_0 -1))^{r^f} - \deltaQ^r (2s_0 -1) \right) = 2 v_r \left( (\deltaQ^r (2s_0 -1))^{r^f} - \deltaQ^r (2s_0 -1) \right) \geq 2.
	\end{equation*}
	If $\nu_r = 2$, then $(\deltaQ^r (2s_0 -1))^2 = 4 \deltaQ^{2r} s_0 (s_0 -1) + \deltaQ^{2r}$ is congruent to $\deltaQ^{2r} \mod r^2$, as $r^2 \mid s_0 (s_0 -1)$. From this, we deduce $(\deltaQ^r (2s_0 -1))^{2r} \equiv \deltaQ^{2r^2} \equiv \deltaQ^{2r} \mod r^2$, as $\Z / r^2 \Z$ has $r(r-1)$ invertible elements. It follows easily that $\vQ((\deltaQ^r (2s_0 -1))^{r^f} - \deltaQ^r (2s_0 -1)) \geq 2$ holds too. Whether $\nu_r = 1$ or $2$, we obtain
	\begin{equation*}
		u^{r^f} - u \equiv \left( \deltaQ^r \left(2 s_0 -1 \right) \right)^{r^f} - \deltaQ^r (2s_0 - 1) - 2 \deltaQ^r \sqt \equiv - 2 \deltaQ^r \sqt \mod \frakQ^2.
	\end{equation*}
	The last term is not congruent to $0 \mod \frakQ^2$, because $\vQ(\deltaQ^r \sqt) = 1$, as we saw above. We conclude that $E(x)$ is Eisenstein, hence irreducible over $\Qone$, so $\gminus$ is irreducible over $\Qr$.
\end{proof}

\subsubsection{Ramification indices of the involved extensions}

As explained in Remark~\ref{rmk: Def Qquad Qtau}, $\Kq$ is a subfield of $\Qqroots$. If one views $\gminus$ as an element of $\Kq[x]$, then its splitting field is $\Kqroots$. Therefore, when considering the base-changed curves $\Cpm / \Kq$, it is important to know the ramification index of $\Kqroots / \Kq$. As we saw above, $\Kqroots$ fits in the tower of extensions $\Kq \subset \Qtau \subset \Kqroots$. We now describe the ramification indices of the intermediate extensions. 

\begin{remark}\label{rmk: nur leq 0}
	Lemma~\ref{lem: s0(s0-1) in Q} shows that $s_0 (s_0 - 1) \in \Q$, and since $\Kr / \Qr$ is totally ramified of degree $\frac{r-1}{2}$, we have $\vr(s_0 (s_0 -1)) = \frac{r-1}{2} \, \nu_r \in \frac{r-1}{2} \, \Z$. Thus, $\nu_r$ is divisible by $r$ if and only if $\vr(s_0 (s_0 -1))$ is divisible by $r$. If $\vr(s_0 (s_0 -1)) \leq r$, the congruence $\vr(s_0 (s_0 -1)) \equiv 0 \mod r$ implies that $\nu_r \leq 0$.
\end{remark}

\begin{lemma}\label{lem: Ramification quadratic}
	\begin{enumerate}[leftmargin=*, itemsep=0pt]
		\item If $\frakq \neq \frakr$, then $\Qtau / \Kq$ is ramified  if and only if $\vq(\deltaQ^{2r} s_0 (s_0-1))$ is odd.
		
		\item If $\frakq = \frakr$, then $\Qtau / \Kr$ is ramified if and only if $\vr(\deltaQ^{2r} s_0 (s_0-1))$ is even. In particular, if $\vr(s_0(s_0-1)) \leq r$ and $\vr(s_0(s_0-1)) \equiv 0 \mod r$, then $\Qtau / \Kr$ is ramified.
	\end{enumerate}
\end{lemma}

\begin{proof}
	The quadratic extension $\Qtau / \Kq$ is defined by $x^2 - (\omega_2 - 2)\deltaQ^{2r} s_0(s_0-1)$. It is ramified if and only if the discriminant of the polynomial above has odd valuation, \ie $\vq((\omega_2 - 2) \deltaQ^{2r} s_0(s_0-1))$ is odd. If $\frakq \neq \frakr$, we have $\vq(\omega_2 -2) = 0$, and if $\frakq = \frakr$, $\vr(\omega_2 -2) = 1$, so the result follows. In particular, if $\frakq = \frakr$, $\vr(s_0(s_0-1)) \leq r$ and $\vr(s_0(s_0-1)) \equiv 0 \mod r$ then Hypothesis~\ref{hyp: val equiv 0 mod r} implies that $\vr(\deltaQ^{2r} s_0 (s_0-1)) = 0$, and we deduce the last claim.
\end{proof}

\begin{remark}\label{rmk: ramif Qone Qr}
	Similarly, when $q = r$, one can check that the quadratic extensions $\Qone / \Qr$ and $\Qr(\deltaQ \alpha_0) / \Qr(\gamma_0)$ are ramified (\cf \ Remark~\ref{rmk: Ram ind quadratics}) if and only if $v_r(\deltaQ^{2r} s_0(s_0-1))$ is odd.
\end{remark}

\begin{theorem}\label{thm: ramif index}
	Assume that $\gminus$ is irreducible over $\Qq$. 
	\begin{enumerate}[leftmargin=*, itemsep=0pt]
		\item If $\frakq \neq \frakr$, then $\Kq / \Qtau$ is (totally) ramified (with ramification degree $r$) if and only if $\vq(s_0(s_0-1)) \not \equiv 0 \mod r$.
		
		\item If $\frakq = \frakr$, then $\Kr / \Qtau$ is (totally) ramified with ramification degree $r$.
	\end{enumerate}
\end{theorem}

\begin{proof}
	Recall that we defined $\Qquad$ and $\Qtau$ as $\Qquad = \Kq(\deltaQ \sqt)$ and $\Qtau = \Kq(\deltaQ \tau \sqt)$, so $\Kq(\zr, \deltaQ \sqt) = \Qquad(\zr) = \Qtau(\zr)$. We describe the different field inclusions in Figure~\ref{fig: Ramif Diagram}.
	
	As discussed in Remark~\ref{rmk: Ram ind quadratics}, the extension $\Kq(\gamma_0) / \Kq$ has degree $1$ or $r$, and is non-trivial if and only if $\Kq(\deltaQ \alpha_0) / \Qquad$ is non-trivial. If $\gminus$ is irreducible over $\Qq$ (or, equivalently, over $\Kq$), then the extensions $\Kq(\deltaQ \alpha_0) / \Qquad$, $\Kq(\zr, \, \deltaQ \alpha_0) / \Qquad(\zr)$, and $\Kqroots / \Qtau$ are all non-trivial, hence of degree $r$.
	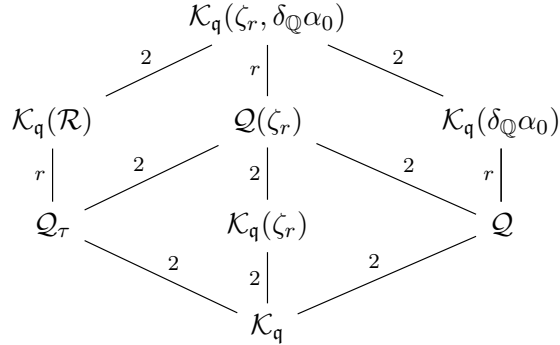
\begin{figure}[h]
		\captionsetup{justification=centering,margin=1cm}
		\begin{center}
			\begin{tikzcd}
				& {\Kq(\zr, \deltaQ \alpha_0)}                                                  &                                                                                           \\
				\Kqroots \arrow[ru, "2", no head]                      & {\Qquad(\zr)} \arrow[u, "r", no head]                            & \Kq(\deltaQ \alpha_0) \arrow[lu, "2"', no head]                                          \\
				\Qtau \arrow[u, "r", no head] \arrow[ru, "2", no head] & \Kq(\zr) \arrow[u, "2", no head]                                               & \Qquad \arrow[u, "r", no head] \arrow[lu, "2"', no head] \arrow[u, no head] \\
				& \Kq \arrow[lu, "2"', no head] \arrow[u, "2", no head] \arrow[ru, "2", no head] &                                                                                          
			\end{tikzcd}
			\caption{\centering Diagram of field inclusions when $\gminus$ is irreducible over $\Kq$. The degree of each extension divides the number appearing on the corresponding line.}
			\label{fig: Ramif Diagram}
		\end{center}
	\end{figure}	

	Studying the ramification indices that appear in Figure~\ref{fig: Ramif Diagram}, the coprimality of $2$ and $r$ gives the equalities 
	\begin{equation*}
		e_{\Kqroots / \Qtau}^{ } \, = \,   e_{\Kq(\zr, \, \deltaQ \alpha_0) / \Qquad(\zr)}^{ } \, = \, e_{\Kq(\deltaQ \alpha_0) / \Qquad}^{ }.  \vspace{0.5em}
	\end{equation*}
	By definition of $\alpha_0$, we have $(\deltaQ \alpha_0)^r = \deltaQ^r (\sqts + \sqtsm)^2$, so $\Kq(\zr, \deltaQ \alpha_0) / \Qquad(\zr)$ is a Kummer extension, obtained by adjoining to $\Qquad(\zr)$ the $r$-th root of $\deltaQ^r (\sqts + \sqtsm)^2$. Denote by $\frakP$ the maximal ideal of the ring of integers of $\Qquad(\zr)$. We treat different cases separately.
	\begin{enumerate}[leftmargin=*]
		\item If $\frakq \neq \frakr$, then Theorem 6.3 in \cite[I \S 6]{Gras03} states that the Kummer extension $\Kq(\zr, \deltaQ \alpha_0) / \Qquad(\zr)$ is ramified if and only if $\vP((\deltaQ \alpha_0)^r ) \not \equiv 0 \mod r$. But we assume that $\gminus$ is irreducible over $\Kq$, so Proposition~\ref{prop: gminus reducible} implies that $\vq(s_0 (s_0-1)) \leq 0$, hence $\vq(\alpha_0) = \pm \vq(s_0) / r$ by Lemma~\ref{lem: valuations alpha0 beta0}. Since $r$ is odd and $e_{\Qquad(\zr) / \Kq}$ is a power of $2$, we deduce that $\Kq(\zr, \deltaQ \alpha_0) / \Qquad(\zr)$ is ramified if and only if $\vq(\deltaQ^r s_0) \not \equiv 0 \mod r \Leftrightarrow \vq(\deltaQ^{2r} s_0(s_0-1)) \not \equiv 0 \mod r \Leftrightarrow \vq(s_0(s_0-1)) \not \equiv 0 \mod r$.
		
		\item If $\frakq = \frakr$, $\vr(s_0(s_0-1)) < 0$ and $\vr(s_0(s_0-1)) \not \equiv 0 \mod r$, then we have $\vP((\deltaQ \alpha_0)^r) \not \equiv 0 \mod r$, just as above. Therefore, \cite[I \S 6 Theorem 6.3]{Gras03} implies that $\Kq(\zr, \deltaQ \alpha_0) / \Qquad(\zr)$ is totally ramified, so we have again $e_{\Kqroots / \Qtau} = r$.
		
		\item If $\frakq = \frakr$ and $0 \leq \vr(s_0(s_0-1)) \leq r$, then $\vr(\alpha_0) = 0$ by Lemma~\ref{lem: valuations alpha0 beta0} \eqref{enum: val alpha0 beta 0 i}, and $\vr(\deltaQ) = 0$ by Hypothesis~\ref{hyp: v(s0(s0-1)) = 0}, so $\vr((\deltaQ \alpha_0)^r) = 0$. We claim that $\Kr(\deltaQ \alpha_0) / \Qquad$ is totally ramified. Indeed, it is defined by $x^r - (\deltaQ \alpha_0)^r$, and $\gminus$ being irreducible over $\Qr$ implies that $[\Kr(\deltaQ \alpha_0) : \Qquad] = r$. But this is not a Galois extension, as $\zr \notin \Qquad$, so it has to be ramified (\cite[Chapter III Theorem 2]{Serre79}). We conclude that $e_{\Kr(\deltaQ \alpha_0) / \Qquad} = e_{\Krroots / \Qtau} = r$.
		
		\item If $\vr(s_0(s_0-1)) < 0$ and $\vr(s_0(s_0-1)) \equiv 0 \mod r$, Hypothesis~\ref{hyp: val equiv 0 mod r} gives $\vr(\deltaQ^{2r} s_0 (s_0-1)) = 0$, so $\vr(\deltaQ^r) = - \vr(s_0) = - \vr(s_0 -1)$. On the other hand, Lemma~\ref{lem: valuations alpha0 beta0} \eqref{enum: val alpha0 beta 0 ii} gives $\vr(\alpha_0) = \pm \vr(s_0)/r$. Up to switching $\alpha_0$ and $\beta_0$, we may assume that $\vr(\alpha_0^r) = \vr(s_0)$ (this is without loss of generality because $\Kr(\deltaQ \alpha_0) = \Kr(\deltaQ \beta_0)$). We obtain $\vr(\deltaQ^r \alpha_0^r) = - \vr(s_0) + \vr(s_0) = 0$, and we conclude that $e_{\Kr(\deltaQ \alpha_0) / \Qquad} = e_{\Krroots / \Qtau} = r$ just as in the previous case.
	\end{enumerate}\vspace{-1.5em}
\end{proof}

\begin{proposition}\label{prop: Krroots ramified nur leq 2}
	If $\nu_r \leq 2$, the extension $\Krroots / \Kr$ is ramified.
\end{proposition}

\begin{proof}
	If $\nu_r \in \lbrace 1 , 2 \rbrace$, or if $\nu_r \leq 0$ and $\nu_r \not \equiv 0 \mod r$, then $\gminus$ is irreducible over $\Qr$ (Propositions~\ref{prop: gminus irreducible} and \ref{prop: irred when nur 1 or 2}). Theorem~\ref{thm: ramif index} implies then that $\Krroots / \Qtau$ is totally ramified. On the other hand, if $\nu_r \leq 0$ and $\nu_r \equiv 0 \mod r$, Lemma~\ref{lem: Ramification quadratic} states that $\Qtau / \Kr$ is ramified. In all cases we have $e_{\Krroots / \Kr} > 1$.
\end{proof}

\subsubsection{Discriminant of the totally ramified extension $\Qr(\gamma_0) / \Qr$}

Until the end of this subsubsection, we focus on the case $\frakq = \frakr$ (\ie $q = r$). Theorem~\ref{thm: ramif index} shows that, when $\gminus$ is irreducible over $\Qr$, then $\Krroots / \Kr$ is wildly ramified. As we will see in \S \ref{sect: Cond & local types}, the wild conductor at $\frakr$ of the $\ell$-adic representation attached to $\Jpm / \Kr$ is non trivial. Recall from Definition~\ref{def: Curve Cpm}that $\Jms / \Qr$ is the base change of $\Jac(\CmssQ)$ to $\Kr$. We will deduce the value of the wild conductor at $\frakr$ of $\rhomsK$ from the one of the $\ell$-adic representation attached to $\Jac(\CmssQ) / \Qr$. In order to describe the latter, we will use Theorem~\ref{thm: Wild cond clusters}, which requires to know the $r$-adic valuation of $\Delta(\Qr(\gamma_0) / \Qr)$.
\medskip 

We begin with a general statement, that we will apply later in our specific context.

\begin{proposition}\label{prop: Valuation discr}
	Let $F / \Q_r$ be a finite extension such that $\zr \notin F$, and let $\vF$ be a valuation on $\overline{F}$ normalised with respect to $F$. Let $u \in \OF \setminus \OF^{r}$, and let $M \coloneqq F(u^{1/r})$. Then
	\begin{equation*}
		\vF (\Delta(M / F)) = \begin{cases}
			r \, e_{F / \Qr} & \text{ if } \vF(u) \equiv 0 \mod r, \\
			r \, e_{F / \Qr} + r-1 & \text{ if } \vF(u) \not \equiv 0 \mod r. \\
		\end{cases}
	\end{equation*}
\end{proposition}

\begin{proof}
	Since $\zr \notin F$, the extension $M / F$ has degree $r$ and is not Galois, so it is totally ramified. Let $\pi_F$ be a uniformizer of $\O_F$ and $\pi_M \coloneqq \pi_F^{1/r}$, which is a uniformizer of $\O_M$. To simplify the notation, we denote by $f \coloneqq f_{F / \Qr}$ the residue degree of $F / \Qr$.
	\begin{enumerate}[leftmargin=*]
		\item Assume that $\vF(u) \equiv 0 \mod r$. Up to multiplying $u$ by an $r$-th power of $\pi_F$, we may assume that $\vF(u) = 0$. The residue field $\OF / (\pi_F)$ has size $r^f -1$, so Fermat's little theorem yields $u^{r^{f}-1} \equiv 1 \mod \pi_F$. Nevertheless, $u^{r^{f}-1} \not \equiv 1 \mod \pi_F^2$, as, if that was not the case, the strong version of Hensel's lemma would imply that $u$ is an $r$-th power in $\OF$. It follows that $u^{r^{f-1}} - u^{1/r}$ is a uniformizer of $\O_M$, so the ring of integers of $M$ is described by $\O_M = \OF \left[u^{r^{f-1}} - u^{1/r} \right] = \OF \left[u^{1/r} \right]$ (see \cite[Chapter III \S 6]{Serre79}). Taking valuations, we obtain the equalities
		\begin{equation*}
			\vF(\Delta(M / F)) = \vF(\disc(x^r - u)) = \vF(\pm r^r u^{r-1}) = r \, e_{F / \Qr}.
		\end{equation*} 
	
		\item Assume now that $\vF(u) \not \equiv 0 \mod r$, and write $\eta \coloneqq \vF(u)$. Since $M/F$ is totally ramified, $u^{1/r} / \pi_M^{\eta - 1}$ is also a uniformizer of $O_M$, and so $\O_M = \O_F \left[u^{1/r} / \pi_M^{\eta - 1} \right]$. Just as above, we obtain
		\begin{equation*}
			\vF(\Delta(M/ F)) = \vF \left(\pm r^r \left( u  / \pi_F^{\eta-1} \right)^{r-1} \right) = r \, \vF(r) + (r-1) \, \vF \left(u / \pi_F^{\eta-1} \right) = r \, e_{F / \Qr} + r-1.
		\end{equation*}
	\end{enumerate}\vspace{-1.5em}
\end{proof}

Recall that we assume that $q = r$, and so we have $\nu_r = v_r(s_0(s_0-1))$. When the polynomial $\gminus$ is reducible over $\Qr$, the extension $\Qr(\gamma_0) / \Qr$ is trivial (\cf \ Proposition~\ref{prop: gminus reducible}). Therefore we focus on the case where $\gminus$ is irreducible over $\Qr$.

\begin{theorem}\label{thm: Discr Kr gamma0}
	The $r$-adic valuation of the discriminant of $\Qr(\gamma_0) / \Qr$ is given by
	\begin{equation*}
		v_r(\Delta(\Qr(\gamma_0) / \Qr)) = \begin{cases}
			r & \text{ if } \nu_r = 2, \\
			\frac{3r-1}{2} & \text{ if } \nu_r = 1, \\
			r & \text{ if } \nu_r \leq 0, \nu_r \equiv 0 \mod r \text{ and } \gminus \text{ is irreducible over } \Qr, \\
			2r-1 & \text{ if } \nu_r < 0 \text{ and } \nu_r \not \equiv 0 \mod r.
		\end{cases}
	\end{equation*}
\end{theorem}

\begin{proof}
	In all the listed cases $\gminus$ is irreducible over $\Qr$, so $\Qr(\gamma_0) / \Qr$ and $\Qr(\deltaQ \alpha_0) / \Qone$ are totally ramified of degree $r$ (see Remark~\ref{rmk: gminus irred Qq Kq} and Theorem~\ref{thm: ramif index}). In order to compute $v_r(\Delta(\Qr(\gamma_0) / \Qr))$, we make use of the diagram introduced in Figure~\ref{fig: Diamond diagram}. With the notation of Remark~\ref{rmk: Ram ind quadratics}, we have $d_0 = d_0' \in \lbrace 1, 2 \rbrace$, and $n_0 = n_0' = r$. Corollary 2.10 in \cite[Chapter III \S 2]{Neukirch99} implies that
	\reqnomode
	\begin{equation}\label{eq: Neukirch discr}
		\Delta(\Qr(\gamma_0) / \Qr)^{d_0}  = \Delta(\Qone / \Qr)^r \, \frac{\Nrm{\Qone}{\Qr} \left( \Delta(\Qr(\deltaQ \alpha_0) / \Qone) \right)}{\Nrm{\Qr(\gamma_0)}{\Qr} \left( \Delta(\Qr(\deltaQ \alpha_0) / \Qr(\gamma_0)) \right)}.
	\end{equation}
	Recall that we denote by $\vQ$ a valuation on $\Qrbar$ normalised with respect to $\Qone$. To simplify the notation, we write $\vgam$ to denote a valuation on $\Qrbar$ normalised with respect to $\Qr(\gamma_0)$. We claim that the $r$-adic valuation of the discriminant of $\Qr(\gamma_0) / \Qr$ is given by	
	\begin{equation}\label{eq: Discr Krgam from almost Kummer}
		v_r(\Delta(\Qr(\gamma_0) / \Qr)) = \begin{cases}
			\vQ(\Delta(\Qr(\deltaQ \alpha_0) / \Qone)) & \text{ if } \Qone / \Qr \text{ is unramified,} \\
			\frac{r-1}{2} + \frac{1}{2} \, \vQ(\Delta(\Qr(\deltaQ \alpha_0) / \Qone)) & \text{ if } \Qone / \Qr \text{ is ramified.}
		\end{cases}
	\end{equation}
	Indeed, if $\Qone / \Qr$ is unramified, then $v_r(\Delta(\Qone / \Qr)) = \vgam(\Delta(\Qr(\deltaQ \alpha_0) / \Qr(\gamma_0))) = 0$, and the claim follows from \eqref{eq: Neukirch discr}. On the other hand, if $\Qone / \Qr$ is ramified, then $v_r(\Delta(\Qone / \Qr)) = 1$, and $\vgam(\Delta(\Qr(\deltaQ \alpha_0) / \Qr(\gamma_0))) = 1$. Since $\Qr(\gamma_0) / \Qr$ is totally ramified, and $v_r$ is normalised with respect to $\Qr$, we deduce that $v_r(\Nrm{\Qr(\gamma_0)}{\Qr} \left( \Delta(\Qr(\deltaQ \alpha_0) / \Qr(\gamma_0)) \right)) = 1$. Moreover, $\Qone / \Qr$ being ramified imposes $d_0 = 2$, so taking $r$-adic valuations in \eqref{eq: Neukirch discr} and simplifying terms yields \eqref{eq: Discr Krgam from almost Kummer}.
	
	Having established the claim, we conclude using Proposition~\ref{prop: Valuation discr} applied to $F = \Qone$ and the element $u \coloneqq (\deltaQ \alpha_0)^r = \deltaQ^r(\sqts + \sqtsm)^2$. We treat different cases separately.
	\begin{enumerate}[leftmargin=*]
		\item If $\nu_r = 2$, then $v_r(\deltaQ) = 0$ by Hypothesis~\ref{hyp: v(s0(s0-1)) = 0}, so one among $v_r(s_0), v_r(s_0-1)$ vanishes, and $v_r(u) = 0$. Remark~\ref{rmk: ramif Qone Qr} shows that $\Qone / \Qr$ is unramified, so $v_r(\Delta(\Qr(\gamma_0) / \Qr)) = \vQ(\Delta(\Qr(\deltaQ \alpha_0) / \Qone)) = r$. 
		
		\item If $\nu_r = 1$, we have again $v_r(u) = 0$, but this time $\Qone / \Qr$ is ramified as $v_r(\deltaQ^{2r} s_0(s_0-1)) = 1$ (see Remark~\ref{rmk: ramif Qone Qr}). It follows that $v_r(\Delta(\Qr(\gamma_0) / \Qr)) = \frac{r-1}{2} + \frac{r}{2} = \frac{3r - 1}{2}$.
		
		\item If $\nu_r \leq 0$, $\nu_r \equiv 0 \mod r$ and $\gminus$ is irreducible over $\Qr$, then Remark~\ref{rmk: nur leq 0} states that $\nu_r \leq 0$. But Hypothesis~\ref{hyp: val equiv 0 mod r} yields $v_r(\deltaQ^{2r} s_0 (s_0-1)) = 0$, so $\Qone / \Qr$ is unramified. Thus we have $v_r(u) \equiv 0 \mod r$, and we conclude that $v_r(\Delta(\Qr(\gamma_0) / \Qr)) = r$.
		
		\item Finally, assume that $\nu_r \leq 0$ and $\nu_r \not \equiv 0 \mod r$. If $\Qone / \Qr$ is unramified, then Proposition~\ref{prop: Valuation discr} gives $v_r(\Delta(\Qr(\gamma_0) / \Qr)) = \vQ(\Delta(\Qr(\deltaQ \alpha_0) / \Qone)) = r + r-1$. On the other hand, if $\Qone / \Qr$ is ramified, then $v_r(\Delta(\Qr(\gamma_0) / \Qr)) = \frac{r-1}{2} + \frac{1}{2} ( 2r + r-1) = 2r-1$. In both cases we get the desired result.
	\end{enumerate}\vspace{-1em}
\end{proof}

\subsection{Conclusion of the proof of Theorems~\ref{thm: reduction types Jms} and \ref{thm: reduction types Jpl}}\label{sect: Conclusion red types}

To conclude this section, we describe the reduction type of the Néron models $\Jmodpmtw$ at $\frakq$, finally proving Theorems~\ref{thm: reduction types Jms} and \ref{thm: reduction types Jpl}. Recall that we denote by $q$ the rational prime lying below $\frakq$, and we let $\nu_q \coloneqq v_q(s_0(s_0-1))$. Recall from Corollary~\ref{cor: Discriminants Wpm} the description of the discriminants $\Delta(\Wpmtw)$. The goal of the twisting parameter $\deltaK$ is to minimise the conductor exponent of $\Jpmtw$, so we take $\deltaK$ to be supported at primes dividing $2, r$ and $\deltaQ^{2r} s_0 (s_0 -1)$.

\begin{proof}[Proof of Theorem~\ref{thm: reduction types Jms}]
	We know from Theorem~\ref{thm: Jpm have RM} that $\Jmstw$ has \RM \ by $\Kgl$, so when $\Jmodmstw$ has bad reduction at $\frakq$ and is not semistable, it has automatically unipotent reduction, according to Proposition~\ref{prop: reduction RM}. If $\frakq$ is even, we use the results from \S \ref{sect: Red even places}. If $\frakq$ is odd, we use the criteria on cluster pictures given in Theorems~\ref{thm: Crit clusters good} and \ref{thm: Crit clusters semist}. By Definition~\ref{def: Twisted Cpm}, $\Cmstw$ is the quadratic twist by $\deltaK$ of $\Cms$, and we describe it by the hyperelliptic equation $(\Hhyp_r^{-})^{(\deltaK)} : y^2 = \deltaK \,  \gminus(x)$. Note that its cluster picture at $\frakq$ is the same as the one of $\Cms$, but the valuation of the leading coefficient for $\Cmstw$ is $\vq(\deltaK)$.	
	\begin{enumerate}[leftmargin=*, itemsep=0pt]
		\item Suppose that $\frakq \mid 2$. If $\nu_2 > 0$, then Proposition~\ref{prop: Toric red Jmodmstw} states that $\Jmodmstw$ has toric reduction at $\frakq$ if and only if $\vq(\deltaK) = 0$. On the other hand, assume that $\nu_2 \leq -8$ and that $\SQ{\deltaK \deltaQ^r (2-4s_0)}$ holds. If $\nu_2 \equiv -8 \mod r$, then $\Jmodmstw$ has good reduction over $\Kq$, and otherwise it attains good reduction over any finite extension of $\Kq$ with ramification index $r$, so it has unipotent reduction at $\frakq$ (see Proposition~\ref{prop: Pot good Cmstw} and Corollary~\ref{cor: Good red Jmstw})
		
		\item Suppose that $\frakq = \frakr$. We have the equality $\vr(s_0(s_0-1)) = \frac{r-1}{2} \, \nu_r$.
		\begin{enumerate}[leftmargin=*, itemsep=0pt]
			\item Assume that $\nu_r > 2$, so that $\vr(s_0(s_0-1)) > r$. Every cluster has at most two odd children, and Propositions~\ref{prop: gamma0 rational}, \ref{prop: gminus reducible} imply that $\Krroots = \Qtau$, so $e_{\Krroots / \Kr} \leq 2$. The only proper cluster strictly containing any twin is $\Rroots$, and $d_{\Rroots} = 1$. Therefore, \cite[Theorem $1.3 \ (iv)$]{Bisatt} implies that every twin is invariant under the action of $I_{\Kr}$, and clearly, $\Rroots$ is $I_{\Kr}$-invariant too. Moreover, $\Rroots$ is the only principal cluster: we have $d_{\Rroots} = 1$, and $\nu_{\Rroots} = \vr(\deltaK) + r d_{\Rroots} = \vr(\deltaK) + r$. Thus, $\Cmstw$ satisfies the semistability criterion (so $\Jmodmstw$ has toric reduction at $\frakr$) if and only if $\vr(\deltaK) = 1$.
			
			\item Assume that $\nu_r \leq 2$, so that $\vr(s_0(s_0-1)) \leq r$. The only proper cluster is $\Rroots$, which has size $r = 2g+1$ (where $g = \frac{r-1}{2}$ is the genus of $\Cms$), so $\Cmstw$ has potential good reduction at $\frakr$. Proposition~\ref{prop: Krroots ramified nur leq 2} states that $\Krroots / \Kr$ is ramified, and then Theorem~\ref{thm: Crit clusters good} implies that $\Jmodmstw$ has unipotent reduction at $\frakr$.
		\end{enumerate}		
		
		\item Suppose now that $\frakq \nmid 2r$. Since $\Kq / \Qq$ is unramified, we have $\vq(s_0(s_0-1)) = \nu_q$. 
		\begin{enumerate}[leftmargin=*, itemsep=0pt]
			\item Assume that $\nu_q > 0$. Every cluster has at most two odd children, and Propositions~\ref{prop: gamma0 rational}, \ref{prop: gminus reducible} imply that $e_{\Kqroots / \Kq}^{ } \leq 2$. Again, the only proper cluster strictly containing any twin is $\Rroots$ (which is invariant under $I_{\Kq}$), and $d_{\Rroots} \in \Z$. Therefore, Theorem $1.3 \ (iv)$ from \cite{Bisatt} implies that every twin is $I_{\Kq}$-invariant. The only principal cluster is $\Rroots$: we have $d_{\Rroots} = 0$, and $\nu_{\Rroots} = \vq(\deltaK) + r d_{\Rroots} = \vq(\deltaK)$. If $\vq(\deltaK) = 0$ then $\Jmodmstw$ has toric reduction at $\frakq$, and otherwise it has unipotent reduction.
			
			\item Assume that $\nu_q \leq 0$ and $\nu_q \equiv 0 \mod r$. The only proper cluster is $\Rroots$, which has size $2g+1$. Hypothesis~\ref{hyp: val equiv 0 mod r} gives $\vq(\deltaQ^{2r} s_0 (s_0-1)) = 0$, and Lemma~\ref{lem: Ramification quadratic} combined with Theorem~\ref{thm: ramif index} imply that $\Kqroots / \Kq$ is unramified. We finally compute \vspace{-0.25em}
			\begin{equation*}
				\nu_{\Rroots} = \vq(\deltaK) + r d_{\Rroots} = \vq(\deltaK) + r \vq(\deltaQ) + \frac{1}{2} \vq(s_0(s_0-1)) = \vq(\deltaK \deltaQ^r (2 - 4s_0)). \vspace{-0.25em}
			\end{equation*}
			The last equality follows from $(2-4s_0)^2 = 4 (4s_0 (s_0-1) + 1)$. Thanks to Theorem~\ref{thm: Crit clusters good}, we conclude that if $\vq(\deltaK \deltaQ^r (2 - 4s_0)) \in 2 \Z$, then $\Jmodmstw$ has good reduction at $\frakq$, and otherwise has unipotent reduction.
			
			\item Assume that $\nu_q \leq 0$ and $\nu_q \not \equiv 0 \mod r$. The only proper cluster is $\Rroots$, which has size $2g+1$ ($g$ being the genus of $\Cms$), so $\Cmstw$ has potential good reduction at $\frakq$. However, Proposition~\ref{prop: gminus irreducible} states that $\gminus$ is irreducible, so Theorem~\ref{thm: ramif index} implies that $\Kqroots / \Kq$ is ramified. Again, we deduce from Theorem~\ref{thm: Crit clusters good} that $\Jmodmstw$ has unipotent reduction at $\frakq$.
		\end{enumerate}
	\end{enumerate} \vspace{-1em}
\end{proof}

We conclude the section by describing the reduction types of the Néron model $\Jmodpltw$.

\begin{proof}[Proof of Theorem~\ref{thm: reduction types Jpl}]
	Recall that, as we manipulate $\Cpl$, we assume that $\deltaQ \in \Z$, so $\vq(\deltaQ) \in \Z$, and $\Kq(\Rroots) = \Kq(\Rroots^+)$. Therefore the description of the splitting field is the same one as for $\Cms$. We follow the same strategy as in the proof of Theorem~\ref{thm: reduction types Jms}. When $\frakq$ is odd, note that the cluster picture of $\Cpltw$ equals that of $\Cpl$, except that the leading coefficient has valuation $\vq(\deltaK)$.
	\begin{enumerate}[leftmargin=*, itemsep=0pt]
		
		\item Suppose that $\frakq \mid 2$. If $\nu_2 > 0$, then Proposition~\ref{prop: Toric red Jmodpltw} states that $\Jmodpltw$ has toric reduction at $\frakq$ if and only if $\vq(\deltaK) = 0$. On the other hand, assume that $\nu_2 \leq -4$ and that $\SQ{\deltaK}$ holds. If $\nu_2 \equiv 0 \mod r$, then $\Jmodpltw$ has good reduction over $\Kq$, and otherwise it attains good reduction over any finite extension of $\Kq$ with ramification index $r$ (see Proposition~\ref{prop: Pot good Cmstw} and Corollary~\ref{cor: Good red Jmstw}). In particular, in the latter case $\Jmodpltw$ has unipotent reduction at $\frakq$.

		\item Suppose that $\frakq = \frakr$. We have the equality $\nu_r = \frac{r-1}{2} \vr(s_0(s_0-1))$. 
		\begin{enumerate}[leftmargin=*, itemsep=0pt]
			\item Assume that $v_r(s_0) > 2$, so that $\vr(s_0) > r$. 
			Every cluster has at most two odd children, and Propositions~\ref{prop: gamma0 rational} and \ref{prop: gminus reducible} give $e_{\Krroots / \Kr}^{ } \leq 2$. The proper clusters are the twins, and the only cluster strictly containing any twin is $\Rroots^+$, whose depth is an integer. Again, \cite[Theorem $1.3 \ (iv)$]{Bisatt} implies that every twin is invariant under the action of $I_{\Kr}$, and $\Rroots^+$ is $I_{\Kr}$-invariant too. Moreover, $\Rroots^+$ is the only principal cluster, and we have $\nu_{\Rroots^+} = \vr(\deltaK) + (r+1) d_{\Rroots^+} = \vr(\deltaK) + (r+1)$. We conclude that $\Cpltw$ satisfies the semistability criterion if and only if $\vr(\deltaK) = 0$. 
			
			\item Assume that $v_r(s_0-1) > 2$, so that $\vr(s_0-1) > r$. Just as in the previous case, we have $e_{\Krroots / \Kr} \leq 2$. Since $\gminus \in \OK[x]$, $\Rroots$ is $I_{\Kr}$-invariant. The only cluster strictly containing any twin is $\Rroots$, whose depth is an integer, so every twin is also $I_{\Kr}$-invariant. Moreover, $\Rroots$ is the only principal cluster, and we have $\nu_{\Rroots} = \vr(\deltaK) + r d_{\Rroots} + d_{\Rroots^+} = \vr(\deltaK) + r$. We conclude that $\Cpltw$ satisfies the semistability criterion (Theorem~\ref{thm: Crit clusters semist}) if and only if $\vr(\deltaK) = 1$.
			
			\item Assume that $\nu_r \leq 2$ so that $\vr(s_0(s_0-1)) \leq r$. The only proper clusters are $\Rroots$ and $\Rroots^+$, whose size is $ \geq 2g+1$, so $\Cpltw$ has potential good reduction at $\frakr$. Again, Lemma~\ref{lem: Ramification quadratic} and Theorem~\ref{thm: ramif index} show that $\Krroots / \Kr$ is ramified. We deduce from Theorem~\ref{thm: Crit clusters good} that $\Jmodpltw$ has unipotent reduction at $\frakr$.
		\end{enumerate}		
		
		\item Suppose now that $\frakq \nmid 2r$. Since $\Kq / \Qq$ is unramified, we have $\nu_q = \vq(s_0(s_0-1))$. 
		\begin{enumerate}[leftmargin=*, itemsep=0pt]
			\item Assume that $\nu_q > 0$. Whether $\vq(s_0) > 0$ or $\vq(s_0-1) > 0$, every cluster has at most two odd children. Again, $\gminus$ is reducible by Proposition~\ref{prop: gminus reducible}, so $\Kqroots = \Qtau$ and $e_{\Kqroots / \Kq} \leq 2$. We deduce from \cite{Bisatt} that every proper cluster is $I_{\Kq}$-invariant. Finally, the only principal cluster is $\Rroots^+$, and $\nu_{\Rroots^+} = \vq(\deltaK)$. We conclude that $\Cpltw$ satisfies the semistability criterion if and only if $\vq(\deltaK) = 0$. 
			
			\item Assume that $\nu_q \leq 0$ and $\nu_q \equiv 0 \mod r$. The only proper cluster is $\Rroots^+$, which has size $2g+2$, so $\Cpltw$ has potential good reduction. Hypothesis~\ref{hyp: val equiv 0 mod r} yields $\vq(\deltaQ^{2r} s_0 (s_0-1)) = 0$, so $\Qtau / \Kq$ is unramified, and $\Kqroots / \Qtau$ too because $\vq(s_0(s_0-1)) \equiv 0 \mod r$ (see Theorem~\ref{thm: ramif index}). This time $\nu_{\Rroots^+} = \vq(\deltaK) + (r+1) (\vq(\deltaQ) + \vq(s_0)/r)$: if $\vq(\deltaK) = 0$, then $\Jmodpltw$ has good reduction at $\frakq$, and otherwise it has unipotent reduction.
			
			\item Assume that $\nu_q \leq 0$ and $\nu_q \not \equiv 0 \mod r$. Just as above, $\Cpltw$ has potential good reduction, but now $r$ divides $e_{\Kqroots / \Kq}$ by Theorem~\ref{thm: ramif index}. We conclude that $\Jmodpltw$ does not have good reduction at $\frakq$, and thus has unipotent reduction.
		\end{enumerate}
	\end{enumerate} \vspace{-1em}
\end{proof}

\begin{corollary}\label{cor: Semistability defect}
	Assume that $\frakq$ is odd, that $\Jpmtw$ has unipotent and potential good reduction at $\frakq$. If $e_{\Kqroots / \Kq} \vq(\deltaK (1 - \zr)^r \deltaQ^r \sqrt{s_0(s_0-1)}) \in 2 \Z$, then $\Jpmtw$ attains good reduction over $\Kqroots$. 
\end{corollary}
 
\begin{proof}
	We use the criterion on cluster pictures for good reduction (Theorem~\ref{thm: Crit clusters good}) to prove that $\Jpmtw$ attains good reduction over $\Kqroots$. The first two items from the mentioned theorem are clearly satisfied, and if $\deltaK$ is chosen as in the statement, the third one holds too.
\end{proof}

% -- % -- % -- % -- % -- % -- % -- % -- % -- % -- % -- % -- % -- % -- % -- % -- % -- % -- % -- % -- % -- % -- % -- % -- % -- % -- % -- % -- % -- % -- % -- % -- % -- % -- % -- % -- % -- % -- % -- 

\section{Properties of the $2$-dimensional representations $\rhojtw$}\label{sect: Main props reps}

In this section, we study the $2$-dimensional representations $\rhojtw$ that arise from the fact that $\Jpmtw$ have real multiplication by $\Kgl$ (see Theorem~\ref{thm: Jpm have RM}). We begin by proving that the compatible system $(\rhojtw)_{\lambda}^{ }$ is modular. After this, we will study its conductor and the inertial local types of the attached \WD-representations at places of bad reduction. We will then establish absolute irreducibility of the residue representations. Finally, we will use level lowering results to obtain a newform giving rise to the representation $\rhojtw$ whose level is supported at primes that we control. 
\medskip 

Throughout the section, we keep using the notation introduced in \S \ref{sect: Common framework}. From now on, we also assume that the following hypothesis holds.

\begin{hypothesis}\label{hyp: place toric red}
	There is at least one prime number $q \neq r$ such that $v_{q}(s_0(s_0-1))  > 0$. 
\end{hypothesis}

\begin{remark}
	Hypotheses~\ref{hyp: place toric red} is restrictive, and implies a loss of generality. However, assuming that it is not satisfied imposes huge restrictions on the coefficients $A, B, C$ and any primitive non-trivial solution $(a, b, c)$ to \eqref{eq: GFE ppr} or \eqref{eq: GFE rrp}. Theorem~\ref{thm: reduction types Jms} and \ref{thm: reduction types Jpl} imply that $\Jpmtw$ have at least one place $\neq \frakr$ of potentially toric reduction.
\end{remark}

\subsection{Modularity}\label{sect: Modularity}

In this subsection, we prove that the compatible system of representations $(\rhopmK)_{\lambda}^{ }$ arises from a Hilbert newform over $\Kgl$ of parallel weight $2$ and trivial character. We refer the reader to \cite{Freitag90} for an introduction on Hilbert modular forms. The discussion below is deeply inspired by Darmon's panorama for propagating modularity among Frey objects of different signatures. Indeed, we are going to use the Frey representation of signature $(r,r,r)$ (used for solving Fermat's last theorem) to establish modularity of $\rhopmK$. The content below is a generalisation of some of the results in \cite[\S 4]{BCDF23} and \cite[\S 2]{Darmon00}. 
\medskip

\begin{remark}\label{rmk: Modularity and twisting}
	If $\Kgl(\sqrt{\deltaK}) = \Kgl$, then $\Jpmtw \simeq \Jpm$ over $\Kgl$, so proving modularity of the latter implies modularity of the former. If $\Kgl(\sqrt{\deltaK}) / \Kgl$ is non-trivial, denote by $\chi_{\deltaK}^{ } : G_{\Kgl} \rightarrow \Qlbar\st$ the associated character of $G_{\Kgl}$. For any finite place $\lambda$ of $\OK$ of good reduction for $\Jpm$, we have an isomorphism \vspace{-0.25em}
	\reqnomode
	\begin{equation}
		\rhojtw \simeq \rho_{\Jpm, \, \lambda} \otimes \chi_{\deltaK}^{ }
	\end{equation}
	In order to prove that $\rhojtw$ is modular, it suffices to show that $\rhojlam$ is modular.
\end{remark}

Recall that $\frakr$ denotes the unique prime ideal of $\OK$ lying above $r$. Since $\Kgl / \Q$ is totally ramified at $r$, the residue field $\OK / \frakr$ is isomorphic to $\Fr$. We are going to treat the particular case $\lambda = \frakr$, and prove that $\rhopmbar : G_{\Kgl} \rightarrow \GL_{2}(\Fr)$ is modular. 

\begin{lemma}\label{lem: reps extend to GQ}
	The representation $\rhomsbar$ extends to an odd representation of $G_{\Q}$.
\end{lemma}

\begin{proof}
	Recall that $\Jpm$ has \RM \ by $\Kgl$, so $\V_r(\Jpm)$ carries a structure of $\Kgl \otimes_{\Q} \Q_r$-module. But $r$ is totally ramified in $\Kgl / \Q$, so $\Kgl \otimes_{\Q} \Q_r \simeq \Kr$. By definition, $\CmssQ$ is defined over $\Q$, so the action of $G_{\Kgl}$ on the $\Kr$-module $\V_r(\Jpm)$ extends to a semilinear $G_{\Q}$-action. Therefore, the action of $G_{\Q}$ on $\T_r(\Jpm) \otimes_{\O_{\frakr}} \Ffr$ is $\Ffr$-linear, and restricts to the action of $G_{\Kgl}$ given by $\rhomsbar$.
\end{proof}

\begin{proposition}\label{prop: rhomsbar irred}
	The restriction of $\rhomsbar$ to $G_{\Q(\zr)}$ is absolutely irreducible.
\end{proposition}

\begin{proof}		
	Specifying $p = r$, Theorem~\ref{thm: Frey rep ppr} implies that $\overline{\rho}_{\Jms(s), \, \frakr} : G_{\Kgl(s)} \rightarrow \GL_{2}(\Ffr)$ is a Frey representation of signature $(r, r, r)$. Consider the Legendre elliptic curve $\Leg(s) / \Kgl(s)$ described by $y^2 = x(x-1)(x - s)$. Darmon proves in \cite[\S 1.3]{Darmon00} that $\overline{\rho}_{\Leg(s), \, r} : G_{\Kgl(s)} \rightarrow \GL_{2}(\Fr)$ is also a Frey representation of signature $(r,r,r)$, and that the latter is unique up to equivalence. Thus, there is a character $\eps : G_{\Kgl} \rightarrow \Frbar\st$ giving the isomorphism
	\reqnomode 
	\begin{equation}\label{eq: isom Frey Legendre}
		\rhomsbar \simeq \overline{\rho}_{\Leg(s_0), \, r} \otimes \eps.
	\end{equation}
	Now $\det \overline{\rho}_{\Leg(s_0), \, r}$ and $\det \rhopmbar$ are both the $\mod r$ cyclotomic character (\cite[Theorem 2.8]{BCDF23}), so $\restr{\eps}{G_{\Q(\zr)}}$ is trivial, and $\eps$ has order at most $2$. To establish the proposition, it suffices to show that the restriction of  $\overline{\rho}_{\Leg(s_0), \, r}$ to $G_{\Q(\zr)}$ is absolutely irreducible.
	\medskip 
	
	We claim that $\Leg(s_0) / \Kgl$ does not have complex multiplication (\CM). We are going to prove this claim by showing that $\Leg(s_0)$ does not have potential good reduction everywhere. Hypothesis~\ref{hyp: place toric red} implies that $\Jms$ has potential toric reduction at least at one place $\frakq \neq \frakr$. Let $L / \Kq$ be a finite extension where $\Jms$ attains toric reduction, and write $I_L$ for its inertia group.	Theorem~\ref{thm: Jpm have RM} states that $\End_{\Kgl(s)}(\Jms(s)) \simeq \OK$, so if we restrict $\rhomsbar$ to $G_L$, then $r$ is \textit{good} in the sense of \cite[II \S 2]{Ribet76}. Lemma $3.5.3$ of \textit{loc. cit.} states that there is an additive character $\psi : G_{L} \rightarrow \Ffr$ such that \vspace{-0.25em}
	\begin{equation}\label{eq: Descript rhomsbar Ribet}
		\restr{\rhomsbar}{G_L} = \begin{pmatrix}
			\overline{\chi_r} & \psi \\
			0 & \mathbbm{1}
		\end{pmatrix}, \vspace{-0.25em}
	\end{equation}
	where $\overline{\chi_r} : G_{L} \rightarrow \Ff_r\st$ is the $\mod r$ cyclotomic character. Since $\Jmodmstw$ has toric reduction at $\frakq$, then $\psi$ has to be ramified. Now $\eps$ is trivial on $G_{\Q(\zr)}$ and $\Q(\zr) / \Kgl$ is unramified outside of $\frakr$, so $\restr{\eps}{I_{L}}$ is trivial, as $\frakq \neq \frakr$. Therefore, the restriction of $\overline{\rho}_{\Leg(s_0), \, r}$ to the inertia group $I_{L}$ is also described by the RHS of \eqref{eq: Descript rhomsbar Ribet}, whose action is unipotent and non-trivial. Grothendieck's inertial criterion implies then that $\Leg(s_0)$ has potential multiplicative reduction at $\frakq$ (see also \cite[Chapter IV \S 10]{Silverman94}).
	\medskip 
	
	To conclude the proof, note the equality $\overline{\rho}_{\Leg(s_0), \, r}(G_{\Q(\zr)}) = \overline{\rho}_{\Leg(s_0), \, r}(G_{\Q}) \cap \SL_2(\Fr)$. Since $\Leg(s_0)$ has no \CM \ and $r \geq 5$, Propositions 3.1 and 4.3 in \cite{Najman24} imply that $\overline{\rho}_{\Leg(s_0), \, r}(G_{\Q(\zr)}) = \SL_2(\Fr)$. Assume by contradiction that the restriction of $\rhomsbar$ to $G_{\Q(\zr)}$ is absolutely reducible. Then its image is isomorphic to a subgroup of upper triangular matrices, or to a non-split Cartan subgroup. In both cases, this image is solvable, contradicting the fact that $\SL_2(\Fr)$ is not solvable for $r \geq 5$.
\end{proof}

\begin{theorem}\label{thm: Modularity Jms}
	The compatible system of representations $(\rho_{\Jms, \, \lambda})_{\lambda}^{ }$ is modular, \ie it arises from a Hilbert newform defined over $\Kgl$.
\end{theorem}

\begin{proof}
	Lemma~\ref{lem: reps extend to GQ} states that $\rhomsbar : G_{\Kgl} \rightarrow \GL_{2}(\Fr)$ extends to an odd representation $\overline{\rho}$ of $G_{\Q}$. By Proposition~\ref{prop: rhomsbar irred}, $\restr{\overline{\rho}}{G_{\Q(\zr)}}$ is absolutely irreducible, so $\overline{\rho}$ is also absolutely irreducible. Serre's conjecture (\cite{ChandrashekharWintenbergerI, ChandrashekharWintenbergerII}) implies that $\overline{\rho}$ is modular, and by cyclic base change, $\rhomsbar$ is modular too (\cite{Langlands80}). Theorem~\ref{thm: reduction types Jms} implies that $\rho_{\Jms, \, \frakr}$ is unramified almost everywhere. Moreover, $\Jms$ is potentially semistable, so $\rho_{\Jms, \, \frakr}$ is de Rham, hence Hodge--Tate, with HT weights $\lbrace 0, 1 \rbrace$. Applying \cite[Theorem 1.1]{ChandrashekharThorne}, we conclude that $\rho_{\Jms, \, \frakr}$ is modular.
\end{proof}

\begin{theorem}\label{thm: Modularity Jpl}
	If $\nu_r > 2$, the compatible system of representations $(\rho_{\Jpltw, \, \lambda})_{\lambda}^{ }$ is modular.
\end{theorem}

\begin{proof}
	Theorem~\ref{thm: Frey rep ppr} states that $\overline{\rho}_{\Jpl(s), \, \frakr}$ is a Frey representation of signature $(r, r, r)$, that is even in the sense of \cite[\S 1.1]{Darmon00}. But the only Frey representation (up to equivalence) of signature $(r,r,r)$ is $\overline{\rho}_{\Leg(s), \, r}$, which is odd in the sense of \cite{Darmon00}. Specialising at $s = s_0$, we deduce that $\rhoplbar$ is reducible. The assumption $\nu_r > 2$ implies that $\Jmodpltw$ has toric reduction at $\frakr$, so $\V_r(\Jpltw)$ is an ordinary representation of $G_{\Kgl}$. The result follows from \cite[\S 4.5, Theorem A]{SkinnerWiles99} (with $k = 2$, and the field denoted by $F(\chi_1 / \chi_2)$ being equal to $\Q(\zr)$).
\end{proof}

\begin{remark}
	The work of Pan \cite{Pan22} seems like a promising way to prove modularity of $\rhoplK$ in broader generality. A result in the style of \cite[Theorem 7.11]{Pan22} would allow to drop $r$-adic assumptions on $s_0(s_0 -1)$ in Theorem~\ref{thm: Modularity Jpl} above. However, in order to apply the mentioned result of Pan, we would need $r$ to be completely split in $\Kgl$, which is not satisfied in our setting. 
\end{remark}

\begin{remark}
	By the work of Carayol \cite{Carayol86}, the modularity of $(\rhopm)_{\lambda}^{ }$ provides another proof for the fact that it is strictly compatible system of Galois representations.
\end{remark}

\begin{example}
	When specialising the values of $s_0, \deltaQ$ as in Proposition~\ref{prop: right choices s0 deltaQ}, one can deduce the modularity of the representations arising from the curve $\Cmsabc$ for the signature $(p,p,r)$, and for $\Crabc$ for the signature $(r,r,p)$. For the signature $(p, p, r)$, the representation arising from $\Cplabc$ is modular as soon as $v_r(Aa^p Bb^p) > 2$. 
\end{example}

\subsection{Conductor and inertial local types}\label{sect: Cond & local types}

We now use the results from \S \ref{sect: WD reps ab var RM} and \S \ref{sect: Reduction types} to compute the Artin conductor of $\rhopmK$ at every finite place of $\Kgl$. Modularity theorems state that the level of the Hilbert newform giving rise to $\rhopmK$ equals the global conductor of such representation. For Diophantine applications, it is therefore crucial to understand this conductor in detail. Along the process, we will describe the local inertial types of the complex Weil--Deligne representations associated to $\rhopmK$. This will be helpful later when discussing irreducibility in \S \ref{sect: Level lowering}.
\medskip

In order to simplify the discussion below, let us introduce some notation.

\begin{definition}
	Fix two finite places $\lambda \neq \frakq$ of $\Kgl$ such that $\frakq \mid 2r \deltaQ^{2r} s_0 (s_0-1)$, and $\lambda \nmid 2r \deltaQ^{2r} s_0 (s_0-1)$. We denote by $\ctamepmK$ and $\cwildpmK$ the tame and wild conductors of $\rhopmK$ restricted to $D_{\frakq} \simeq G_{\Kq}$, and by $\cdpmK \coloneqq \ctamepmK + \cwildpmK$. We define the global conductor of $\rhopmK$ as 
	\begin{equation*}
		\Cond{\rhopmK} \coloneqq \prod_{\frakq} \frakq^{\cdpmK}. 	\vspace{-0.25em}
	\end{equation*}
\end{definition}

Strict compatibility of the system $(\rhopm)_{\lambda}^{ }$ implies that the Artin conductor at $\frakq$ does not depend on $\lambda$, so neither does the global conductor $\Cond{\rhopmK}$. It is convenient to choose $\lambda$ coprime with $\frakq$ and of good reduction for $\Jmodpm$, so we assume that this is indeed the case.
\medskip 

We now describe the Artin conductor $\cdpmK$, depending on $\frakq$ and the other involved parameters. Recall that we let $q$ denote the rational prime lying below $\frakq$, and we write $\nu_q = v_q(s_0 (s_0-1))$.

\begin{theorem}\label{thm: Conductor rhomsK}
	The value of the Artin conductor $\cdmsK$ is described in Table~\ref{table: Cond Jmstw}. 
	\begin{table}[h!]
		\renewcommand{\arraystretch}{1.25}
		\captionsetup{justification=centering,margin=1.5cm}
		\centering
		\begin{tabular}{!{\vrule width 1.25pt}  P{1.25cm}  !{\vrule width 1.25pt}  P{0cm} P{4cm}  !{\vrule width 1.25pt}  P{4.25cm}  !{\vrule width 1.25pt}  P{1.25cm}  !{\vrule width 1.25pt}}
			\noalign{\hrule height 1.25pt}
			Place $\frakq$                     & \multicolumn{2}{P{6.5cm}!{\vrule width 1.25pt}}{Behaviour of $\nu_q$ and $\gminus \in \Qq[x]$}      & Condition on $\deltaK$       & $\cdmsK$                             \\ \noalign{\hrule height 1.25pt}
			\multirow{3}{*}{$\frakq \mid 2$}   & \multicolumn{2}{P{6.5cm}!{\vrule width 1.25pt}}{$\nu_2 > 0$}                    & $\vq(\deltaK) = 0$                     & $1$                                    \\ \cline{2-5} 
			& \multicolumn{2}{P{6.5cm}!{\vrule width 1.25pt}}{$\nu_2 \leq -8$ and $\nu_2 \equiv -8 \mod r$}                                                                                & $\SQ{\deltaK \deltaQ^r (2-4s_0)}$                        & $0$                                     \\ \cline{2-5} 
			& \multicolumn{2}{P{6.5cm}!{\vrule width 1.25pt}}{$\nu_2 \leq -8$ and $\nu_2 \not \equiv -8 \mod r$}                                                                                &   $\SQ{\deltaK \deltaQ^r (2-4s_0)}$     & $2$                                     \\ \noalign{\hrule height 1.25pt}
			\multirow{6}{*}{$\frakq = \frakr$} & \multicolumn{2}{P{6.5cm}!{\vrule width 1.25pt}}{$\nu_r > 2$}                                                                     & $\vr(\deltaK) = 1$       & $1$                                    \\ \cline{2-5} 
			& \multicolumn{2}{P{6.5cm}!{\vrule width 1.25pt}}{$\nu_r = 2$}                                                                     & $\emptyset$       & $3$                                    \\ \cline{2-5} 
			& \multicolumn{2}{P{6.5cm}!{\vrule width 1.25pt}}{$\nu_r = 1$}   & $\emptyset$   & $2 + \frac{r+1}{2}$                                    \\ \cline{2-5} 
			& \multicolumn{1}{P{2.5cm}|}{$\nu_r \leq 0$ and} & $\gminus$ reducible over $\Qr$ & $\emptyset$                                      & $2$                                    \\ \cline{3-5} 
			& \multicolumn{1}{l|}{\ $\nu_r \equiv 0 \mod r$}    & $\gminus$ irreducible over $\Qr$                         & $\emptyset$                                      & $3$ \\ \cline{2-5} 
			& \multicolumn{2}{P{6.5cm}!{\vrule width 1.25pt}}{$\nu_r \leq 0$ and $\nu_r \not \equiv 0 \mod r$}                                     & $\emptyset$                                      & $2+r$                                  \\ \noalign{\hrule height 1.25pt}
			\multirow{3}{*}{$\frakq \nmid 2r$} & \multicolumn{2}{P{6.5cm}!{\vrule width 1.25pt}}{$\nu_q > 0$}                                                                       & $\vq(\deltaK) = 0$                     & $1$                                    \\ \cline{2-5} 
			& \multicolumn{2}{P{6.5cm}!{\vrule width 1.25pt}}{$\nu_q \leq 0$ and $\nu_q \equiv 0 \mod r$}                                          & $\vq(\deltaK) = 0$ & $0$                                    \\ \cline{2-5} 
			& \multicolumn{2}{P{6.5cm}!{\vrule width 1.25pt}}{$\nu_q\leq 0$ and $\nu_q \not \equiv 0 \mod r$}                                     & $\emptyset$                            & $2$                                    \\ \noalign{\hrule height 1.25pt}
		\end{tabular}

		\caption{Table describing the Artin conductor of the $\lambda$-adic representation $\rhomsK$ at $\frakq$ in terms of $\nu_q, \gminus$ and $\deltaK$.}
		\label{table: Cond Jmstw}
	\end{table}
\end{theorem}

\begin{proof}
	We use Theorem~\ref{thm: reduction types Jms} to describe the reduction type of the Néron model $\Jmodmstw$, and Proposition~\ref{prop: Tame cond lambda} to recover from this the tame part of the conductor. As explained in Corollary~\ref{cor: Semistability defect}, the choice of $\deltaK$ done in the Figure~\ref{fig: Reduction Jmodmstw} minimises the semistability defect of $\Jmstw$.
	\begin{enumerate}[leftmargin=*, itemsep=0pt, topsep=5pt]
		\item Assume that $\frakq$ is even.
		\begin{enumerate}[leftmargin=*, itemsep=0pt, topsep=0pt]
			\item If $\nu_2 > 0$ and $\vq(\deltaK) = 0$ , then $\Jmodmstw$ has toric reduction at $\frakq$. Proposition~\ref{prop: Tame cond lambda} gives $\ctamemsK = 1$, and $\cwildmsK = 0$ as $\Jmstw$ is semistable. We deduce that $\cdmsK = 1$.
			
			\item If $\nu_2 \leq -8$, $\nu_2 \equiv -8 \mod r$ and $\SQ{\deltaK \deltaQ^r (2-4s_0)}$, then $\Jmodmstw$ has good reduction at $\frakq$, so $\cdmsK = 0$ by the Néron--Ogg--Shafarevich criterion \cite{SerreTate}.
			
			\item If $\nu_2 \leq -8$, $\nu_2 \not \equiv -8 \mod r$ and $\SQ{\deltaK \deltaQ^r (2-4s_0)}$ holds, then $\Jmstw$ attains good reduction over a finite extension of $\Kq$ with ramification index $r$. Now $r$ being odd, this is a tame extension of $\Kq$, so Lemma~\ref{lem: Wild cond lambda} yields $\cwildmsK = 0$, and therefore $\cdmsK = 2$.
		\end{enumerate}
		
		 \item Assume that $\frakq = \frakr$. 
		 \begin{enumerate}[leftmargin=*, itemsep=0pt, topsep=0pt]
		 	\item If $\nu_r > 2$ and $\vr(\deltaK) = 1$, then $\Jmodmstw$ has toric reduction at $\frakr$, so $\cdmsrK = \ctamemsrK = 1$.
		 	
		 	\item If $\nu_r \leq 2$, then $\Jmodmstw$ has unipotent reduction at $\frakr$, so $\ctamemsrK = 2$. 
		 	
		 	If $\gminus$ is reducible over $\Qr$, then $\Jmodmstw$ attains good reduction over $\Qtau$, which is a quadratic extension of $\Kr$. Now $\Qtau / \Kr$ is tame, so $\cwildmsrK = 0$, and thus $\cdmsrK = 2$.
		 	
		 	If $\gminus$ is irreducible over $\Qr$, then $\Jmodmstw$ attains good reduction over a wild extension of $\Kr$. The wild conductor of $\rho_{\Jmstw, \, \ell}$ at $\frakr$ is described in Theorem~\ref{thm: Wild cond clusters}. In particular, it does not depend on $\deltaK$, but only on the roots of $\gminus$, so we may choose $\deltaK = 1$. Since $\Jms / \Kr$ is the base change of $\Jac(\CmssQ) / \Qr$, and $\Kr / \Qr$ is a tame extension, Lemma~\ref{lem: Behaviour nwild tame ext} relates the wild conductor of the $\ell$-adic representations attached to the two Jacobians. Combining it with \eqref{eq: Conductor l adic lambda adic} gives
		 	\reqnomode
		 	\begin{equation}\label{eq: relating wild conds}
		 		\frac{r-1}{2} \, \condwild{\restr{\rho_{\Jac(\CmssQ), \, \ell }}{G_{\Qr}}} = \condwild{\restr{\rho_{\Jmstw, \, \ell} }{G_{\Kr}}} = \frac{r-1}{2} \, \cwildmsrK
		 	\end{equation}
		 	To describe the left-most term, we use Theorem~\ref{thm: Wild cond clusters}. Since $\gminus$ is irreducible over $\Qr$, there is a single orbit in $\Rroots$ under the action of $G_{\Qr}$, and a representative of this is $\gamma_0$, so $[\Qr(\gamma_0) : \Qr] = r$. Theorem~\ref{thm: ramif index} combined with Remarks~\ref{rmk: Ram ind quadratics} and \ref{rmk: gminus irred Qq Kq} imply that $\Qr(\gamma_0) / \Qr$ is totally ramified, so its residue degree is $1$. Then \eqref{eq: relating wild conds} and Theorem~\ref{thm: Wild cond clusters} give
		 	\begin{equation*}
		 		\cwildmsrK =  v_r \left( \Delta(\Qr(\gamma_0) / \Qr) \right) - r + 1 .
		 	\end{equation*}
	 		From the description of $v_r(\Delta(\Qr(\gamma_0) / \Qr))$ given in Theorem~\ref{thm: Discr Kr gamma0}, we obtain 
	 		\begin{equation*}
	 			\cwildmsrK = \begin{cases}
	 				r - r + 1 = 1 & \text{ if } \nu_r = 2, \\
	 				\frac{3r-1}{2} - r + 1 = \frac{r+1}{2} & \text{ if } \nu_r = 1, \\
	 				r -r + 1 = 1 & \text{ if } \nu_r \leq 0, \nu_r \equiv 0 \mod r \text{ and } \gminus \text{ is irreducible over } \Qr, \\
	 				2r-1 - r +1 = r & \text{ if } \nu_r < 0 \text{ and } \nu_r \not \equiv 0 \mod r.
	 			\end{cases}
	 		\end{equation*}
		 \end{enumerate}
	 
	 	\item Assume that $\frakq \nmid 2r$. 
	 	\begin{enumerate}[leftmargin=*, itemsep=0pt, topsep=0pt]
	 		\item If $\nu_q > 0$ and $\vq(\deltaK) = 0$, then $\Jmodmstw$ has toric reduction at $\frakq$, so $\cdmsK = \ctamemsK = 1$.
	 		
	 		\item If $\nu_q \leq 0$, $\nu_q \equiv 0 \mod r$ and $\vq(\deltaK \deltaQ^r (2-4s_0)) \in 2 \Z$, then $\Jmodmstw$ has good reduction at $\frakq$, so the Artin conductor is trivial $\cdmsK = 0$.
	 		
	 		\item If $\nu_q \leq 0$, and $\nu_q \not \equiv 0 \mod r$ then $\Jmodmstw$ has unipotent reduction at $\frakq$, so $\ctamemsK = 2$. Moreover, $\Jmstw$ attains good reduction over $\Kqroots$ or a quadratic extension of it. Both of them are tame, as $\frakq \nmid 2r$, so Lemma~\ref{lem: Wild cond lambda} yields $\cwildmsK = 0$, hence $\cdmsK = 2$. 		
	 	\end{enumerate}
	\end{enumerate}\vspace{-1.5em}
\end{proof}

\begin{theorem}\label{thm: Conductor rhoplK}
	The value of the Artin conductor $\cdplK$ is described in Table~\ref{table: Cond Jpltw}. 
	\begin{table}[h!]
		\renewcommand{\arraystretch}{1.25}
		\captionsetup{justification=centering,margin=1.5cm}
		\centering
		\begin{tabular}{!{\vrule width 1.25pt}  P{1.25cm}  !{\vrule width 1.25pt}  P{0cm} P{4cm}  !{\vrule width 1.25pt}  P{4.25cm}  !{\vrule width 1.25pt}  P{1.25cm}  !{\vrule width 1.25pt}}
			\noalign{\hrule height 1.25pt}
			Place $\frakq$                     & \multicolumn{2}{P{6.5cm}!{\vrule width 1.25pt}}{Behaviour of $\nu_q$ and $\gminus \in \Qq[x]$}      & Condition on $\deltaK$       & $\cdplK$                             \\ \noalign{\hrule height 1.25pt}
			\multirow{3}{*}{$\frakq \mid 2$}   & \multicolumn{2}{P{6.5cm}!{\vrule width 1.25pt}}{$\nu_2 > 0$}                                                                                & $\vq(\deltaK) = 0$                        & $1$                                     \\ \cline{2-5} 
			& \multicolumn{2}{P{6.5cm}!{\vrule width 1.25pt}}{$\nu_2 \leq -4$ and $\nu_2 \equiv 0 \mod r$}                                                                                & $\SQ{\deltaK}$                        & $0$                                     \\ \cline{2-5} 
			& \multicolumn{2}{P{6.5cm}!{\vrule width 1.25pt}}{$\nu_2 \leq -4$ and $\nu_2 \not \equiv 0 \mod r$}                                                                                &   $\SQ{\deltaK}$     & $2$                                     \\ \noalign{\hrule height 1.25pt}
			\multirow{7}{*}{$\frakq = \frakr$} & \multicolumn{2}{P{6.5cm}!{\vrule width 1.25pt}}{$\vr(s_0) > 2$}                                                                     & $\vr(\deltaK) = 0$       & $1$                                    \\ \cline{2-5} 
			& \multicolumn{2}{P{6.5cm}!{\vrule width 1.25pt}}{$\vr(s_0-1) > 2$}                                                                     & $\vr(\deltaK) = 1$       & $1$                                    \\ \cline{2-5} 
			& \multicolumn{2}{P{6.5cm}!{\vrule width 1.25pt}}{$\nu_r = 2$}                                                                     & $\emptyset$       & $3$                                    \\ \cline{2-5}
			& \multicolumn{2}{P{6.5cm}!{\vrule width 1.25pt}}{$\nu_r = 1$}    & $\emptyset$    & $2 + \frac{r+1}{2}$                                    \\ \cline{2-5} 
			& \multicolumn{1}{P{2.5cm}|}{$\nu_r \leq 0$ and} & $\gminus$ reducible over $\Qr$ & $\emptyset$                                      & $2$                                    \\ \cline{3-5} 
			& \multicolumn{1}{l|}{\ $\nu_r \equiv 0 \mod r$}    & $\gminus$ irreducible over $\Qr$                         & $\emptyset$                                      & $3$ \\ \cline{2-5}  
			& \multicolumn{2}{P{6.5cm}!{\vrule width 1.25pt}}{$\nu_r \leq 0$ and $\nu_r \not \equiv 0 \mod r$}                                     & $\emptyset$                                      & $2+r$                                  \\ \noalign{\hrule height 1.25pt}
			\multirow{3}{*}{$\frakq \nmid 2r$} & \multicolumn{2}{P{6.5cm}!{\vrule width 1.25pt}}{$\nu_q > 0$}                                                                       & $\vq(\deltaK) = 0$                     & $1$                                    \\ \cline{2-5} 
			& \multicolumn{2}{P{6.5cm}!{\vrule width 1.25pt}}{$\nu_q \leq 0$ and $\nu_q \equiv 0 \mod r$}                                          & $\vq(\deltaK) = 0$ & $0$                                    \\ \cline{2-5} 
			& \multicolumn{2}{P{6.5cm}!{\vrule width 1.25pt}}{$\nu_q\leq 0$ and $\nu_q \not \equiv 0 \mod r$}                                     & $\emptyset$                            & $2$                                    \\ \noalign{\hrule height 1.25pt}
		\end{tabular}
		\caption{Table describing the Artin conductor of the $\lambda$-adic representation $\rhoplK$ in terms of $\nu_q, \gminus$ and $\deltaK$.}
		\label{table: Cond Jpltw}
	\end{table}
\end{theorem}

\begin{proof}
	Recall that, when dealing with the curve $\Cpltw$, we assume that $\deltaQ \in \Z$. The defining polynomial of $\Cpltw$ is obtained from the one of $\Cmstw$ by adjoining the linear factor $(x + 2 \deltaQ)$. The extra root $\gamma_r = -2 \deltaQ$ is rational, so, whenever $\frakq$ is odd, Theorem~\ref{thm: Wild cond clusters} implies that $\cwildplK = \cwildmsK$ (see \cite[Remark 2.10]{ACIKMM}). Again, we treat different cases separately.
	\begin{enumerate}[leftmargin=*, itemsep=0pt]
		\item Assume that $\frakq$ is even. 
		\begin{enumerate}[leftmargin=*, itemsep=0pt, topsep=0pt]
			\item If $\nu_2 > 0$ and $\vq(\deltaK) = 0$, then $\Jmodpltw$ has toric reduction at $\frakq$. Proposition~\ref{prop: Tame cond lambda} gives $\ctameplK = 1$, and $\cwildplK = 0$ as $\Jpltw$ is semistable. We deduce that $\cdplK = 1$.
			
			\item If $\nu_2 \leq -4$, $\nu_2 \equiv 0 \mod r$ and $\SQ{\deltaK}$ holds, then $\Jmodpltw$ has good reduction at $\frakq$, so $\cdplK = 0$.
			
			\item If $\nu_2 \leq -4$, $\nu_2 \not \equiv 0 \mod r$ and $\SQ{\deltaK}$ holds, then $\Jmodpltw$ has unipotent reduction at $\frakq$, so $\ctameplK = 2$. However, $\Jpltw$ attains good reduction over an extension of $\Kq$ with ramification index $r$, which is tame. Lemma~\ref{lem: Wild cond lambda} yields $\cwildplK = 0$, and thus $\cdplK = 2$.
		\end{enumerate}
		
		\item Assume that $\frakq = \frakr$.
		\begin{enumerate}[leftmargin=*, itemsep=0pt, topsep=0pt]
			\item If $v_r(s_0) > 2$ and $\vr(\deltaK) = 0$, then $\Jmodpltw$ has toric reduction at $\frakr$. Similarly, if $v_r(s_0-1)> 2$ and $\vr(\deltaK) = 1$, then $\Jmodpltw$ has toric reduction at $\frakr$. In both cases we have $\cdplrK = 1$.
			
			\item If $\nu_r \leq 2$, then $\Jmodpltw$ has unipotent reduction at $\frakr$, so $\ctameplrK = 2$. By the discussion at the beginning of the proof, we have $\cwildplrK = \cwildmsrK$, so the result follows from Theorem~\ref{thm: Conductor rhomsK}.
		\end{enumerate}

		\item Assume that $\frakq \nmid 2r$.
		\begin{enumerate}[leftmargin=*, itemsep=0pt, topsep=0pt]
			\item If $\nu_q > 0$ and $\vq(\deltaK) = 0$, then $\Jmodpltw$ has toric reduction at $\frakq$, so $\cdplK = \ctameplK = 1$.
			
			\item If $\nu_q \leq 0$, $\nu_q \equiv 0 \mod r$ and $\vq(\deltaK) = 0$, then $\Jmodpltw$ has good reduction at $\frakq$, so $\cdplK = 0$.
			
			\item If $\nu_q \leq 0$, $\nu_q \not \equiv 0 \mod r$, and $\vq(\deltaK) = 0$, then $\Jmodpltw$ has unipotent reduction at $\frakq$, so $\ctameplK = 2$. As explained above, we have $\cwildplK = \cwildmsK = 0$, so $\cdplK = 2$.
		\end{enumerate}\vspace{-1em}
	\end{enumerate}
\end{proof}

\begin{remark}
	The description of the conductor exponents in Tables~\ref{table: Cond Jmstw} and \ref{table: Cond Jpltw} generalises the results of \cite[Theorems 3.9, 4.24 and 4.27]{ACIKMM}. In \textit{loc. cit.} are computed the Artin conductors at odd places  of the $\ell$-adic representations attached $\Jac(\Crabc) / \Q$ and $\Jac(\Cpmabc) / \Q$, just as those of their base change to $\Kgl$ (this corresponds to choosing trivial coefficients $A = B = C = 1$). In \cite{ACIKMM}, the tame conductor is computed exploiting the combinatorial data of the cluster picture, using \cite[Theorem 12.3]{Hyperuser}. Following the same line of arguments as in \cite{ACIKMM}, the authors of \cite{CazorlaVillagra} compute the Artin conductors for a parametric family of hyperelliptic curves, whose values can be directly deduced from Theorem~\ref{thm: Conductor rhomsK} above.

	The approach that we follow here is different, as we use cluster pictures to describe the reduction types of $\Jpmtw$. We deduce the tame conductors thanks to Proposition~\ref{prop: Tame cond lambda}, without having to compute the sets $U$ and $V$ from \cite[Theorem 12.3]{Hyperuser}. The author believes that the strategy in here is more synthetic and beneficial for our Diophantine purposes, as knowing the reduction types of $\Jmodpmtw$ has more applications than just computing Artin conductors (see Proposition~\ref{prop: rhomsbar irred} or Theorem~\ref{thm: Finiteness}).
\end{remark}

Knowing Theorems~\ref{thm: Conductor rhomsK} and \ref{thm: Conductor rhoplK}, we can describe the global conductor of the strictly compatible systems of representations $(\rhopmK)_{\lambda}$. In order to simplify the notation, we introduce:

\begin{definition}\label{def: n2 ntor nunip}
	Denote by $\ntwopm$ the Artin conductor $\cdpmK$ at any even place of $\Kgl$. Let $\ntor$ the square-free product of the primes $\frakq \nmid 2r$ such that $\nu_q > 0$. Let $\nunip$ the square-free product of the primes $\frakq \nmid 2r$ such that $\nu_q < 0$ and $\nu_q \not \equiv 0 \mod r$.
\end{definition}

\begin{corollary}\label{cor: Global conductor}
	The global conductor of the strictly compatible system $(\rhopmK)_{\lambda}$ is given by
	\begin{equation*}
		\Cond{\rhopmK} = 2^{\ntwopm} \, \frakr^{\cdpmrK} \, \ntor^{ } \, \nunip^2.
	\end{equation*}
\end{corollary}

\begin{example}\label{examp: Conductor Fermat}
	In the style of Example~\ref{examp: Reduction types Fermat}, we describe the Artin conductor $\cdpmK$ for the specific choice of $s_0, \deltaQ$ done in Proposition~\ref{prop: right choices s0 deltaQ}, when $\frakq \nmid 2r$. Again, one could deduce in the same way $\cdpmK$ when $\frakq \mid 2$ or $\frakq = \frakr$, using Theorems~\ref{thm: Conductor rhomsK} and \ref{thm: Conductor rhoplK}. 
	\begin{enumerate}[leftmargin=0em, itemsep=0pt]
		\item[] \boxppr Assume that $(a, b, c)$ is a primitive non-trivial solution to \eqref{eq: GFE ppr}. If $\deltaK$ is well chosen, the Artin conductor at $\frakq$ of the $2$-dimensional $\lambda$-adic representations attached to $\Jac(\Cpmabc)^{(\deltaK)}$ is \vspace{-0.5em}
		\begin{equation*}
			\mathfrak{n} \left( \restr{\rho_{\Jac(\Cpmabc)^{(\deltaK)}, \, \lambda }}{D_{\frakq}} \right) = \begin{cases}
				1 & \text{ if } q \mid ABab, \\
				0 & \text{ if } q \mid c \text{ and } q \nmid C, \\
				2 & \text{ if } q \mid C .
			\end{cases}
		\end{equation*}

		\item[] \boxrrp Assume that $(a, b, c)$ is a primitive non-trivial solution to \eqref{eq: GFE rrp}. If $\deltaK$ is well chosen, the Artin conductor at $\frakq$ of the $2$-dimensional $\lambda$-adic representations attached to $\Jac(\Crabc)^{(\deltaK)}$ is \vspace{-0.5em}
		\begin{equation*}
			\mathfrak{n} \left( \restr{\rho_{\Jac(\Crabc)^{(\deltaK)}, \, \lambda }}{D_{\frakq}} \right) = \begin{cases}
				1 & \text{ if } q \mid Cc, \\
				0 & \text{ if } q \mid ab \text{ and } q \nmid AB, \\
				2 & \text{ if } q \mid AB.
			\end{cases}
		\end{equation*}
	\end{enumerate}
\end{example}

Now that we have computed the Artin conductors of $\rhopmK$, we describe the local inertial type of the associated complex \WD-representations. In order to simplify the notation, we introduce:

\begin{definition}\label{def: WD reps Jpmtw}
	Let $\lambda$ and $\frakq$ be as above. For any complex embedding $\iota : \Klam \hookrightarrow \C$, we denote by $\rhiotpmqK$ the complex \WD-representation \vspace{-0.25em}
	\begin{equation*}
		\rhiotpmqK \coloneqq \W_{\iota} \left(\restr{\rhopmK}{W_{\frakq}} \right) : W_{\frakq} \longrightarrow \GL_2(\C).
	\end{equation*}
\end{definition}

For any $x \in \Z$ coprime to $r$, we denote by $\ord\st(x \mod r)$ the order of $q$ in $(\Z / r \Z)\st$.

\begin{proposition}\label{prop: Local types Jpm}
	Let $\deltaK$ be as described in Theorems~\ref{thm: Conductor rhomsK}, \ref{thm: Conductor rhoplK}. If $\Jpmtw$ has potential good reduction at $\frakq$, let $\deltaK$ be as in Corollaries~\ref{cor: Semist defect q even}, \ref{cor: Semistability defect}.	The inertial local type of $\rhiotpmqK$ is given as follows.
	\begin{enumerate}[leftmargin=*, itemsep=0pt, topsep=5pt]
		\item Assume that $\frakq \mid 2$. If $\nu_2 > 0$, then $\rhiotpmqK$ is Steinberg. If $\nu_2 \leq -8$ and $\nu_2 \not \equiv -8 \mod r$, then $\rhiotmsqK$ is principal series if $\ord\st(2 \mod r)$ is odd, and supercuspidal otherwise. If $\nu_2 \leq -4$ and $\nu_2 \not \equiv 0 \mod r$, then $\rhiotplqK$ is principal series if $\ord\st(2 \mod r)$ is odd, and supercuspidal otherwise.
		
		\item Assume that $\frakq = \frakr$. If $\nu_r > 2$, then $\rhiotpmqK$ is Steinberg, and if $\nu_r \leq 2$, then $\rhiotpmqK$ is principal series.
		
		\item Assume that $\frakq \nmid 2r$. If $\nu_q > 0$, then $\rhiotpmqK$ is Steinberg. If $\nu_q \leq 0$ and $\nu_q \not \equiv 0 \mod r$, then $\rhiotpmqK$ is principal series if $\ord\st(q \mod r)$ is odd, and supercuspidal otherwise.
	\end{enumerate}
	Moreover, whenever $\rhiotpmqK$ is supercuspidal, it is non-exceptional and arises as the induction of a character on the unramified quadratic extension of $\Kq$.
\end{proposition}

\begin{proof}
	We use Proposition~\ref{prop: local inert types RM}. If $\Jmodmstw$ has potential toric reduction at $\frakq$, then $\rhiotpmqK$ is Steinberg, so this rules out the case $\nu_r > 2$ and $\nu_q > 0$ for $q \neq r$. Assume now that $\Jpmtw$ has potential good reduction at $\frakq$. Then $\rhiotpmqK$ is principal series if and only if the prime to $q$ part of the semistability defect $\sd_{\Jpmtw / \Kq}$ divides $|\Ffq\st|$. One can check (see \cite[\S 1]{DarmonMestre}) that the residue degree $[\Ffq : \Ff_q]$ equals $\ord\st(q^2 \mod r)$, so we obtain $|\Ffq\st| = q^{\ord\st(q^2 \mod r)} -1$. On the other hand, we have the equality 
	\begin{equation*}
		\ord\st(q^2 \mod r) = \frac{\ord\st(q \mod r)}{\gcd(2, \, \ord\st(q \mod r))}.
	\end{equation*}
	By definition, $r$ divides $q^{\ord\st(q \mod r)} -1$, so we deduce that $r$ divides $|\Ffq\st|$ if and only if $\ord\st(q \mod r)$ is odd. Knowing this, we treat different cases separately:
	\begin{enumerate}[leftmargin=*]
		\item Assume that $\frakq \mid 2$. If $\nu_2 \leq -8$ and $\nu_2 \not \equiv -8 \mod r$, then Corollary~\ref{cor: Semist defect q even} states that the semistability defect of $\Jmstw$ equals $r$. If $\nu_2 \leq -4$ and $\nu_2 \not \equiv 0 \mod r$, then Corollary~\ref{cor: Semist defect q even} gives $\sd_{\Jpltw / \Kq} = r$ too.	By the discussion above, $\rhiotpmqK$ is principal series if and only if $\ord\st(2 \mod r)$ is odd.
		
		\item Assume that $\frakq = \frakr$ and $\nu_r \leq 2$. Corollary~\ref{cor: Semistability defect} states that the semistability defect of $\Jpmtw / \Kr$ is 
		\begin{equation*}
			\sd_{\Jpmtw / \Kr} = e_{\Krroots / \Kr} \in \lbrace 2, r, 2r \rbrace.
		\end{equation*}
		Therefore, its prime to $r$ part of $\sd_{\Jpmtw / \Kr}$ is $1$ or $2$. Since $\Kr / \Qr$ is totally ramified, we have $|\Ffr\st| = |\Ff_r\st| = r-1$, which is divisible by $1$ and $2$. We conclude that $\rhiotpmrK$ is principal series.
		
		\item Assume that $\frakq \nmid 2r$, $\nu_q \leq 0$ and $\nu_q \not \equiv 0 \mod r$. Then $\Jpmtw$ attains good reduction over $\Kqroots$ (Corollary~\ref{cor: Semistability defect}), which is a tame extension of $\Kq$ with ramification index $r$ or $2r$ (Theorem~\ref{thm: ramif index}). Since $\frakq \nmid 2$, then $|\Ffq\st|$ is even. Therefore, $\rhiotmsqK$ is principal series if and only if $r$ divides $|\Ffq\st|$, which by the discussion above happens if and only if $\ord\st(q \mod r)$ is odd.
	\end{enumerate}
	
	To prove the last statement, assume that $\frakq \neq \frakr$ and that $\rhiotpmqK$ is supercuspidal. The semistability defect of $\Jpmtw$ equals the order of $\rhiotpmqK(I_{\frakq})$ (see Remark~\ref{rmk: Charact semist defect}). Now $\sd_{\Jpmtw / \Kq}$ is either $r$ or $2r$, so the projective image of $\rhiotpmqK$ cannot be $A_4$ nor $S_4$, and the latter \WD-representation is therefore non-exceptional. By Remark~\ref{rmk: Supercuspidal WD reps}, $\rhiotpmqK$ is the induction of a character from a quadratic extension of $\Kq$. Since $\sd_{\Jpmtw / \Kq} = |\rhiotpmqK(I_{\Kq})|$ is always coprime to the residue characteristic (see Corollaries~\ref{cor: Semist defect q even} and \ref{cor: Semistability defect}), we deduce that such quadratic extension is unramified.
\end{proof}

\subsection{Level lowering}\label{sect: Level lowering}

As illustrated in Example~\ref{examp: Conductor Fermat}, when choosing the parameters $s_0, \deltaQ$ as in Proposition~\ref{prop: right choices s0 deltaQ}, the global conductor $\Cond{\rhopmK}$ depends on the putative solution to the considered equation. Such global conductor equals the level of the newform which gives rise to the system $(\rhopm)_{\lambda}^{ }$. Its level depending on the solution is problematic for the elimination step, as one has to perform computations with parameters $a, b, c$ that are not meant to exist. To solve this issue, we apply level lowering results to $\rhopmK$. This way the level of the newform is supported at primes that are independent of the solution, allowing for numerical computations in the elimination step.
\medskip

All the content of sections \S \ref{sect: Reduction types} and \S \ref{sect: Main props reps} so far depended only on the prime number $r$, which we fixed since the beginning. Recall that we aim at solving the infinite families of generalised Fermat equations \eqref{eq: GFE ppr}$_p$ and \eqref{eq: GFE rrp}$_p$. Thus, from now on, we let $p > r$ be a prime number and $\frakp \mid p$ a place of $\Kgl$ dividing $p$. We are mostly interested in the $2$-dimensional  residual representation $\rhopmpbar : G_{\Kgl} \rightarrow \GL_2(\Ffp)$, which corresponds to the particular case $\lambda = \frakp$.

\subsubsection{Absolute irreducibility}\label{sect: Absolute irreducibility}

We begin by proving that, under some $q$-adic condition on $s_0 (s_0 -1)$, the residual representation $\rhopmpbar$ is absolutely irreducible. In Proposition~\ref{prop: rhomsbar irred} we proved that this is the case for $\rhopmbar$, so the same holds for its twist $\overline{\rho}_{\Jpmtw, \, \frakr}$. The particular case $\frakp = \frakr$ being already treated, we focus now on the generic situation $\frakp \neq \frakr$.

\begin{proposition}\label{prop: Irred residue rep q}
	Assume that there is some odd prime $q \neq r$ such that $\nu_q \leq 0$, $\nu_q \not \equiv 0 \mod r$, and $\ord\st(q \mod r)$ is even. Then for any odd place $\frakp$ of $\Kgl$, the residual representation $\rhopmpbar$ is absolutely irreducible.
\end{proposition}

\begin{proof}
	As explained above, we focus on the case $\frakp \neq \frakr$. Since $\rhopmpbar$ is odd and $\Kgl$ is totally real, then $\rhopmpbar$ is absolutely irreducible if and only if it is irreducible. The assumptions on $\nu_q$ imply that $\Jpmtw$ has unipotent and potential good reduction at any $\frakq$ above $q$. The \WD-representation $\rhiotpmqK$ is simply given by 
	\reqnomode
	\begin{equation}\label{eq: rhiotpmqK obtained from tensoring}
		\rhiotpmqK = \left(\restr{\rho_{\Jpmtw, \, \frakp}}{W_{\Kq}} \otimes_{\iota}^{ } \C , 0\right)
	\end{equation}
	(the tensor product is taken along the embedding $\iota : \Kgl_{\frakp} \hookrightarrow \C$). The image of $I_{\Kq}$ through $\rhiotpmqK$ has order $r$ or $2r$, and by \eqref{eq: rhiotpmqK obtained from tensoring}, the same holds for $\restr{\rho_{\Jpmtw, \, \frakp}}{W_{\Kq}}$. Now Proposition~\ref{prop: Local types Jpm} implies that $\rhiotpmqK$ (and thus $\restr{\rho_{\Jpmtw, \, \frakp}}{W_{\Kq}}$) is an irreducible induction of a character from the quadratic unramified extension of $\Kq$. Since $\frakp \nmid 2r$, its reduction $\restr{\rhopmpbar}{W_{\Kq}}$ is also an irreducible induction of a character. The latter being a restriction of $\rhopmpbar$, we obtain the desired result.
\end{proof}

\begin{proposition}\label{prop: Irred residue rep 2}
	Assume that $q = 2$, $\ord\st(2 \mod r)$ is even, and fix some place $\frakp \nmid 2r$. If $\nu_2 \leq -8$ and $\nu_2 \not \equiv -8 \mod r$ then $\rhomspbar$ is absolutely irreducible. If $\nu_2 \leq -4$ and $\nu_2 \not \equiv 0 \mod r$, then $\rhoplpbar$ is absolutely irreducible.
\end{proposition}

\begin{proof}
	Corollary~\ref{cor: Semist defect q even} states that, under the cited conditions, $\rhiotpmqK(I_{\Kq})$ has order $r$. Proposition~\ref{prop: Local types Jpm} implies that $\rhiotpmqK$ is an irreducible induction of a character from the quadratic unramified extension of $\Kq$. Since $\frakp \nmid 2r$, we conclude just as in the proof of Proposition~\ref{prop: Irred residue rep q}.
\end{proof}

\subsubsection{Finiteness}

Level lowering results also require to better understand the local behaviour of $\rhopmpbar$. Namely, we need it to be \textit{finite} at every $\frakq \mid p$, meaning that it arises from a finite flat group scheme over $\Oq$. The next result, which builds upon the content of \cite[\S 7]{BCDF23}, shows that this is indeed the case.

\begin{theorem}\label{thm: Finiteness}
	Let $p$ be a rational prime, and $\frakp \mid p$. Let $q \neq 2, r$ be a rational prime such that $\nu_q > 0$ and $\nu_q \equiv 0 \mod p$. Let $\frakq$ be a finite place of $\Kgl$ above $q$.
	\begin{enumerate}[itemsep=0pt, topsep=5pt]
		\item If $\frakq \nmid p$, then $\rhopmpbar$ is unramified at $\frakq$.
		
		\item If $\frakq \mid p$, then $\restr{\rhopmpbar}{D_{\frakq}}$ is finite at $\frakq$.
	\end{enumerate}
\end{theorem}

\begin{proof}
	We apply \cite[Theorem 7.5]{BCDF23} to the Jacobian $\Jac(\Cpm(s)^{(\deltaQ \deltaK)})$ defined over the function field $\Kq(s)$. Specialising $s = s_0$ yields the well-known $\Jpmtw$. Theorem~\ref{thm: Jpm have RM} combined with Lemma~\ref{lem: quad twist has RM} show that $\Jac(\Cpm(s)^{(\deltaQ \deltaK)})$ has \RM \  by $\Kgl$. The assumption $\nu_q > 0$ implies that, among $\vq(s_0)$ and $\vq(s_0-1)$, exactly one is zero and the other is positive and divisible by $p$. In the notation of \cite[\S 7]{BCDF23}, we let $t_1 \coloneqq s_0$, and then $z_0 \coloneqq 0$ if $\vq(s_0) > 0$, and $z_0 \coloneqq 1$ if $\vq(s_0-1) > 0$, so that $\vq(\pi_0(t_1)) > 0$ and $\vq(\pi_0(t_1)) \equiv 0 \mod p$. Moreover, the three other conditions listed in \cite[page 39]{BCDF23} are satisfied: the first one is due to the fact that $\Cpm(s)$ are \textit{Mumford curves} over $\Kq[[s]]$ and $\Kq[[s-1]]$ (see \cite[Theorem 1.10]{Darmon00}). The second condition follows by considering the discriminants of the models defining $\Cpm(s)^{(\deltaQ \deltaK)}$ (to obtain them it suffices to replace $s_0$ by $s$ in Corollary~\ref{cor: Discriminants Wpm}). The third condition follows from $\Jpmtw$ having toric reduction at $\frakq$, as $\nu_q > 0$ (see \S \ref{sect: Reduction types}). Theorem $7.5$ in \cite{BCDF23} then states that $\rhopmpbar$ is unramified (resp. finite) at $\frakq$ if $\frakq \nmid p$ (resp. $\frakq \mid p$).
\end{proof}

\subsubsection{Applying level lowering}

We are now going to apply level lowering results to the modular compatible system $(\rhopmK)_{\lambda}^{ }$. The newform obtained after this will have a suitable level for our Diophantine purposes. As we will see in Example~\ref{examp: Level lowering}, for the specific choice of $s_0, \deltaQ$ done in Proposition~\ref{prop: right choices s0 deltaQ}, such level will be supported at primes depending on the considered Diophantine equation, but not on its solution. 
\medskip 

Recall the notation $\ntwopm$ and $\nunip$ introduced in Definition~\ref{def: n2 ntor nunip}. For simplicity, we define:

\begin{definition}
	Denote by $\ntornotp$ the square-free product of the primes $\frakq \nmid 2r$ such that $\nu_q > 0$ and $\nu_q \not \equiv 0 \mod p$. 
\end{definition}

We now state our level lowering result. We are going to apply a theorem of Breuil--Diamond \cite{BreuilDiamond14}, that has the advantage of preserving inertial local types. In order to ensure that the considered residual representations are absolutely irreducible, we assume the following hypothesis.

\begin{hypothesis}\label{hyp: place unip red}
	There is at least one prime number $q' \nmid 2r$ such that $\nu_{q'} < 0$, $\nu_{q'} \not \equiv 0 \mod r$, and $\ord\st( q' \mod r)$ is even.
\end{hypothesis}

By Proposition~\ref{prop: Irred residue rep q}, Hypothesis~\ref{hyp: place unip red} implies that there is at least one place of unipotent reduction, and that $\rhopmpbar$ is absolutely irreducible. 	

\begin{theorem}\label{thm: Level lowering}
	Assume that $\frakp \nmid 2r$ and that $\nu_p \geq 0$. When considering $\rhoplpbar$, suppose that $\nu_r > 2$.
	Then, there is a Hilbert newform $g$ over $\Kgl$ of parallel weight $2$, trivial character and level $2^{\ntwopm} \, \frakr^{\cdpmrK} \, \ntornotp^{ } \, \nunip^2$, satisfying the following:
	\begin{enumerate}[itemsep=0pt]		
		\item $\rhopmpbar \simeq \overline{\rho}_{g, \frakP}$ for some place $\frakP \mid p$ in the field of coefficients $K_g$.
		
		\item We have the field inclusion $\Kgl \subset K_g$.
		
		\item If there is some prime $q \nmid 2r$ such that $\nu_q > 0$ and $\nu_q \not \equiv 0 \mod p$, then $g$ does not have complex multiplication.
	\end{enumerate}
\end{theorem}

\begin{proof}
	By Theorems~\ref{thm: Modularity Jms} and \ref{thm: Modularity Jpl}, the representation $\rhoplpbar$ is modular. The assumption on $\frakp$ implies that $\Jmodpmtw$ has good or toric reduction at $\frakp$. Since we assume that Hypothesis~\ref{hyp: place unip red} holds, there is at least one place $\frakq' \mid q'$ at which $\Jmodpmtw$ has unipotent reduction. Moreover, $\rhopmpbar$ is absolutely irreducible by Hypothesis~\ref{hyp: place unip red}, and Lemma $6.9$ in \cite{BCDF23} implies that $\rhopmpbar$ restricted to $G_{\Kgl(\zp)}$ is absolutely irreducible. Since we suppose that $r \geq 5$ and $p > r$, we necessarily have $p \neq 5$. The first claim follows from combining Theorem~\ref{thm: Finiteness} with Theorem $3.2.2$ in \cite{BreuilDiamond14}. 
	\begin{comment}
	If $p = 5$, we claim that $\rhopmpbar(G_{\Kgl(\zp)}) \not \simeq \operatorname{PSL}_2(\Ff_5)$. Indeed, $|\operatorname{PSL}_2(\Ff_5)| = 60$, and Corollary~\ref{cor: Semistability defect} implies that $\rhopmpbar(I_{\frakq '})$ has order $r$ or $2r$, as $\rho_{\Jpmtw, \, \frakp}(I_{\frakq '})$ does not intersect the kernel of the reduction $\mod \frakp$. Since $[\Kgl(\zeta_{5}) : \Kgl] = 4$, the claim follows, and we can also apply \cite[Theorem 3.2.2]{BreuilDiamond14}. We conclude that $\rhopmpbar \simeq \overline{\rho}_{g, \frakP}$, where the level of $g$ is as claimed in the statement.
	\medskip
	\end{comment}
	
	Next, we prove that $\Kgl \subset K_g$. As explained above, $\Jmodpmtw$ has unipotent reduction at $\frakq '$. Consider the complex \WD-representation $\rhiotpmqprK : W_{\frakq'} \rightarrow \GL_2(\C)$ introduced in Definition~\ref{def: WD reps Jpmtw}, which is isomorphic to $\W(\restr{\rho_{g, \frakP}}{D_{\frakq '}}) \otimes \C$ (the tensor product is taken along an embedding $(K_g)_{\frakP} \hookrightarrow \C$). Proposition~\ref{prop: Local types Jpm} implies that $\rhiotpmqprK$ is supercuspidal, and its restriction to $I_{\frakq '}$ can be written as $\delta \oplus \delta^{-1}$, where $\delta : I_{\frakq '} \rightarrow \C\st$ has order $r$ or $2r$.  Since $\Q(\zeta_{2r})^{+} = \Q(\zr)^{+}$, \cite[Proposition 8.4]{BCDF23} gives the desired inclusion. 
	\medskip
	
	Finally, if the condition on the third item is satisfied, then $\rho_{g, \frakP}$ has a Steinberg prime, so $g$ cannot have complex multiplication.
\end{proof}

\begin{example}\label{examp: Level lowering}
	Once again, we specialise Theorem~\ref{thm: Level lowering} to the particular Frey objects obtained by setting $s_0, \deltaQ$ as in Proposition~\ref{prop: right choices s0 deltaQ}. For any $x \in \Z$, define $\radst(x)$ to be the product of all prime ideals of $\OK$ dividing $x$ that are coprime to $2r$. Let $\ntwopm$ and $\cdpmrK$ be the Artin conductors at any even place and at $\frakr$ respectively, as depicted in Tables~\ref{table: Cond Jmstw} and \ref{table: Cond Jpltw}.
	\begin{enumerate}[leftmargin=0em, itemsep=0pt]
		\item[] \boxppr  Assume that $(a, b, c)$ is a primitive non-trivial solution to \eqref{eq: GFE ppr}. Let $s_0, \deltaQ$ be as in Proposition~\ref{prop: right choices s0 deltaQ}. The level of the newform $g$ given by Theorem~\ref{thm: Level lowering} is 
		\begin{equation*}
			\N_g = 2^{\ntwopm} \, \frakr^{\cdpmrK} \radst(AB) \, \radst(C)^2.
		\end{equation*}
		
		\item[] \boxrrp Assume that $(a, b, c)$ is a primitive non-trivial solution to \eqref{eq: GFE rrp}. Let $s_0, \deltaQ$ be as in Proposition~\ref{prop: right choices s0 deltaQ}. The level of the newform $g$ given by Theorem~\ref{thm: Level lowering} is 
		\begin{equation*}
			\N_g = 2^{\ntwoms} \, \frakr^{\cdmsrK} \radst(C) \, \radst(AB)^2.
		\end{equation*}
	\end{enumerate}
\end{example}

% -- % -- % -- % -- % -- % -- % -- % -- % -- % -- % -- % -- % -- % -- % -- % -- % -- % -- % -- % -- % -- % -- % -- % -- % -- % -- % -- % -- % -- % -- % -- % -- % -- % -- % -- % -- % -- % -- % -- 

\section{Solving families of GFEs of signatures $(p, p, r)$ and $(r,r,p)$}\label{sect: Solving instances GFEs}

In this section we explain how to effectively perform the elimination step to solve infinite families of generalised Fermat equations. We begin by summarising the discussion above to solve specific families of equations. We then explain how to discard isomorphisms of Galois representations by comparing traces of Frobenius. To conclude, we present a \Magma \ package that allows to perform the elimination step for several choices of $r, A, B, C$, therefore allowing to solve many families of GFEs.

\subsection{Specialising $s_0$ and $\deltaQ$ for solving specific families of equations}

The whole content of sections \S \ref{sect: Reduction types} and \S \ref{sect: Main props reps} was done for generic values of $s_0$ and $\deltaQ$, although we left a trail of examples after every substantial result. We are now going to fix the values of these parameters as displayed in Table~\ref{table: Main values}. We begin by summarising the strategy of the modular method to solve an instance of a family of GFEs. We focus on the signature $(p, p, r)$, but the discussion for $(r,r,p)$ would follow the same lines.

\subsubsection*{The modular method in practice}

Let $r, A, B, C$ be three fixed integers as in Definition~\ref{def: Parameters GFE}, and let $p > r$ be any prime number. Assume that there exists a primitive non-trivial solution $(a, b, c)$ to the generalised Fermat equation
\begin{equation}\tag{$\Eppr$}
	Ax^p + By^p = Cz^r.
\end{equation}

We proceed in various steps.

\begin{enumerate}[topsep=5pt, itemsep=0pt]
	\item To $(a, b, c)$ correspond two Frey hyperelliptic curves $\Cpmabc / \Q$. By Theorem~\ref{thm: Jpm have RM}, the base-changed Jacobians $\Jpmabc \coloneqq \Jac(\Cpmabc \times \Kgl)$ have \RM \ by $\Kgl$. To them, we associate $2$-dimensional representations $\rho_{\Jpmabc, \, \lambda} : G_{\Kgl} \rightarrow \GL_2(\Klam)$, for every $\lambda \in \Spec(\OK)$.
	
	\item Theorem~\ref{thm: Modularity Jms} implies that the system of Galois representations $(\rho_{\Jmsabc, \, \lambda})_{\lambda}^{ }$ is modular. If $r \mid ab$, or $v_r(AB) > 2$, then Theorem~\ref{thm: Modularity Jpl} also gives the modularity of $(\rho_{\Jplabc, \, \lambda})_{\lambda}^{ }$. Consider a finite place $\frakp$ dividing $p$. By Theorem~\ref{thm: Level lowering}, the residual representation $\overline{\rho}_{\Jpmabc, \, \frakp}$ arises from a Hilbert newform over $\Kgl$ of parallel weight $2$ and level $\N^{\pm} \coloneqq 2^{\ntwopm} \, \frakr^{\cdpmrK} \radst(AB) \, \radst(C)^2$.
	
	\item Let us work with the curve $\Cms(a, b, c)^{(\deltaK)}$ (one could alternatively work with $\Cpl(a, b, c)^{(\deltaK)}$). We numerically compute the space of Hilbert newforms of parallel weight $2$ and level $\N^-$, which we denote by $\Snew{\N^-}$. Following the content of subsection~\ref{sect: Comparing traces}, we prove that, for any $g \in \Snew{\N^-}$, and any place $\frakP \mid p$ in the field of coefficients $K_g$, we have $\overline{\rho}_{\Jpmabc, \, \frakp} \not \simeq \overline{\rho}_{g, \, \frakP}$. This contradicts the previous point, so we conclude that $(a, b, c)$ does not exist.
\end{enumerate}

\subsection{Discarding isomorphisms by comparing traces of Frobenius}\label{sect: Comparing traces}

We now explain how to effectively discard isomorphisms of Galois representations as above, both from a theoretical and algorithmic point of view. For simplicity, we keep working with the curve $\Cms(a, b, c)^{(\deltaK)}$, and we write $J \coloneqq \Jpmabc$, $\N \coloneqq \N^-$, and we keep using the notation from the previous subsection.

\begin{definition}
	Let $g$ be a newform in $\Snew{\N}$. We say that we have \textit{eliminated} the pair $(p, g)$ if we show that, for any place $\frakP \mid p$ in the coefficient field $K_g$, we have 
	\reqnomode
	\begin{equation}\label{eq: Isom to negate}
		\overline{\rho}_{J, \, \frakp} \not \simeq \overline{\rho}_{g, \, \frakP}.
	\end{equation}
\end{definition}

To prove the non-existence of primitive non-trivial solutions to \eqref{eq: GFE ppr}, one has to eliminate all pairs $(p, g)$. If we do so for $p$ greater than a certain bound, we talk about an \textit{asymptotic result}. Let us fix a newform $g \in \Snew{\N}$, and a place $\frakP \mid p$ in its field of coefficients. For any finite place $\frakq$ of $\Kgl$, we consider the traces of Frobenius
\begin{equation*}
	\aqJ \coloneqq \Trace \left( \overline{\rho}_{J, \, \frakp}(\Frob_{\frakq}) \right), \qquad \text{ and } \qquad \aqg \coloneqq \Trace \left( \overline{\rho}_{g, \, \frakP}(\Frob_{\frakq}) \right)
\end{equation*}
assuming that the Néron model of $J$ has good reduction at $\frakq$ in the former case, and that $\frakq \nmid \N$ in the latter case. If the isomorphism $\overline{\rho}_{J, \, \frakp} \simeq \overline{\rho}_{g, \, \frakP}$ holds then $\aqJ$ and $\aqg$ are equal in $\Fpbar$. To test if such an equality holds, we would need to fix embeddings $\Ff_{\frakp} \hookrightarrow \Fpbar$ and $\Ff_{\frakP} \hookrightarrow \Fpbar$, which are not canonical. In order to get rid of this ambiguity, we consider the divisibility relationship:
\reqnomode 
\begin{equation}\label{eq: test for discarding isoms}
	p \mid  \gcd_{\sigma \in \Gal(\Kgl / \Q)} \Nrm{K_g}{\Q} (a_{\sigma(\frakq)} (J) - a_{\sigma(\frakq)} (g)).
\end{equation}
If we had $\overline{\rho}_{J, \, \frakp} \simeq \overline{\rho}_{g, \, \frakP}$, then \eqref{eq: test for discarding isoms} would hold for every $\frakq$. Therefore, as soon as we find some $\frakq$ such that $p$ does not divide the RHS of \eqref{eq: test for discarding isoms}, then the pair $(p, g)$ is eliminated. But the Jacobian $J = \Jpmabc$ is not meant to exist, so how to compute the RHS of \eqref{eq: test for discarding isoms}. We use Remark~\ref{rmk: Curves 2 paramts}, which explains that the curves $\Cpmabc$ admit models that depend only on $a$ and $c$.  Moreover, the curve $\Cpmabc$ is the base change of a curve defined over $\Q$, so the $\aqJ$'s depend only on the congruence classes of $a, c \mod q$. Therefore, it suffices to compute the RHS of \eqref{eq: test for discarding isoms} for $a$ and $c$ ranging through all possible values of $\Ff_{q}$ that do not yield singular curves over $\Ff_q$.
\medskip

Recall that $G_{\Q}$ acts on $\Snew{\N}$ \textit{via} its action on Fourier coefficients. An orbit under this action is called a Hecke constituent. For any $g \in \Snew{\N}$, we denote by $[g]$ its Hecke constituent. Let $q$ be a prime number and $\frakq \subset \OK$ a place above $q$. Assume first that the Néron model of $J$ has good reduction at $\frakq$. Following \cite[\S 9.4]{BCDF23}, we introduce 
\begin{equation}
	\Ngood \coloneqq \prod_{\frakq' \mid q} \ \gcd_{\sigma \in \Gal(\Kgl / \Q)}  \Nrm{K_g}{\Q} (\sigma( a_{\frakq'}(J)) - a_{\sigma(\frakq)}(g)).
\end{equation}
As explained in \textit{loc. cit.}, this quantity depends only on the Hecke constituent $[g]$, and not on the choice of a representative. If $p \nmid \Ngood$, then every pair $(p, g')$ is eliminated, for $g' \in [g]$. 

Suppose now that the Néron model of $J$ has toric reduction at $\frakq$. Assuming that we had $\overline{\rho}_{J, \, \frakp} \not \simeq \overline{\rho}_{g, \, \frakP}$, then $p$ would divide 
\begin{equation*}
	\Mtoric \coloneqq \gcd_{\sigma \in \Gal(\Kgl / \Q)} \Nrm{K_g}{\Q} \left( a_{\sigma(\frakq)}(g)^2 - \left( \Nrm{\Kgl}{\Q}(\frakq) + 1 \right)^2 \right).
\end{equation*}

In practice, to eliminate a pair $(p, g)$ we first range through the different $g$'s, and for each of these, we bound the set of $p$'s that are not yet eliminated. When implementing this discussion, we proceed as follows. First, we initialise $\Snew{\N}$, and compute a basis of  eigenforms: \Magma \ returns a list of the corresponding Hecke constituents. For each of these, we check if the field $\Kgl$ is included in $K_g$ (\cf \ Theorem~\ref{thm: Level lowering}). If not, we can already eliminate the form $g$ for every prime $p$. If $\Kgl \subset K_g$, we initialise the set of $p$'s to eliminate as an infinite set of prime numbers. We run a loops on ascending $q$'s, and on appropriate pairs $(a, c) \in \Ff_q^2$. At every stage, we compute the set of prime divisors of $\Ngood \, \Mtoric$, and we intersect this set with the one obtained for the previous $q$. The remaining $p$'s are those for which the pair $(p, g)$ is not yet eliminated.
\medskip 

This procedure is completely effective, and the author has written a \Magma \ package that performs this elimination process. This package is available at \cite{Git}. To illustrate the explicitness of this elimination process, we specialise the parameter $r = 5$ and solve families of equations of signature $(p, p, 5)$ and $(5, 5, p)$. In this case, the field $\Kgl$ equals $\Q(\sqrt{5})$. We denote by $\frakr_5^{ }$ the unique prime ideal in $\OK$ above $5$. More precisely, we prove the following asymptotic results:

\begin{theorem}\label{thm: Example 1 pp5}
	Let $p > 71$ be any prime number. There are no primitive non-trivial solutions $(a, b, c) \in \Z^3$ to the generalised Fermat equation 
	\begin{equation*}
		7 x^p +  y^p = 3 z^5 \vspace{-0.5em}
	\end{equation*}
	that satisfy $10 \mid ab$.
\end{theorem}

\begin{proof}
	We fix $r = 5, A = 7, B = 1$ and $C = 3$. If $2 \mid ab$, the Artin conductor of $\rho_{\Jmsabc, \, \frakp}$ at any even place is $1$. If $5 \mid ab$ and the twisting parameter $\deltaK$ has $\frakr_{5}^{ }$-adic valuation $1$, the Artin conductor at $\frakr_{5}^{ }$ is also $1$. Assuming that $10 \mid ab$, we get $\N = 2^1 \, 3^2 \, \frakr_{5}^{1} \, 7^1$: for this level, $\Snew{\N}$ has dimension $680$, and there are $101$ Hecke constituents. The output of the elimination process can be found in the file {\fontfamily{lmtt}\selectfont "expp5.txt"} in \cite{Git}.
\end{proof}

\begin{theorem}\label{thm: Example 1 55p}
	Let $p > 41$ be any prime number. There are no primitive non-trivial solutions $(a, b, c) \in \Z^3$ to the generalised Fermat equation 
	\begin{equation*}
		x^5 + 7 y^5 = z^p \vspace{-0.5em}
	\end{equation*}
	that satisfy $10 \mid c$. 
\end{theorem}

\begin{proof}
	This time we set $r = 5 , A = 1, B = 7$ and $C = 1$. If $10 \mid c$, the level given by Theorem~\ref{thm: Level lowering} is $\N = 2^{1} \, \frakr_5^{1} \, 7^2$. For this level, $\Snew{\N}$ has dimension $471$, and there are $71$ Hecke constituents. The output of the elimination process can be found in the file {\fontfamily{lmtt}\selectfont "ex55p.txt"} in \cite{Git}.
\end{proof}

\begin{remark}
	If one drops the assumption $10 \mid c$ in Theorem~\ref{thm: Example 1 55p}, there exist primitive non-trivial solutions to the equation. For $p = 2$, there is $(2, -1, 5)$, and for $p = 3$, there is $\pm (1, 1, 2)$. This is problematic, as to these solutions correspond well-defined Jacobians, which have certain associated Hilbert newforms of the level given by Theorem~\ref{thm: Level lowering}. In this case, in the elimination step, there is some $g$ for which we cannot bound the set of $p$'s such that $(p, g)$ is eliminated. 
	
	Nevertheless, in our particular case, we manage to overcome this difficulty. Indeed, to the mentioned solutions correspond newforms whose level is greater than $\N = 2^{1} \, \frakr_5^{1} \, 7^2$. Therefore, these newforms do not live in the same space as the one where we perform the elimination. 
\end{remark}

\begin{remark}
	Consider the GFE $(\Errp) :  Aa^r + Bb^r = Cc^p$ for some choice of parameters $r, A, B, C$. If we assume that $r \mid c$, then $p$ being odd implies that $Aa^r + Bb^r \equiv 0 \mod r^2$, so Lemma~\ref{lem: r2 divides} gives $A^{r-1} \equiv B^{r-1} \mod r^2$. If the latter congruence is not satisfied, then it is straightforward to see that there are no primitive non-trivial solutions $(a, b, c)$ satisfying $r \mid c$. The choice of coefficients done in Theorem~\ref{thm: Example 1 55p} satisfies such a congruence, so solving the displayed equation is not elementary, and requires the full power of Darmon's program.
\end{remark}

As explained in Remark~\ref{rmk: Curves 2 paramts}, the Frey curves depend only on two parameters among $a, b, c$. On the other hand, the level provided by Theorem~\ref{thm: Level lowering} depends on the radical of the coefficients, but not on the valuations of the latter. For instance, for the signature $(p, p, r)$, the LHS in isomorphism ~\eqref{eq: Isom to negate} is independent of $Bb^p$, whereas the RHS depends only on the primes dividing $B$. 
\medskip 

This allows for a very interesting phenomenon, that we call \textit{``Pay for 1 equation, get n for free."} As a direct consequence of Theorems~\ref{thm: Example 1 pp5} and \ref{thm: Example 1 55p}, we obtain:

\begin{theorem}\label{thm: Example 2 pp5}
	Let $p > 71$ be any prime number. For any $i \in \Iintv{1,4}$ and $j \in \lbrace 3, 4 \rbrace$, there are no primitive non-trivial solutions to the generalised Fermat equation
	\begin{equation*}
		7 x^p + 2^i 5^j y^p = 3 z^5.
	\end{equation*}
\end{theorem}

\begin{proof}
	For this choice of coefficients, the level given by Theorem~\ref{thm: Level lowering} is $\N = 2^1 \, 3^2 \, \frakr_{5}^{1} \, 7^1$, just as in the proof of Theorem~\ref{thm: Example 1 pp5}. Since the Jacobian $\Jms(a, b, c)$ does not depend on $b$, the elimination process from {\fontfamily{lmtt}\selectfont "expp5.txt"} in \cite{Git} shows the non-existence of solutions for any of the considered equations.
\end{proof}

\begin{theorem}\label{thm: Example 2 55p}
	Let $p > 41$ be any prime number. For any $i \in \Iintv{1,4}$ and $j \in \Iintv{2,4}$, there are no primitive non-trivial solutions to the generalised Fermat equation
	\begin{equation*}
		x^5 + 7 y^5 = 2^i 5^j z^p.
	\end{equation*}
\end{theorem}

\begin{proof}
	Just as in the proof of Theorem~\ref{thm: Example 1 55p}, the level here is $\N = 2^{1} \, \frakr_5^{1} \, 7^2$. This time, the Jacobian depends only on $a, b$, so the elimination process from {\fontfamily{lmtt}\selectfont "ex55p.txt"} in \cite{Git} allows to conclude.
\end{proof}

\subsection{Computational aspects of the elimination step}

The elimination process has been implemented in a \Magma \ package, that is available at \cite{Git}. We encourage the reader to use this package, freely adapt it, and report any bugs or suggested improvements! This package allows to perform the elimination step for many choices of parameters $r, A, B, C$, thus solving infinite families of GFEs in great generality. 
\medskip 

There are some computational and theoretical drawbacks about this package. First of all, it does not allow to solve any kind of GFE of signature $(p, p, r)$ or $(r, r, p)$. There are 2 limitations for this:
\begin{itemize}
	\item First, our description of the Artin conductor (Theorems~\ref{thm: Conductor rhomsK} and \ref{thm: Conductor rhoplK}) at even places is not complete, so for certain choices of coefficients we cannot compute the level of the newforms to eliminate. 
	
	\item Second, we do not prove absolute irreducibility of the residual representations for all choices of parameters $r, A, B, C$ (Propositions \ref{prop: Irred residue rep q} and \ref{prop: Irred residue rep 2}). Indeed, it requires some local conditions, that not all choices of coefficients satisfy. Nevertheless, the user can still employ our \Magma \ package by setting a certain parameter to be true. The user may provide a theoretical argument different than the one given in \S \ref{sect: Absolute irreducibility}, or have a Diophantine result conditional to absolute irreducibility.
\end{itemize}

\begin{remark}
	Assume the user runs our package for a choice of parameters whose corresponding equations has a primitive non-trivial solution. To this solution corresponds a well-defined Jacobian, which is associated with an existing newform. Therefore, there exists a newform $g$ for which we do not manage to bound the set of primes $p$'s such that $(p, g)$ has not been eliminated yet. Therefore, the elimination process does not end for this newform $g$, and the script should run indefinitely, never exiting the loop. 
	
	However, we have considered this scenario, and we introduce some timers forcing the process to stop after a finite amount of time (or a finite number of steps). The user can manually set these parameters, modifying the maximal amount of time spent trying to eliminate each newform, or the maximal amount of $q$'s to be used. This way, the script always terminates.
\end{remark}

Computing spaces of Hilbert newforms is computationally expensive, and the complexity grows very quickly with the level. Many of the newforms are often eliminated because their field of coefficients does not contain the field $\Kgl$. For the other ones, the execution time for bounding the set of primes grows with the degree of the field of coefficients of the form. For some given forms, we observe that certain primes are difficult to eliminate. To eliminate such pairs $(p, g)$, one could use refined versions of the elimination step as in \cite[\S 9.11]{BCDF23}. This is work in in progress. 
\medskip

All the output files and running times can be found in \cite{Git}. We hope this \Magma \ package will be useful for proving other Diophantine results regarding families of generalised Fermat equations!

% -- % -- % -- % -- % -- % -- % -- % -- % -- % -- % -- % -- % -- % -- % -- % -- % -- % -- % -- % -- % -- % -- % -- % -- % -- % -- % -- % -- % -- % -- % -- % -- % -- % -- % -- % -- % -- % -- % -- 

\bibliographystyle{alpha}
{\footnotesize \bibliography{main}}

\normalsize\vfill
\noindent\rule{7cm}{0.5pt}

\smallskip
\noindent
\textsc{Martin Azon} (\textit{martin.azon@uca.fr}) --
{\sc Laboratoire de Mathématiques B. Pascal, Université Clermont Auvergne,} 
Campus des Cézeaux. 3 place Vasarely, TSA 60026 CS 60026, 63178 Aubière Cedex (France). 
	
\end{document}